\documentclass[a4paper,10pt]{amsart}

\usepackage{systeme}
\usepackage{graphicx}
\usepackage[]{units}
\usepackage{color}
\usepackage[dvipsnames]{xcolor}
\usepackage{mathrsfs}
\usepackage{bm}
\usepackage[ansinew]{inputenc}
\usepackage[breaklinks,hidelinks]{hyperref}
\usepackage{float}
\usepackage{amsmath}
\usepackage{isomath}
\usepackage{amsthm}
\usepackage{amssymb}
\usepackage{amssymb}

\usepackage{mathtools}

\usepackage[foot]{amsaddr}

\theoremstyle{definition}

\theoremstyle{remark}

\usepackage[margin=2.7cm]{geometry}

\usepackage{systeme}
\usepackage{graphicx}
\usepackage[]{units}
\usepackage{color}
\usepackage{mathrsfs}
\usepackage{bm}
\usepackage[ansinew]{inputenc}
\usepackage{float}
\usepackage{amsmath}
\usepackage{verbatim}

\usepackage[hidelinks]{hyperref}
\numberwithin{equation}{section}

\newcommand{\ldb}{\mathopen{\lbrack\!\lbrack}}
\newcommand{\rdb}{\mathclose{\rbrack\!\rbrack}}

\newcommand{\blue}{\color{black}}

\begin{document}
\title[Projection-based ROM\lowercase{s} for a C\lowercase{ut}FEM method in parametrized domains]{Projection-based reduced order models for a cut finite element method in parametrized domains}

\author{Efthymios N. Karatzas\textsuperscript{1,2,*}}
\address{\textsuperscript{1}SISSA, Mathematics Area, mathLab, Via Bonomea 265, Trieste, 34136, Italy.}
\address{\textsuperscript{2}
Department of Mathematics, School of Applied Mathematical and Physical Sciences, 
National Technical University of Athens, Zografou 15780, Greece.}
\thanks{\textsuperscript{*}Corresponding author}
\email{efthymios.karatzas@sissa.it}


\author{Francesco Ballarin\textsuperscript{1}}
\email{francesco.ballarin@sissa.it}

\author{Gianluigi Rozza\textsuperscript{1}}
\email{gianluigi.rozza@sissa.it}





\date{\today}

\begin{abstract}
This work presents a reduced order modeling technique built on a high fidelity embedded mesh finite element method. Such methods, and in particular the CutFEM method, are attractive in the generation of projection-based reduced order models thanks to their capabilities to seamlessly handle large deformations of parametrized domains {\blue{and in general to handle topological changes}}. The combination of embedded methods and reduced order models allows us to obtain fast evaluation of parametrized problems, avoiding remeshing as well as the reference domain formulation, often used in the reduced order modeling for boundary fitted finite element formulations. The resulting novel methodology is presented on linear elliptic and Stokes problems, together with several test cases to assess its capability. The role of a proper extension and transport of embedded solutions to a common background is analyzed in detail.
\end{abstract}
\maketitle
\section{Introduction and Motivation}\label{sec:intro}
A wide variety of numerical methods and computational libraries for the solution of problems governed by partial differential equations is nowadays available.
However, there are still many cases in which either the solution of the governing equations or the solution of associated inverse problems becomes impractical or unfeasible using standard discretization techniques, such as the Finite Element Method (FEM).
For instance, complicated topology of the problem or a complex geometry may pose a challenge in the discretization of complex phenomena, and ultimately affect the quality of the resulting simulation. Moreover, repeated queries to the underlying solver in the context of an iterative solution of inverse problems may result in unbearably large computational times.
Such situations occur, for example, when a large number of different configurations are in need of being tested, such as in uncertainty quantification, optimal control and shape optimization.
The  overall objective of this manuscript is to investigate how the recently introduced unfitted mesh finite element methods
may be used combined with reduced order modeling techniques for parametrized partial differential problems. 
Indeed, unfitted methods are very useful in cases characterized by complex geometrical configurations; however, as their Finite Element counterparts, they usually require large computational efforts. A combination with reduced order methods, able to widely decrease the overall computational time, would result in a very compelling methodology to be possibly applied in several different fields.

Classical embedded/immersed methods provide simple, efficient, and robust numerical schemes for solving PDE in general domains \cite{MiIa05,MoDoBe99,RoHeCo15}. 
Since these early works, several improvements have been made, for instance for what concerns the rate of spatial accuracy near embedded boundaries.
Recent improvements go under the names of Ghost-Cell finite difference methods, Cut-Cell finite volume approach, Immersed Interface, Ghost Fluid, Volume Penalty methods, for which we refer to the review paper \cite{MiIa05} and references within. 
In particular, for what concerns incompressible flows in arbitrary smooth domains, the Immersed Boundary Smooth Extension method has shown high-order convergence for the incompressible Navier-Stokes equations \cite{SteGuTho17}.

More in detail, extended mesh finite element methods using cut elements are examined in \cite{Bu_et_all_14,BuHa14} for stationary Stokes flow systems, as well as for Navier-Stokes. An analysis for high Reynolds numbers, independent of the local Reynolds, has been carried out in \cite{BuErnFe17,BuFe07,BuFeHa06}. XFEM approaches in 2D and 3D Navier-Stokes are reported in \cite{SchWa14}. Higher Reynolds number aerodynamics problems in unbounded domains, thin vortex ring and Lattice Green function Immersed Boundary methods are studied in \cite{LiCo16,MeLiYuCoJa17}. Furthermore, embedded and immersed methods have been used in solving fluid structure interaction problems, see e.g. \cite{Bu07,BuFe07_2,BuFe09,BuFe14,CoFo15,GeWa08,KaBhaGriDo16,RoHeCo15,TaCo07,TaCo08,Wa11,WaLeMaMcGaFa14}.
More recently a new embedded boundary method called Shifted Boundary method was introduced in \cite{MaSco17_2}, with good approximation properties and experimental results in several fields, including Navier-Stokes systems and fluid structure interaction applications \cite{MaSco17_2,MaSco17_3,SoMaScoRi17}.
The main idea of this  method is that boundary conditions are imposed on the boundary of a surrogate domain whose geometry is properly chosen avoiding cut cells.

First attempts to apply Reduced Order Methods (ROMs) in the context of viscous flows and Stokes / Navier-Stokes systems can be found in \cite{ITO1998403,peterson1989}. 
ROMs based on FEM full order approximations have been used to treat several problems based on linear elliptic equations \cite{Rozza2008229}, linear parabolic equations \cite{grepl2005} and even non-linear problems \cite{Grepl2007,Veroy2003}.
Stability of the resulting reduced order systems is often an issue, see for more details \cite{ballarin2015supremizer,Caiazzo2014598,Gerner2012,RoHuMa13,RoVe07,HiAlStaBaRo17}, and for transient systems see for instance \cite{Akhtar2009,Bergmann2009516,taddei2017,Iollo2000,Sirisup2005218}.

The main novelty of this work is the combination of the cut element finite element method (CutFEM), which we will consider as high fidelity method, to projection-based model reduction techniques. In particular, this allows to overcome a reference domain formulation customarily employed in ROMs for parametrized geometries built on conforming discretizations, see e.g. \cite{BallarinRozza2016,BeOhPaRoUr17,Laura,LeRa18,Manzoni2017,Rozza2009,RoHuMa13,Rozza2008229,RoVe07} and references therein.
Even though such reference domain formulation avoids remeshing when updating the parametric domain, the choice of the transformation map is usually problem-dependent \cite{BallarinJCP,BallarinRozza2016,LeRa18}, suitable parameter space dimensionality reduction techniques \cite{BallarinDAmarioPerottoRozza2017,BallarinManzoniRozzaSalsa2013,Marco2,Marco1} must be employed to identify the most relevant shape parameters and preserve good quality meshes, and often limited only to small parametric deformations. In contrast, the combination of ROMs with a CutFEM approach results in a novel methodology capable of handling large parametric deformations as well.

To the best of the authors' knowledge this combination is not particularly investigated. We mention here the approach proposed in \cite{BaFa2014} for classical embedded methods and model reduction, in which two regions are separated by an evolving in time interface, where a snapshot compression problem as a weighted low-rank approximation is formulated. We also mention our on-going work on combination of the Shifted Boundary method and model reduction techniques for small parametric variations \cite{KaratzasStabileAtallahScovazziRozza2018,KaratzasStabileNouveauScovazziRozza2018,KaratzasStabileNouveauScovazziRozzaNS2019,KaratzasRozzaCH2019}.
The present work touches several key points of growing interest in the model reduction community, especially since it shares challenges to be tackled for the efficient reduction of advection dominated problems. In particular, a snapshots preprocessing is advocated by \cite{BeIoRi2018,CaMaSta19,IolloLombardi2014,Naira2017,Reiss2015,Welper2017} in order to improve the efficiency of the resulting reduced model for advection dominated problems by providing a better (smaller) representation of the reduced basis space. Such snapshots preprocessing procedure, based in our case on extension and transportation, will be pivotal in this work in order to deal with large deformations.

The work is organized as it follows: in Section \ref{sec:HF} the abstract formulation of the problem is introduced, as well as the embedded method used for the high fidelity problem of Darcy flow pressure model and of the steady Stokes equations. The reduced order methodology is discussed in details in Section \ref{sec:ROM}, whilst in Section \ref{sec:num_exp} the proposed ROM technique is tested on several numerical benchmarks. 
Finally in Section \ref{sec:conclusions} conclusions and perspectives are drawn, highlighting the directives for future improvements and developments.
\section{High fidelity CutFEM approximation} \label{sec:HF}
We introduce in this section the abstract parametrized model problem, which has the form:
find $u(\mu)\in V(\mu)$ such that
\begin{equation}\label{eqn:abstract}
a_\gamma(u(\mu), v;\mu) = \ell_\gamma(v;\mu) \text{ in } V(\mu),
\end{equation}
where $a_\gamma(\cdot, \cdot; \mu)$ is the weak form of an operator defined on a domain $\mathcal D (\mu)\subset \mathbb R ^d$,  for $d= 2,3$, while $\ell _\gamma(\cdot,\mu)$ is the right hand side of the system of equations related to the forcing term. 
The forms $a_\gamma(\cdot, \cdot; \mu)$, $\ell _\gamma(\cdot,\mu)$, the domain $\mathcal D (\mu)$, as well as the Sobolev space $V(\mu)$, depend on the parameter $\mu \in \mathcal K$, being $\mathcal K \subset \mathbb R^{K}$ the set of possible outcomes, which we assume to be a compact set in $\mathbb R^{K}$, $K \in \mathbb{N}$. Since CutFEM discretizations usually require suitable penalty and stabilization procedures, for the sake of exposition we will denote by $a_\gamma(\cdot, \cdot; \mu)$ and $\ell _\gamma(\cdot,\mu)$ the forms with penalty and stabilization, where $\gamma$ represents the penalty and stabilization coefficients, see e.g. \cite{Bu_et_all_14} and references therein.
 
We start with a sketch description of the continuous strong form and the weak formulation used for the problems under consideration. The CutFEM formulation \eqref{eqn:abstract} will be used for the high fidelity simulation employed during the training of the reduced order model, which will be introduced in Section \ref{sec:ROM}.

\subsection{The Darcy flow pressure model}\label{sec:HFDarcy}

For any $\mu\in \mathcal K
$,  let $\mathcal D(\mu)\subset {\mathbb R}^d$ be  a bounded domain 
depending on $\mu$,  
 with boundary $\Gamma(\mu)$. We consider the following model problem in $\mathcal D(\mu)$: for any $\mu \in \mathcal K$, find $u(\mu):{\overline{\mathcal D}(\mu)} \to \mathbb R$ such that
\begin{equation}\label{P:Poisson}
\begin{cases}
 -\Delta  u(\mu)=g(\mu), & \text{in $\mathcal D(\mu)$},
 \\
  u(\mu)=g_D(\mu), & \text{on $\Gamma_D(\mu)$},
  \\
  {\bf{n}_\Gamma}\cdot \nabla u(\mu)=g_N(\mu), &\text{on $\Gamma_N(\mu)$},
\end{cases}
\end{equation}
 and $\Gamma_D(\mu)$, $\Gamma_N(\mu)$ are non-overlapping parts of $\Gamma(\mu) = \Gamma_D(\mu) \cup \Gamma_N(\mu)$ where Dirichlet and Neumann boundary conditions are applied,  and $g(\mu)$,   $g_N(\mu)$, $g_D(\mu)$ are given functions in $\mathcal D(\mu)$ and on the boundaries $\Gamma_N(\mu)$, $\Gamma_D(\mu)$  respectively, \cite{BraGhaLyWu09}.

We start with the classical weak formulation in the spatial domain: for any $\mu\in \mathcal K$, find $u(\mu)\in V_{g_D}(\mu)=\left\{w\in H^1\left(\mathcal D(\mu)\right) \text{ with } w|_{\Gamma_D(\mu)}=g_D(\mu)\right\}$ such that
\begin{eqnarray}
  \left( \nabla u(\mu) ,\nabla \upsilon(\mu)\right) = \left(g(\mu) , \upsilon(\mu) \right)
  + \left(g_N(\mu) , \upsilon(\mu) \right)_{\Gamma_N(\mu)}, \quad\forall v(\mu)\in V_0(\mu).\label{P:ContWeakForm}
\end{eqnarray}

\begin{figure}\centering
(a)\\
  \includegraphics[width=.79\textwidth]{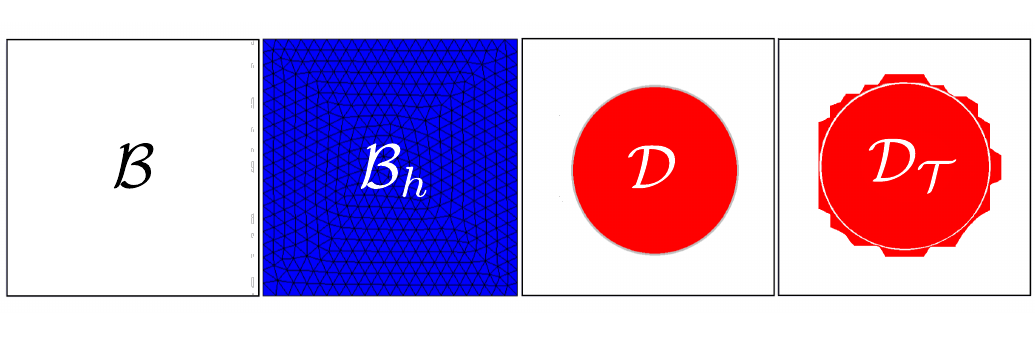}
\newline
(b)\\
  \includegraphics[width=0.67\textwidth]{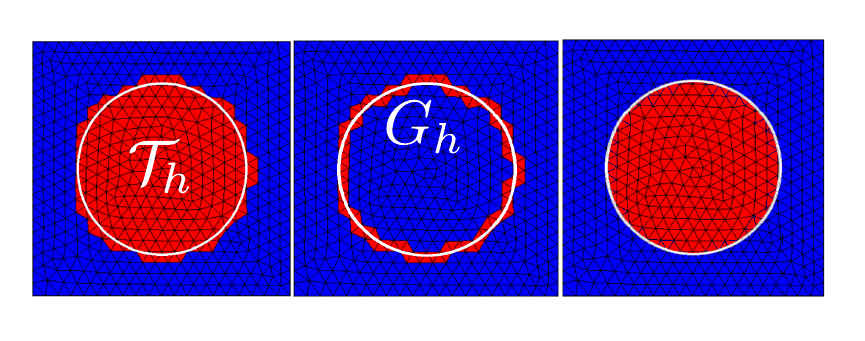}
    \includegraphics[width=0.3\textwidth]{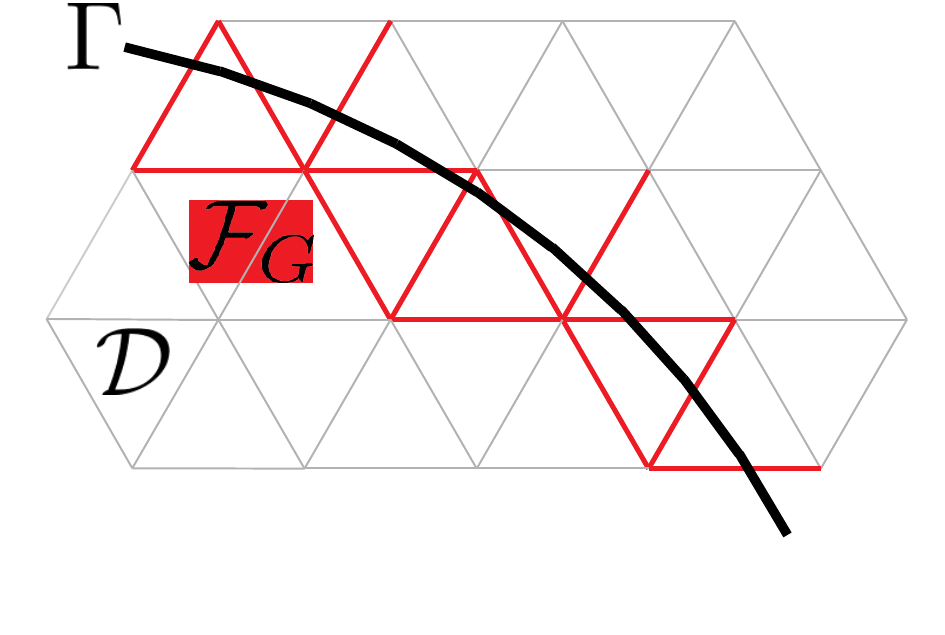}
  \caption{Geometric quantities employed in the definition of the CutFEM discretization: (a) the background geometry, the background mesh, the geometry of a disk, the surrogate geometry, and (b) the unfitted/extended geometry mesh, the cut elements, the cut mesh {\blue{and the set of element faces $\mathcal{F}_G$ associated with ${G}_h$}}.}
  \label{SurrogateMesh}
\end{figure}

We denote by $\mathcal{B}$ the background domain, and by $\mathcal{B}_h$ its corresponding mesh. The continuous boundary value problem is next formulated on a domain ${\mathcal D}_{\mathcal{T}}(\mu)$ that contains ${\mathcal D}(\mu) \subset {\mathcal D}_{\mathcal{T}}(\mu)$, while its mesh $\mathcal{T}_h(\mu)$ is not fitted to the domain boundary and ${\mathcal D}_{\mathcal{T}}(\mu) \subset \mathcal{B}$ and $\mathcal{T}_h(\mu) \subset \mathcal{B}_h$ for all $\mu \in \mathcal{K}$. Let also ${G}_h(\mu):=\{K\in \mathcal{T}_h(\mu): K\cap\Gamma(\mu)\neq\emptyset\}$ be the set of elements that are intersected by the interface. We remark that $\mathcal{T}_h(\mu)$, ${G}_h(\mu)$ and ${\mathcal D}_{\mathcal{T}}(\mu)$ depend on $\mu$ through $\mathcal{D}(\mu)$ (or its boundary), while the background domain $\mathcal{B}$ and its mesh $\mathcal{B}_h$ \emph{do not} depend on $\mu$.

Furthermore, the set of element faces $\mathcal{F}_G(\mu)$ associated with ${G}_h(\mu)$, is defined as follows: for each face $F\in \mathcal{F}_G(\mu)$, there exist two simplices $K \neq K'$ such that $F=K\cap K'$ and at least one of the two is a member of ${G}_h(\mu)$. Note that boundary faces of $\mathcal{T}_h(\mu)$ are excluded from $\mathcal{F}_G(\mu)$.
On a face $F \in \mathcal{F}_G(\mu)$, $F=K\cap K'$, the jump of the gradient of $v \in C^0(\overline{\mathcal{D}}_{\mathcal{T}})$ is defined by $\ldb{{\bf n}_F}\cdot\nabla  v\rdb={\bf n}_F\cdot\nabla v |_K -{\bf n}_F\cdot \nabla v|_{K'}$, where ${\bf n}_F$ denotes the outward pointing unit normal vector to $F$. Finally, let $h_k$ be the diameter of $K$, $h=\max _{K\in \mathcal{T}_h(\mu)} h_k$.

The CutFEM discretization is as follows. We seek a discrete solution $u_h(\mu)$ in the finite element space
\begin{equation}
V_h(\mu)=\left \{  \upsilon \in C^0({\overline {\mathcal D}}_{\mathcal T}(\mu))\,:\,\upsilon |_K  \in P^1(K), \, \forall K\in \mathcal{T}_h(\mu)    \right \},
\label{eqn:V_poisson}
\end{equation}
such that
\begin{eqnarray}\label{CutFem_Discrete_Form_mu}
  &&\left( \nabla u_h(\mu) ,\nabla \upsilon_h(\mu)\right)_{{\mathcal D}(\mu)}-\left({\bf n}_\Gamma\cdot \nabla u_h(\mu),\upsilon_h(\mu)\right)_{\Gamma_D(\mu)}
- \left(u_h(\mu),{\bf n}_\Gamma \cdot \nabla\upsilon_h(\mu)\right)_{\Gamma_D(\mu)}
  \nonumber \\
  && \quad+\left(\gamma_Dh^{-1}u_h(\mu), \upsilon_h(\mu)\right )_{\Gamma_D(\mu)}
   +\left(\gamma_N h {\bf n}_\Gamma\cdot \nabla u_h(\mu),{\bf n}_\Gamma\cdot \nabla \upsilon_h(\mu) \right)_{\Gamma _N(\mu)}
 + { j(u_h(\mu), \upsilon_h(\mu);\mu)  }
 \nonumber \\
  &&\qquad= \left(g
  {(\mu)}, \upsilon_h(\mu)\right )_{{\mathcal D}{(\mu)}} +\left(g_D(\mu),\gamma_D h^{-1} \upsilon_h(\mu) -  {\bf n}_\Gamma\cdot \nabla \upsilon_h(\mu)\right)_{\Gamma _D{(\mu)}} +\big(g_N(\mu), \upsilon_h(\mu)
\nonumber \\
&&\qquad\qquad\qquad\qquad\qquad\qquad\qquad\qquad\qquad\qquad\qquad\qquad\qquad\qquad +\gamma_N h {\bf n}_\Gamma\cdot \nabla \upsilon_h(\mu)\big)_{\Gamma _N{(\mu)}} ,
\end{eqnarray}
 \begin{eqnarray*}
 \text{ where}\qquad j(u_h(\mu), \upsilon_h(\mu);\mu)=\sum_{F\in \mathcal{F}_G(\mu)}\left(\gamma_1 h\ldb{{\bf n}_F}\cdot\nabla  u_h(\mu)\rdb,\ldb{{{\bf n}_F}\cdot\nabla} \upsilon_h(\mu)\rdb\right)_F,
  \end{eqnarray*}
  for all test functions $\upsilon_h(\mu)\in V_h(\mu)$ and  for $\gamma_D \blue{> 0}$, $\gamma_N \blue{\geq 0}$, and $\gamma_1 \blue{> 0}$  penalty parameters, see for instance \cite{BuHa11}.
 The stabilization term $j(u_h(\mu), \upsilon_h(\mu);\mu)$, which depends on $\gamma_1$, extends  the coercivity from the physical domain ${\mathcal D}(\mu)$ to the extended mesh domain ${\mathcal D}_{\mathcal{T}}(\mu)$, while the penalty terms involving coefficients $\gamma_D$ and $\gamma_N$ account for a Nitsche weak imposition of boundary conditions.
Assuming that for the Nitsche terms we use the notation
 \begin{eqnarray*} 
b_\gamma (u_h(\mu),\upsilon_h(\mu);\mu):=- \left(u_h(\mu),{\bf n}_\Gamma \cdot \nabla\upsilon_h(\mu)\right)_{\Gamma_D(\mu)}
+\left(\gamma_Dh^{-1}u_h(\mu), \upsilon_h(\mu)\right )_{\Gamma_D(\mu)}
 \\+\left(\gamma_N h {\bf n}_\Gamma\cdot \nabla u_h(\mu),{\bf n}_\Gamma\cdot \nabla \upsilon_h(\mu) \right)_{\Gamma _N(\mu)}
,
  \end{eqnarray*}
the variational form \eqref{CutFem_Discrete_Form_mu}, using the classical linear and bilinear forms, can be equivalently expressed as: find $u_h(\mu) \in V_h(\mu)$ such that
\begin{eqnarray}
a_\gamma({u_h(\mu)},{ \upsilon _h(\mu)}; \mu) = \ell _\gamma(\upsilon _h(\mu); \mu), \quad\forall \upsilon _h(\mu) \in V_h(\mu),
\label{eqn:hf_lapl}
\end{eqnarray}
where
\begin{eqnarray}
a_\gamma({u_h(\mu)},{ \upsilon _h(\mu)}; \mu) := a({u_h(\mu)},{ \upsilon _h(\mu)}; \mu) + b_\gamma ({u_h(\mu)},{ \upsilon _h(\mu)}; \mu) + j(u_h(\mu), \upsilon_h(\mu);\mu),
\end{eqnarray}
\begin{eqnarray}
\ell _\gamma(\upsilon _h(\mu);\mu)&=& (g(\mu), \upsilon _h(\mu))    +\left({g_N(\mu), \upsilon_h(\mu) + b_\gamma (g_D(\mu),\upsilon _h(\mu))}
 \right)_{{\Gamma _N{(\mu)}}} ,\nonumber
\end{eqnarray}
and
\begin{eqnarray}
a({u_h(\mu)},{ \upsilon_h(\mu)}; \mu)={\left( \nabla u_h(\mu) ,\nabla \upsilon_h(\mu)\right)_{{\mathcal D}(\mu)}}
{-\left({\bf n}_\Gamma\cdot \nabla u_h(\mu),\upsilon_h(\mu)\right)_{\Gamma_D(\mu)}}. 
\end{eqnarray}
We refer to \cite{BuHa11} for more details concerning the CutFEM discretization of elliptic problems.

\subsection{Steady Stokes problem}\label{sec:HFStokes}
The strong form of the stationary Stokes flow system of equations with Dirichlet and Neumann
 boundary conditions, geometrically parametrized by $\mu$,  is given by: for any $\mu \in \mathcal{K}$, find velocity $u(\mu): {\overline{\mathcal D}(\mu)} \to \mathbb R^d$ and pressure $p(\mu): {\overline{\mathcal D}(\mu)} \to \mathbb R$ such that
\begin{equation*}
\begin{cases}
-\nabla\cdot (2\nu{ \varepsilon (u(\mu))} - p(\mu){ I}) = g(\mu), &\text{ in } {\mathcal{D}(\mu)}, 
\\
\nabla \cdot { u(\mu)} = 0, &\text{ in }{\mathcal{D}(\mu)}, 
\\ u(\mu) = g_D(\mu), &\text{ on } \Gamma_D(\mu),
\\
(2\nu {\varepsilon (u(\mu))} - p(\mu) I) \cdot {\bf n}_\Gamma = 0, &\text{ on } \Gamma _N(\mu),
\end{cases}
\end{equation*}
where ${ \varepsilon ({ u})} = 1/2( \nabla   { u}+ \nabla { u}^T )$ is the velocity strain tensor (i.e., the symmetric gradient of the velocity), $\nu$ is the viscosity, $\nabla \cdot u$ denotes the divergence of $u$, $g$ a body force, and $ g_{D}$ the value of the velocity on the Dirichlet boundary.
The first equation represents the conservation of the linear momentum of the fluid, while the second equation is the incompressibility condition and enforces the mass conversation, see e.g. \cite{MaSco17_2} and references therein.

As in the Darcy case, the high fidelity CutFEM formulation is based on three ingredients: (i) discrete FE spaces, (ii) Nitsche weak imposition of Dirichlet boundary conditions, and (iii) suitable stabilization terms.
For what concerns the first topic, 
using similar notation as in the elliptic case for geometrical quantities,
we introduce the discrete spaces ${V}_h(\mu)$ and $Q_h(\mu)$, for the velocity and the pressure, respectively, as piecewise linear spaces as in \cite{BuHa11}:
\begin{align*}
V_h(\mu) &= \left\{ \upsilon_h \in (C^0 (\overline{\mathcal{D}}_{\mathcal T}(\mu)))^{d} : \upsilon_h|_K \in ({P}^1 (K))^{d}, \forall K \in \mathcal{T}_h(\mu)\right\},\\
Q_h(\mu) &= \left\{ w_h \in C^0 (\overline{\mathcal{D}}_{\mathcal T}(\mu)) : w_h|_K \in  {P}^1 (K), \forall K \in \mathcal{T}_h(\mu)\right \},
\end{align*}
Alternative choices are possible, such as, for example, discontinuous Galerkin spaces, see e.g. \cite{AntoGeorg_etall16,GeSu05, Cangiani_et_all}. For shortness in the remainder of this Section we will drop the suffix $h$ when referring to discrete functions $\upsilon_h$, as well as sometimes omit the parameter $\mu$ in our notation.

The Nitsche discrete weak formulation with focus on the Dirichlet boundary conditions of the embedded geometry, takes the form 
\begin{eqnarray*}
                        \nu(\nabla u , \nabla \psi) - (p, \nabla \cdot \psi)  +( \nabla \cdot u,  \xi)  - (\nu n_\Gamma \cdot  u - pn, \psi)_{\Gamma _D}
-(u, \nu n_\Gamma \cdot  \psi + \xi n, \psi)_{\Gamma _D} 
\\
        + \sum _{K\in G_h}\int_{\Gamma_K}\nu {\gamma_D }{h^{-1}_K}u \psi
 =( g, \psi) - ( g_D, \nu n_\Gamma\cdot \nabla \psi +\xi n )_{\Gamma _D} 
        + \sum _{K\in G_h}\int_{\Gamma_K}\nu {\gamma_D }{h^{-1}_K}g_D \psi,
\end{eqnarray*}
being $\gamma_D$ a penalty coefficient.
For the sake of notation, let us introduce
\begin{eqnarray}
b_\gamma (\bm u,\bm \phi) &=& -( u, \nu {\bf n}_\Gamma\cdot \nabla \psi +\xi {\bf n}_\Gamma)_{\Gamma _D} + \sum _{K\in G_h}\int_{\Gamma_K}\nu {\gamma_D }{h^{-1}_K}u \psi,
\nonumber\\
a({\bm u},{ \bm \phi})&=& ({\blue{\nu}} \nabla u, \nabla \psi) - (p, \nabla \cdot \psi)  +( \nabla \cdot u,  \xi)  - (\nu {\bf n}_\Gamma \cdot\nabla  u - p{\bf n}_\Gamma, \psi)_{\Gamma _D} ,
\nonumber\\
\ell _\gamma(\bm \phi)&=& (g, \psi)  + b_\gamma (g_D,\bm \phi),
\nonumber
\end{eqnarray}
where ${\bm u}=(u, p)$ and ${\bm \phi}=(\psi, \xi)$ denote the pair of velocity and pressure functions.

Finally, introducing suitable stabilization terms, an extended mesh weak form can be expressed using a symmetric form for the velocity and pressure, as follows: find ${\bm u(\mu)}=(u(\mu), p(\mu)) \in \bm{V}_h(\mu) := {V}_h(\mu) \times Q_h(\mu)$  such that
\begin{equation}\label{eqn:weak Stokes}
a_\gamma (\bm u(\mu),\bm \phi(\mu);\mu) =  \ell _\gamma(\bm \phi(\mu);\mu), \quad\forall {\bm \phi(\mu) = (\psi(\mu), \xi(\mu))} \in \bm{V}_h(\mu).
\end{equation}
where
\begin{equation*}
a_\gamma (\bm u,\bm \phi) := a({\bm u},{\bm\phi}) + b_\gamma (\bm u,\bm \phi) + j_u(u, \psi; \mu) - j_p(p, \xi; \mu),
\end{equation*}
The terms $j_u$ and $j_p$ account for the stabilization due to the equal order FE spaces, and are expressed as follows:
\begin{align*}
  j_u(u,\psi;\mu)&=\sum_{i=1}^2\sum_{F\in \mathcal{F}_G(\mu)} \left(\gamma_{1,u} h_K \ldb \nabla u_i(\mu) {\bf n}_F\rdb, \ldb \nabla \psi_i(\mu) {\bf n}_F \rdb\right)_F,\\
  j_p(p,\xi;\mu)&=\sum_{F\in \mathcal{F}_G(\mu)} \left(\gamma_{1,p} h^3 _K \ldb \nabla  p(\mu)\cdot{\bf n}_F\rdb,\ldb \nabla \xi(\mu)\cdot{\bf n}_F\rdb\right)_F,
\end{align*}
for all test functions $ \bm\phi = (\psi, \xi) \in  { V}_h \times Q_h$, and for $\gamma_{1,u}$, $\gamma_{1,p}$  positive penalty parameters.
We refer to \cite{BuHa11} for further details on CutFEM discretization for Stokes problems.
\section{POD--Galerkin reduced order method}\label{sec:ROM}
In this Section we obtain a projection-based reduced order model (ROM) built on the high fidelity discretization introduced in Section \ref{sec:HF}. Following the ROM procedure of \cite{HeRoSta16}, a two stage procedure will be used, the offline and the online.
During the offline stage, one examines the solution manifold to construct a reduced basis that well approximates the manifold through a small number of basis functions. 
This may involve the solution of a large number of high fidelity problems, which may be expensive to query. 
In contrast, during the online stage, a Galerkin projection of the problem \eqref{eqn:abstract} onto the space spanned
by the reduced basis is required. During this stage each solution for a new value of $\mu$ entails
substantially reduced cost. 

\subsection{The Darcy flow pressure model}\label{subsec_POD_theory}
Let us start from the elliptic case introduced in Section \ref{sec:HFDarcy}. 
In order to generate the reduced basis space we will employ a compression by means of a Proper Orthogonal Decomposition (POD) \cite{HeRoSta16}, even though several other options (e.g. based on greedy procedures, such as the certified Reduced Basis method (RB) \cite{ChinestaEnc2017,HeRoSta16,Kalashnikova_ROMcomprohtua,quarteroniRB2016,Rozza2008229} or the Proper Generalized Decomposition (PGD) \cite{ChinestaEnc2017,Chinesta2011,Dumon20111387}) are available.

The offline phase of the POD--Galerkin ROM consists in an exploration of the solution manifold, obtained by querying the CutFEM high fidelity solver for $M \gg 0$ values of the parameter, and a successive compression to a basis of size $N < M$. The procedure starts by collecting the $M$ high fidelity solutions in the so-called snapshots matrix $\mathcal{S}$, defined as
\begin{equation}
{\mathcal{S}} = [\widehat{u}(\mu^1),\dots,\widehat{u}(\mu^{M})] \in \mathbb{R}^{\mathcal{N}_h\times M}.
\label{eq:snapmat}
\end{equation}
Here $\mu^1, \hdots, \mu^M$ are randomly selected values in $\mathcal{K}$ and ${u}(\mu^i)$ corresponding solutions of \eqref{eqn:hf_lapl}, $i = 1, \hdots, M$. 
We note that, since the solution ${u}(\mu^i)$ of \eqref{eqn:hf_lapl} is sought in the $\mu$-dependent FE space \eqref{eqn:V_poisson}, at least a suitable extension $\widehat{u}(\mu^i)$ should be carried out to provide snapshots defined on the (common) background mesh $\mathcal{B}_h$. 
Such extensions define a snapshot on the ($\mu$-independent) background FE space
\begin{equation*}
\widehat{V}_h=\left \{  \upsilon \in C^0({\overline {\mathcal B}})\,:\,\upsilon |_K  \in P^1(K), \, \forall K\in \mathcal{B}_h  \right \},
\end{equation*}
denoting by $\mathcal{N}_h$ its dimension.
We will return to this topic in Sections \ref{sec:extension} and \ref{sec:transport}.

Afterwards, following e.g. \cite{HeRoSta16,Kunisch2002492}, a compression by POD is carried out. In practice, this derives the following eigenvalue problem:
\begin{gather}
{\mathcal{C}}{Q} = {Q}{\Lambda} ,\qquad\text{being }\mathcal{C}_{ij} = \langle{\widehat{u}(\mu^i),\widehat{u}(\mu^j)}\rangle_{L^2{({\mathcal T}_{h})}},\, \mbox{ for } i,j = 1,\dots,M ,\nonumber
\end{gather}
where ${\mathcal{C}}$ is the correlation matrix obtained starting from the snapshots ${\mathcal{S}}$, $\Lambda$ is a diagonal matrix collecting eigenvalues on the diagonal, and ${Q}$ is an orthogonal matrix of the corresponding eigenvectors. 
The resulting reduced space is then spanned (possibly after a $L^2$ normalization) by the columns of the matrix obtained as product between $\mathcal{S}$ and the first $N$ columns of $Q$. The basis functions are denoted by $\varphi_i$, $i = 1, \hdots, N$, as well as
\begin{equation*}
\widehat{V}_N = \text{span}\{\varphi_1, \hdots, \varphi_N\}
\end{equation*}
denotes the obtained $N$-dimensional space, which will replace the high fidelity space $V_h(\mu)$ in all online computations. We remark that $\widehat{V}_N$ is a parameter independent space, composed of basis functions defined on the whole background mesh; the former assumption will be relaxed in Section \ref{sec:transport}.

During the online stage, a reduced solution $u_N(\mu)$ of the form
\begin{equation}\label{eq:aprox_fields}
{u_N}(\mu) = \sum_{i=1}^{N} \alpha_i(\mu) \ {\varphi_i},
\end{equation}
is sought\footnote{As the basis functions are global, ${u_N}(\mu)$ is defined on the whole background mesh, understanding that its value outside of ${\mathcal D}_{\mathcal{T}}(\mu)$ is not interesting and can be discarded during the analysis of the numerical results.}. The unknown coefficients $\bm{\alpha} = [\alpha_1, \hdots, \alpha_N] \in \mathbb{R}^N$ are obtained through a Galerkin projection of the governing equations onto the reduced basis space, as follows: for any $\mu \in \mathcal{K}$, find $\bm{\alpha} = \bm{\alpha}(\mu) \in \mathbb{R}^N$ such that
\begin{equation}
\sum_{i=1}^{N} \alpha_i(\mu) a_\gamma(\varphi_i, \varphi_j; \mu) = \ell_\gamma(\varphi_j; \mu), \quad \forall j = 1, \hdots, N.
\label{eq:online_system}
\end{equation}
Provided that $N \ll \mathcal{N}_h$, this linear system is usually very inexpensive to solve. The overall efficiency of the online stage usually comes also from an inexpensive assembly of the left-hand and right-hand side of this linear system, that can be obtained owing to affinity assumptions. Such assumptions, although natural in many applicative cases, can be approximately regained (if lacking) through suitable hyper-reduction procedures \cite{BARRAULT2004667}. We do mention that numerical examples in Section \ref{sec:num_exp} do not fulfill affinity assumptions. However, being this a preliminary work on the combination of unfitted methods and ROMs, we do not aim at the utmost efficiency at this point. Thus, we will allow inefficient assembly of the left-hand and right-hand side by projecting $\mu$-dependent forms during the online stage. Further perspectives and improvements in this direction will be discussed in the conclusion in Section \ref{sec:conclusions}.

\subsection{Steady Stokes problem}\label{subsec_POD_theory_Stokes}
Let us now discuss the reduction of the Stokes problems presented in Section \ref{sec:HFStokes}.
During the offline stage, as in the previous discussion, we collect snapshots from the CutFEM high fidelity discretization. Since the problem is characterized by two unknowns, following e.g. \cite{ballarin2015supremizer,stabile_stabilized} we define two snapshots matrices, namely
\begin{gather*}
{\mathcal{S}_u} = [\widehat{u}(\mu^1),\dots,\widehat{u}(\mu^{N_s})] \in \mathbb{R}^{\mathcal{N}_h^u\times M},\qquad
{\mathcal{S}_p} = [\widehat{p}(\mu^1),\dots,\widehat{p}(\mu^{N_s})] \in \mathbb{R}^{\mathcal{N}_h^p\times M},
\end{gather*}
and ultimately perform two separate compressions by means of POD. We remark that two different extension operators, both denoted by $\widehat{\cdot}$ above, may be used for velocity and pressure snapshots, respectively. 
Here $\mathcal{N}_h^u$ and $\mathcal{N}_h^p$ denote the dimension of the background FE spaces
\begin{align*}
\widehat{V}_h &= \left\{ \upsilon_h \in (C^0 (\overline{\mathcal{B}}))^{d} : \upsilon_h|_K \in ({P}^1 (K))^{d}, \forall K \in \mathcal{B}_h\right\},\\
\widehat{Q}_h &= \left\{ w_h \in C^0 (\overline{\mathcal{B}}) : w_h|_K \in  {P}^1 (K), \forall K \in \mathcal{B}_h\right \},
\end{align*}
respectively.
Let us denote by $\varphi_i^u$, $i = 1, \hdots, N$ the first $N$ velocity basis functions, as well as by $\varphi_j^p$, $j = 1, \hdots, N$ the first $N$ pressure basis functions.

However, as well known in the reduced basis approximation of saddle point problems, the resulting basis may lead to inaccurate results, especially for what concerns the pressure approximation. Therefore, we resort to a supremizer enrichment procedure \cite{RoVe07,Gerner2012,ballarin2015supremizer}, which we summarize in the following. From a practical standpoint, this requires the solution of $M$ additional Laplace problems during the offline phase which, in strong form, read as: find $s(\mu^i) \in V_h(\mu^i)$ such that
\begin{equation}\label{eq:sup_problem}
\begin{cases}
- \Delta s(\mu^i) = \nabla p(\mu^i) & \mbox{ in } \mathcal D(\mu^i),\\
s(\mu^i) = 0 & \mbox{ on } \Gamma_D(\mu^i),\\
\nabla s(\mu^i) \cdot {\bf n}_{\Gamma} = 0 & \mbox{ on } \Gamma_N(\mu^i),\\
\end{cases}
\end{equation}
for all $i = 1, \hdots, M$.
We note that, with the only exception of being a vector problem, \eqref{eq:sup_problem} is a special instance of \eqref{P:Poisson}. Therefore, we can use the CutFEM discretization introduced in Section \ref{sec:HFDarcy} for its approximation, component by component. The obtained snapshots are then collected in 
\begin{gather*}
{\mathcal{S}_s} = [\widehat{s}(\mu^1),\dots,\widehat{s}(\mu^{M})] \in \mathbb{R}^{\mathcal{N}_h^u\times M},
\end{gather*}
where the same extension (and the same extended FE space) used for velocity is also used for supremizers; correspondingly, $N$ basis functions are generated, denoted by $\varphi_i^s$, $i=1, \hdots, N$ in the following. The resulting reduced basis spaces are therefore obtained as
\begin{align*}
\widehat{V}_N = \text{span}\{\varphi_1^u, \hdots, \varphi_N^u, \varphi_1^s, \hdots, \varphi_N^s\},\qquad
\widehat{Q}_N = \text{span}\{\varphi_1^p, \hdots, \varphi_N^p\}.
\end{align*}

Using the combined velocity-pressure notation, we can express the combined reduced basis space as
\begin{align*}
\widehat{\bm{V}}_N = \text{span}\{\varphi_1^u, \hdots, \varphi_N^u, \varphi_1^s, \hdots, \varphi_N^s, \varphi_1^p, \hdots, \varphi_N^p\};
\end{align*}
as in the following we will not be interested in differentiating the notation between velocity, supremizer and pressure, we equivalently write this space as
\begin{align*}
\widehat{\bm{V}}_N = \text{span}\{\varphi_1, \hdots, \varphi_{3 N}\},
\end{align*}
owing to a suitable renumbering of the basis functions. Thus, during the online stage we seek a $3N$-dimensional solution $\bm{u}_N(\mu)$ of the form
\begin{equation*}
{\bm{u}_N}(\mu) = \sum_{i=1}^{3 N} \alpha_i(\mu) \ {\varphi_i}
\end{equation*}
such that
\begin{equation*}
\sum_{i=1}^{3N} \alpha_i(\mu) a_\gamma(\varphi_i, \varphi_j; \mu) = \ell_\gamma(\varphi_j; \mu), \quad \forall j = 1, \hdots, 3 N.
\end{equation*}

Finally, we remark that alternative approaches to supremizer enrichment are available, see e.g. \cite{shafqat,Caiazzo2014598,Giacomo2018,stabile_stabilized}. In particular, one that relies on the underlying $P^1/P^1$ (Cut) FEM stabilization \cite{shafqat,Giacomo2018} seemed attractive in view of providing a less intrusive reduced order model, allowing to neglect the supremizer enrichment stage.
However, numerical experiments carried out in the preparation of this work have shown a sensible deterioration of the accuracy when relying only on the high-fidelity stabilization without supremizer enrichment due to the weaker stabilization employed. Although, stronger CutFEM stabilizations are available in the literature, \cite{BuErnFe17,Bu07}.

\subsection{Snapshots extension to the background mesh}\label{sec:extension}
Let us now discuss a few practical options for what concerns snapshot extension. 
The first two possibilities that we investigate are to carry out \emph{trivial extensions}. 
The first option, which will be called \emph{zero extension} in the following, is to extend the snapshot to zero in $\mathcal{B} \setminus \mathcal{D}(\mu)$. A slightly different option is to extend snapshots to zero in $\mathcal{B} \setminus \mathcal{D}_{\mathcal{T}}(\mu)$, as the CutFEM solution is naturally defined up to the boundary of $\mathcal{D}_{\mathcal{T}}(\mu) \supset \mathcal{D}(\mu)$. This second options will be named \emph{natural smooth extension}, as this name is usually employed in the CutFEM literature \cite{BeBuHa09}.

The third method we propose is an \emph{harmonic extension}, see e.g. \cite{BallarinRozzaMaday2017}, which, for the scalar case of Section \ref{subsec_POD_theory}, reads: for any $\mu \in \mathcal{K}$, assuming the corresponding solution $u(\mu)$ of \eqref{P:Poisson} to be known, find $u_c(\mu): \mathcal{B} \to \mathbb{R}$ such that
\begin{equation}
\begin{cases}
 -\Delta  u_c(\mu) = 0, & \text{in $\mathcal{B} \setminus \mathcal{D}(\mu)$},
 \\
  u_c(\mu)=u(\mu), & \text{on $\partial \mathcal{D}(\mu)$},
  \\
  u_c(\mu)=0, &\text{on $\partial\mathcal{B}$},
\end{cases}
\label{eq:harmonic}
\end{equation}
and then define
\begin{equation*}
\widehat{u}(\mu) = 
\begin{cases}
 u(\mu), & \text{in $\mathcal{D}(\mu)$},
 \\
 u_c(\mu), & \text{in $\mathcal{B} \setminus \mathcal{D}(\mu)$}.
\end{cases}
\end{equation*}
Being another elliptic problem, the harmonic extension problem \eqref{eq:harmonic} is discretized as described in Section \ref{sec:HFDarcy}, up to a suitable replacement of domain, boundaries and problem coefficients. A similar approach may be used for the Stokes problem in Section \ref{sec:HFStokes}, acting component by component.

\subsection{Snapshots transportation on the background mesh}\label{sec:transport}
Similarly to problems characterized by a strong hyperbolic nature, the methodology described up to now may suffer from a slowly decreasing Kolmogorov $n$-width. To exemplify, let us consider the natural smooth extension, and let $\mu^1 \neq \mu^2$ be two parameters such that $\overline{\mathcal{D}}_{\mathcal{T}}(\mu^1) \cap \overline{\mathcal{D}}_{\mathcal{T}}(\mu^2) = \emptyset$. Then, the snapshot corresponding to $\mu^1$ bears no useful information for the reduced basis representation of the phenomena associated to $\mu^2$, as the support of the solutions do not intersect and the extension is trivial. Even though a non-trivial extension (such as the harmonic one) may alleviate this drawback to some extent, we propose in this Section to combine extension with an additional snapshots preprocessing based on a transportation, as advocated in {{\cite{BeIoRi2018,CaMaSta19,IolloLombardi2014,Naira2017,Reiss2015,Welper2017}}.

Let $\boldsymbol{\tau}(\mu): \mathcal{B} \to \mathcal{B}$ be a bijective map chosen such that $\mathcal{D}(\mu)$ is mapped into $\mathcal{D}(\overline{\mu})$, and let $\boldsymbol{\tau}^{-1}(\mu)$ denote its inverse with respect to the spatial coordinates. Here $\overline{\mu}$ is a fixed parameter, e.g. the barycenter of $\mathcal{K}$. We will discuss how to devise $\boldsymbol{\tau}(\mu)$ from $\mathcal{D}(\mu)$ later on in this Section.

Let $u(\mu)$ be the snapshot obtained having solved \eqref{P:Poisson} (the same reasoning goes for Stokes), and extend it to $\widehat{u}(\mu)$ according to the one of the methods introduced in the previous Section. By composition with $\boldsymbol{\tau}(\mu)$ the following transported (and extended) snapshot is defined:
\begin{equation*}
\widehat{u}_{\boldsymbol{\tau}}(\mu) = \widehat{u}(\mu) \circ \boldsymbol{\tau}(\mu).
\end{equation*}
Such transported snapshots $\widehat{u}_{\boldsymbol{\tau}}(\mu^i)$ can then be employed in the definition of the snapshots matrix \eqref{eq:snapmat} in place of the extended ones $\widehat{u}(\mu^i)$, $i = 1, \hdots, M$. Since all transported snapshots are centered around $\mathcal{D}(\overline{\mu})$, it is expected that a more effective reduction (i.e., faster decay for POD singular values) is obtained.

However, as a result of this transformation, the reduced basis space needs now to depend on $\mu$. Indeed, for any online query associated to a new parameter $\mu$, basis functions need to be transported back from a neighborhood of $\mathcal{D}(\overline{\mu})$ to the corresponding neighborhood of $\mathcal{D}(\mu)$. In particular, the $\mu$-dependent reduced basis space is now defined as
\begin{equation*}
\widehat{V}^{\boldsymbol{\tau}}_N(\mu) = \text{span}\{\varphi_1 \circ \boldsymbol{\tau}^{-1}(\mu), \hdots, \varphi_N \circ \boldsymbol{\tau}^{-1}(\mu)\},
\end{equation*}
and inverse transported basis functions $\varphi_j \circ \boldsymbol{\tau}^{-1}(\mu)$, $j = 1, \hdots, N$ are to be employed in the solution of the reduced linear system \eqref{eq:online_system}, as well as in the representation \eqref{eq:aprox_fields} of the reduced order solution.
We do remark that such inverse transport introduces an additional non-affinity in the problem formulation; as discussed in Section \ref{subsec_POD_theory}, such issue will be a topic of future developments.

Let us now go back to the choice of $\boldsymbol{\tau}(\mu)$. As customary in CutFEM methods, $\mathcal{D}(\mu)$ is obtained through a level set method{\blue\footnote{\blue Although, we remark that the use of a level-set domain description is not an inherent limitation of CutFEMs, see for example \cite{JONSSON2017366} where a parametric description of the domain boundary is used.}}. 
Being the level sets functions employed in the numerical experiments simple algebraic expressions, it is a matter of simple algebraic manipulations to devise the associated maps $\boldsymbol{\tau}(\mu)$. In future, in case of lacking explicit representation of $\boldsymbol{\tau}(\mu)$, we seek to resort to optimal transportation procedures \cite{BeCaCuNePe15,BeIoRi2018,SolGoesPeyCuButNguyenDuGuibas15}, as proposed in \cite{IolloLombardi2014,NoBaRo19}.

Before concluding this Section, we finally stress the difference with a reference domain formulation, as both methods revolve around the idea of introducing a fixed parameter $\overline{\mu}$. In the reference domain formulation, the domain $\mathcal{D}(\overline{\mu})$ associated to fixed parameter $\overline{\mu}$ is the only one employed in the computations, and the differential problems associated to any $\mu$ is recast on $\mathcal{D}(\overline{\mu})$ through suitable pull backs, which are often limited to small deformations.
In contrast, the proposed transport procedure requires only a postprocessing of the CutFEM solution, each obtained on its own domain $\mathcal{D}(\mu)$, without any change to the high fidelity solver. The introduction of $\overline{\mu}$ is thus mandatory in the reference domain approach and it is strongly intertwined with the snapshots computation, whilst in our case it is associated only to a desirable postprocessing carried out in order to improve the decay of the Kolmogorov $n$-width.

\section{Numerical experiments}\label{sec:num_exp}
\subsection{The Darcy flow pressure model}\label{sec:num_exp_scalar}
The proposed reduced order technique is numerically tested in this Section on a case characterized by large deformations. The domain $\mathcal{D}(\mu) \subset \mathbb{R}^2$ is a parametrized ellipse, defined through the level set function
\begin{equation*}
\phi(x,y;\mu_1, \mu_2, \mu_3, \mu_4) = \mu^2_2(x-\mu_3)^2 +\mu^2_1(y-\mu_4)^2 -\mu^2_1 \mu^2_2 R^2,
\end{equation*}
where the reference radius $R = 0.05$, the length of the axes of the ellipse is parametrized by $(\mu_1, \mu_2) \in [0.3, 1.8]^2$, while the position of the center of the ellipse is parametrized by $(\mu_3, \mu_4) \in [-0.85, 0.85]^2$. A corresponding background domain $\mathcal{B} = [-1.2, 1.2]^2$ is chosen so that the ellipse is strictly contained in $\mathcal{B}$ for any $\mu = (\mu_1, \mu_2, \mu_3, \mu_4)$ in the parametric range $\mathcal{K} = [0.3, 1.8]^2 \times [-0.85, 0.85]^2$.
The value $\overline{\mu} = (1, 1, 0, 0)$, corresponding to a circle of radius $R$ centered in the origin, is chosen for what concerns the transport method introduced in Section \ref{sec:transport}, and the transportation mapping\footnote{We do note that this expression does not map $\mathcal{B}$ in itself. However, as the goal of the mapping is to transport snapshots in a $R$-neighborhood of the origin, we are satisfied with the proposed mapping as it provides properly mapped values in such neighborhood. Mapped snapshots will be extended to zero for points which are not in $\boldsymbol{\tau}(\mu)(\mathcal{B})$. This is trivially through in the zero and natural smooth extensions, while it is in agreement with boundary conditions on $\partial\mathcal{B}$ for the harmonic extension. As an alternative, one could flip the role of extension and transportation, first transporting snapshots and then extending the transported snapshots; we omit this case, as it only affects the harmonic extension case with negligible impact on the overall results and discussion.}  $\boldsymbol{\tau}(x,y;\mu)=(\mu_1x +\mu_3, \mu_2y + \mu_4)$ together with the inverse mapping $\boldsymbol{\tau}^{-1}(x,y;\mu)=\left(\frac{x-\mu_3}{\mu_1}, \frac{y-\mu_4}{\mu_2}\right)$ is considered{{\footnote{\blue
In particular, we transport the elliptical geometries (corresponding to the snapshots for various parameters $\mu$ during the offline stage) to circular ones with fixed radius. Indeed, the transformation $x\mapsto \mu_1x +\mu_3$, $y\mapsto \mu_2y + \mu_4$ guarantees that the ellipse $\mu^2_2(x-\mu_3)^2 +\mu^2_1(y-\mu_4)^2 -\mu^2_1 \mu^2_2 R^2 = 0$ is mapped into the circle $x^2 + y^2 - R^2 = 0$. Afterwards, we collect the transported snapshots and we construct the basis. Then, for every new (in the online stage) geometry we apply the inverse trasportation $x\mapsto ({x-\mu_3})/{\mu_1}$, $y\mapsto ({y-\mu_4})/{\mu_2}$, i.e., we transport basis functions defined on the circular geometry to the original elliptical geometry.}}
{\blue{, see \autoref{fig:Transport_Ellipse}}}.
\begin{figure} \centering
  \includegraphics[width=0.4\textwidth]{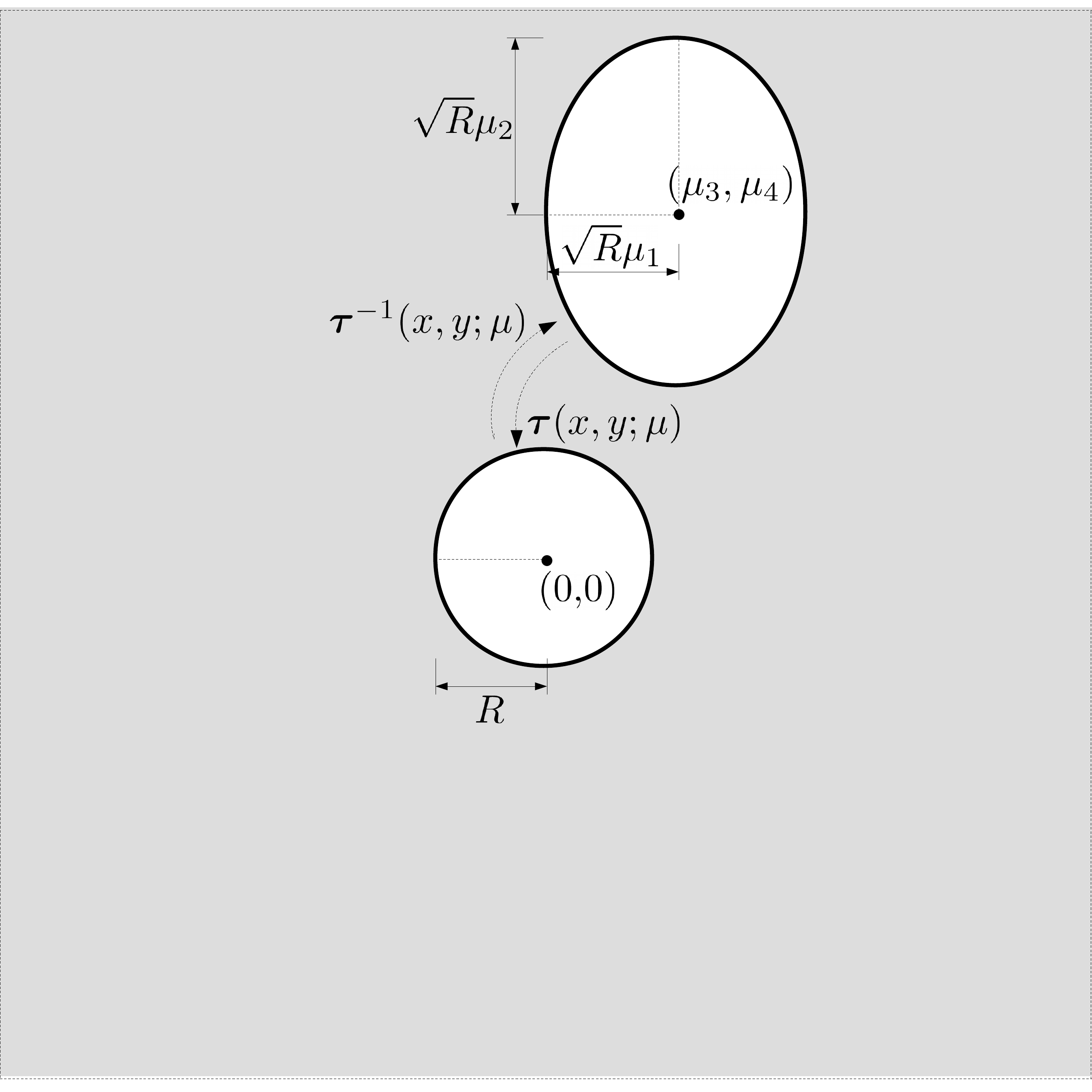}
  \caption{{\blue{{\emph{Darcy flow test case}}: 
The basic concept of the transportation mapping: ellipse transportation to a circle with fixed radius (through $\boldsymbol{\tau}(x,y;\mu)$), and reverse circle transportation  to the original ellipse geometry (through $\boldsymbol{\tau}^{-1}(x,y;\mu)$).
.}}}
  \label{fig:Transport_Ellipse}
\end{figure}
This example is characterized by large parametric deformations (e.g., translation of any coordinate of the center of the ellipse can reach values up to $1700\%$ its radius) which a boundary fitted method with reference domain formulation (obtained e.g. by creating a boundary fitted mesh to $\mathcal{D}(\overline{\mu})$ on $\mathcal{B}$ and applying a pull back to the differential problem) would hardly handle.
The strong formulation of the problem is as in \eqref{P:Poisson} for $g(x, y; \mu) = 20$, $g_D(x, y; \mu) = 0.5 + xy$, $\Gamma_D(\mu) \equiv \partial\mathcal{D}(\mu)$. For the background domain discretization we used mesh step size $h=0.05$, corresponding to 2806 degrees of freedom. {\blue{The parameter values $\gamma_D$, $\gamma_N$, $\gamma_1$ under consideration are $10$, $0$, $0.1$ respectively.}}

\begin{figure} \centering
  \includegraphics[width=0.7\textwidth]{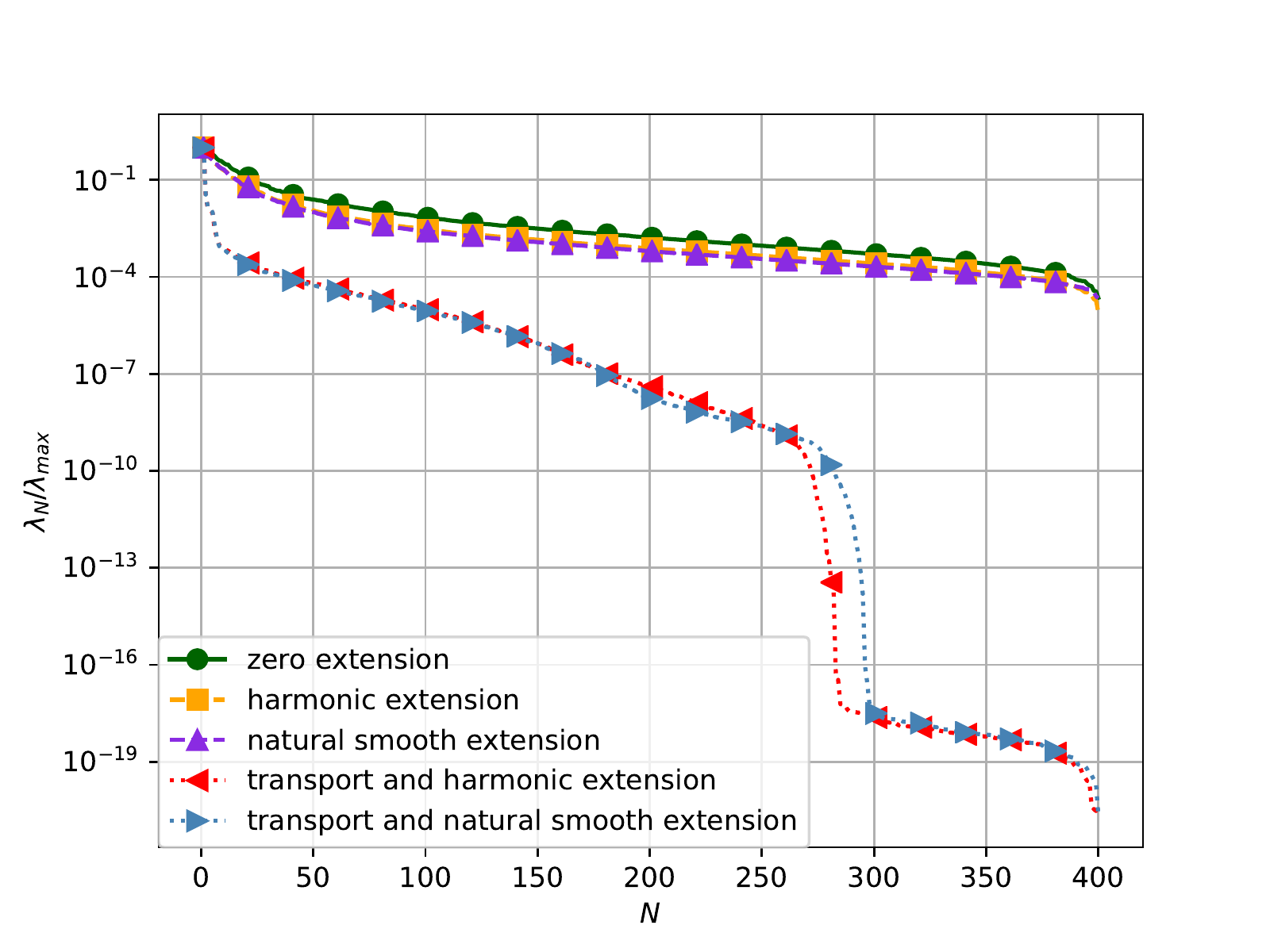}
  \caption{{\emph{Darcy flow test case}}: The POD eigenvalues decay (normalized to the maximum eigenvalue) for a training set of 400 snapshots is reported against the number of modes.}
  \label{fig:Poisson eigs}
\end{figure}

\begin{figure}
\centering
\includegraphics[width=0.32\textwidth]{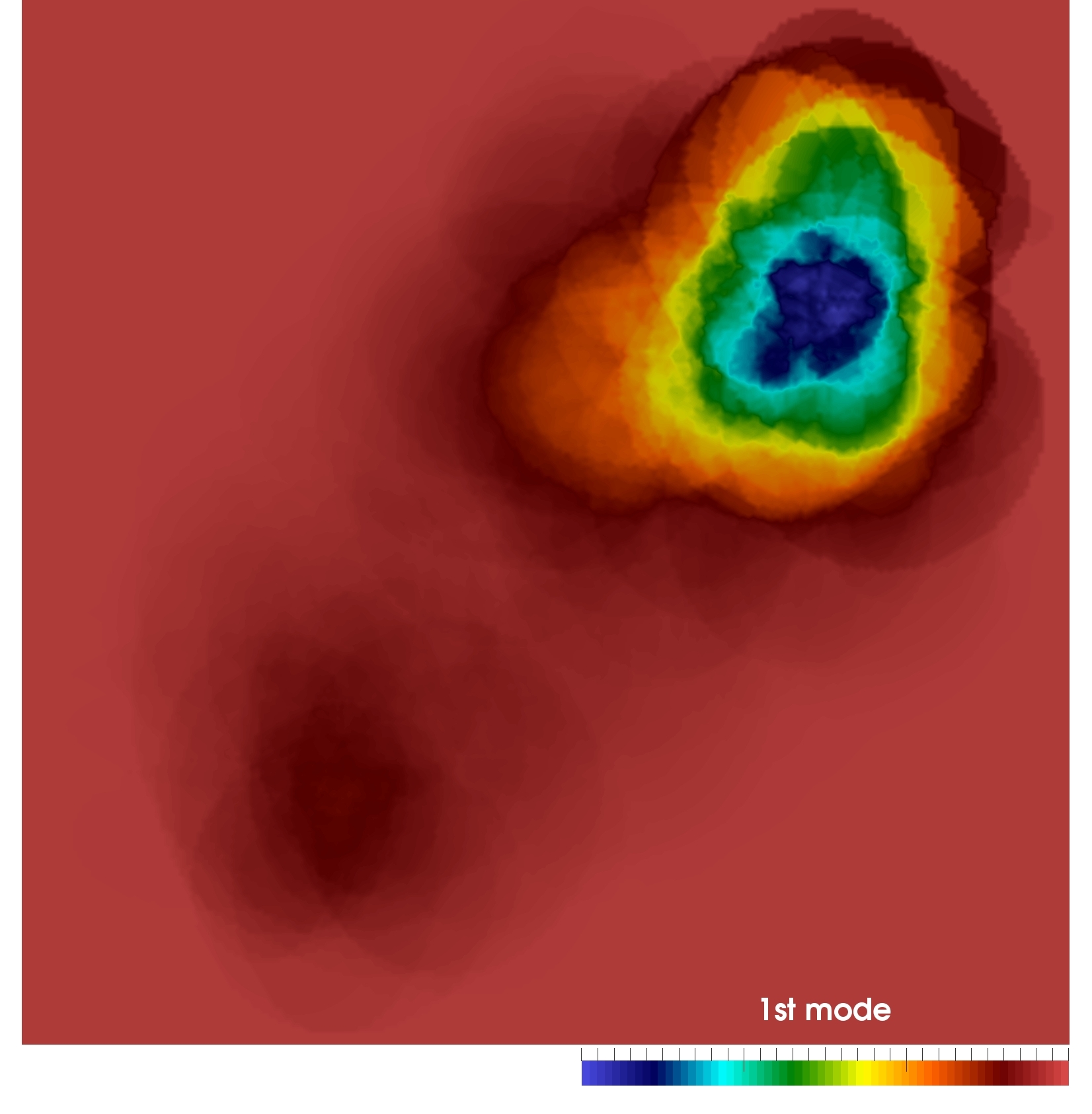}
\includegraphics[width=0.32\textwidth]{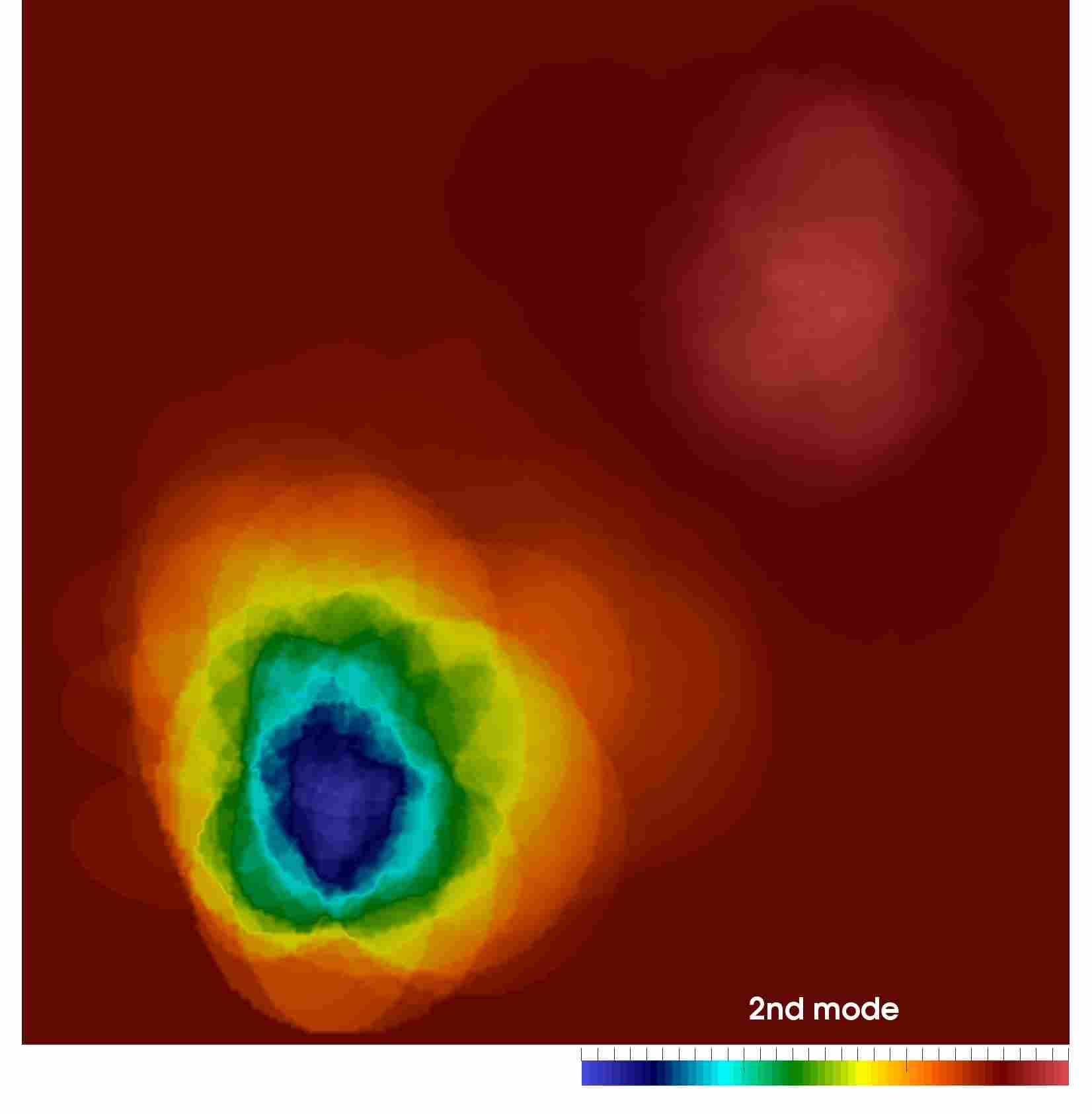}
\includegraphics[width=0.32\textwidth]{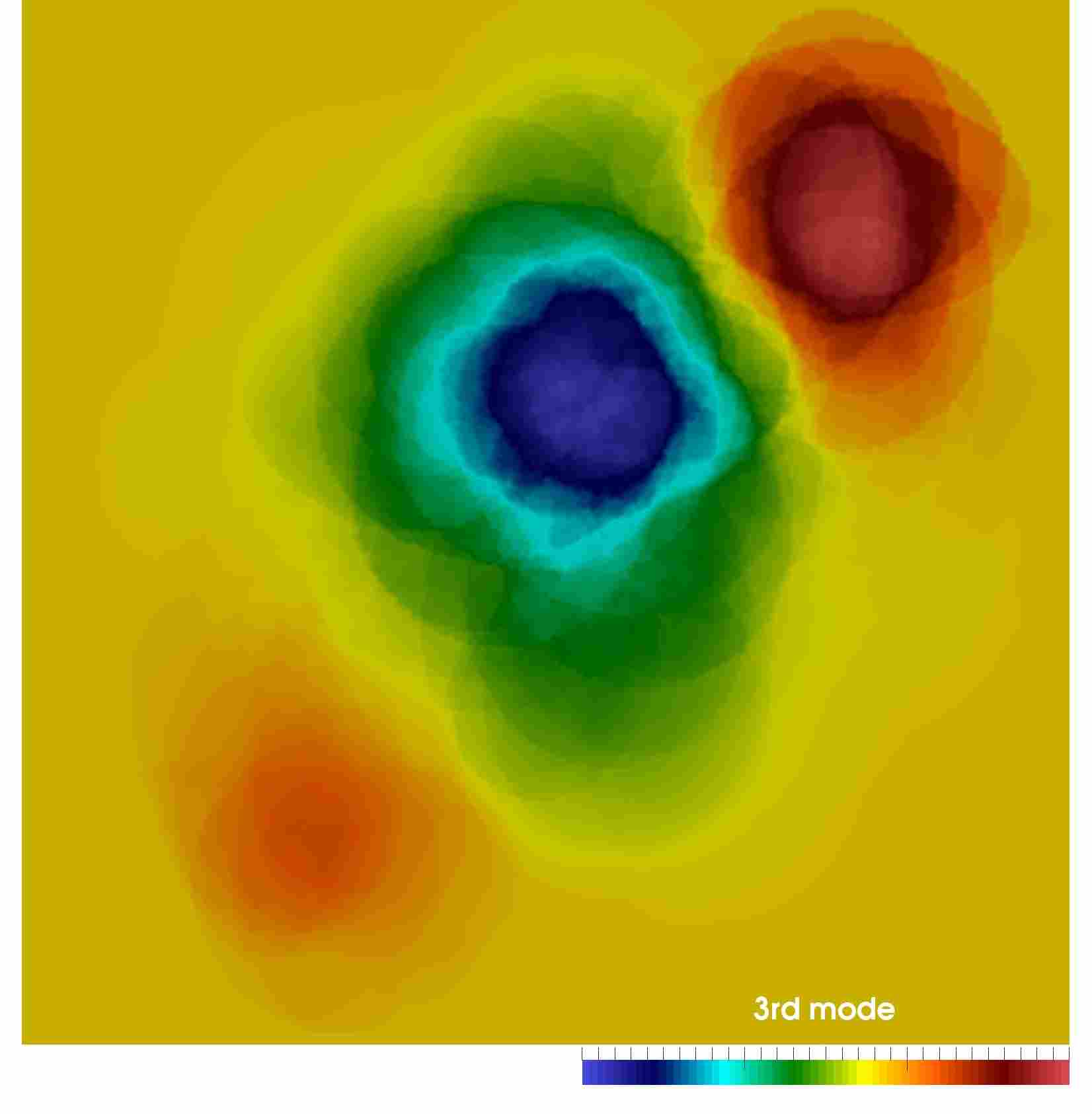}\\
\includegraphics[width=0.32\textwidth]{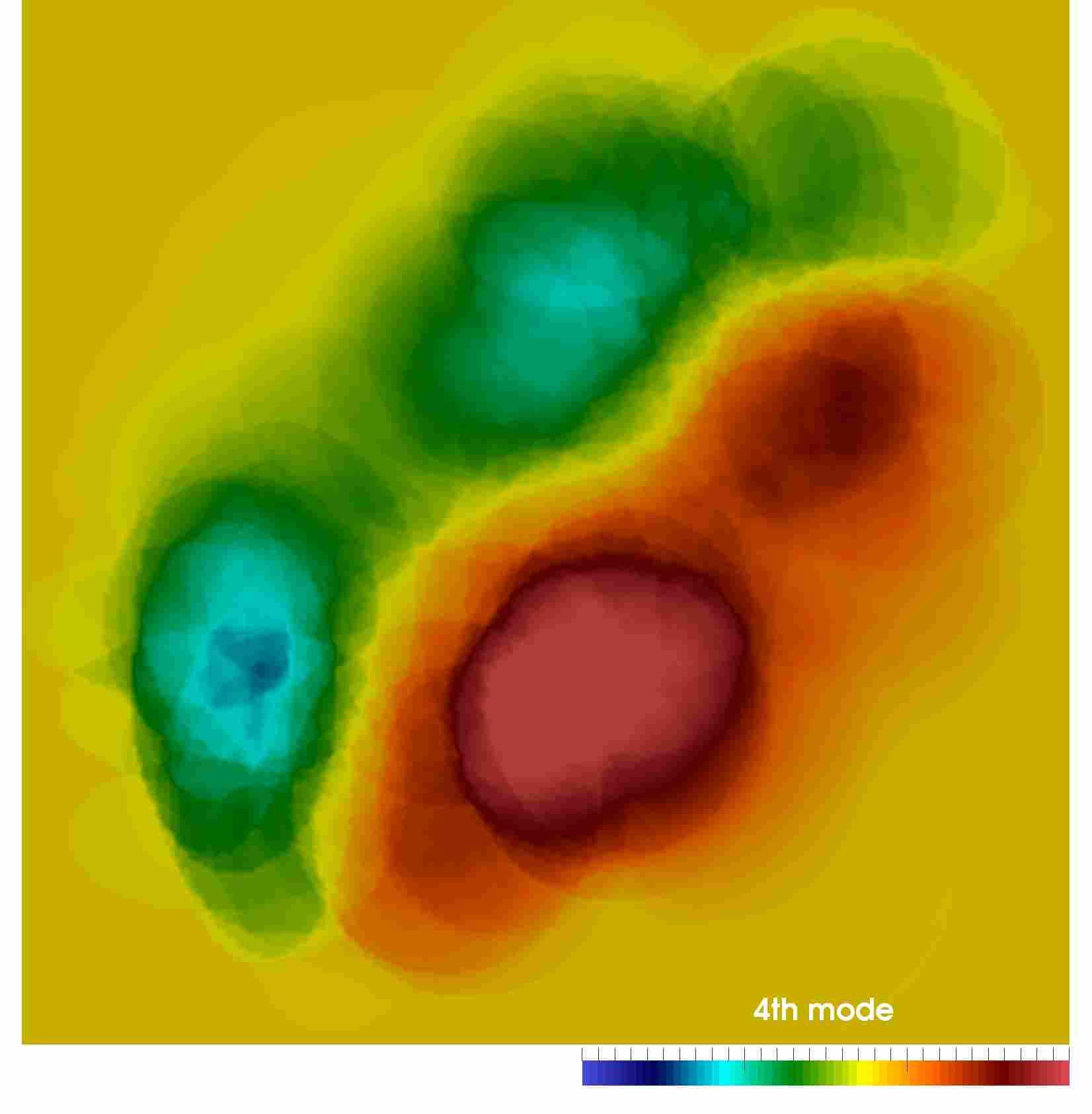}
\includegraphics[width=0.32\textwidth]{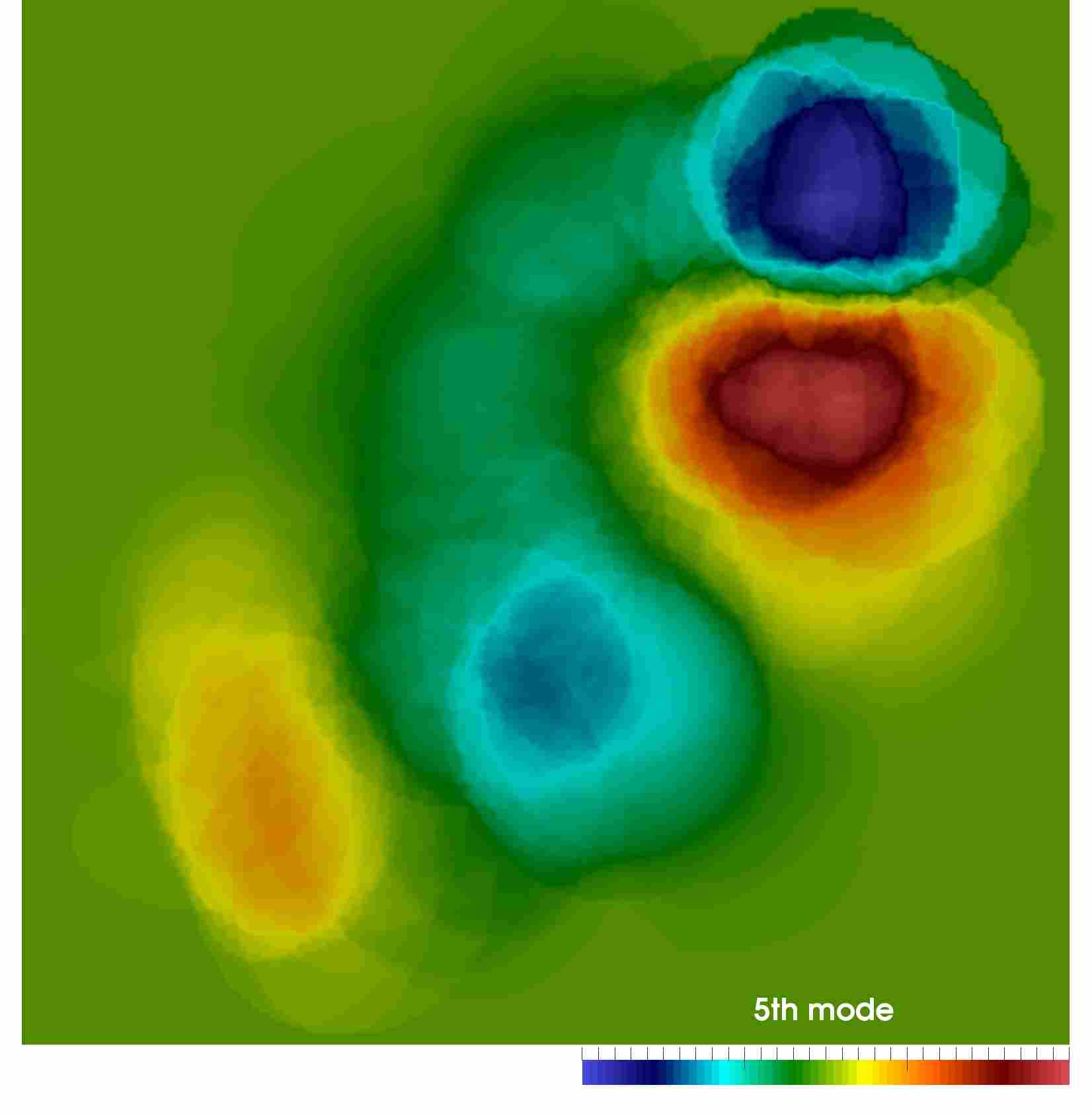}
\includegraphics[width=0.32\textwidth]{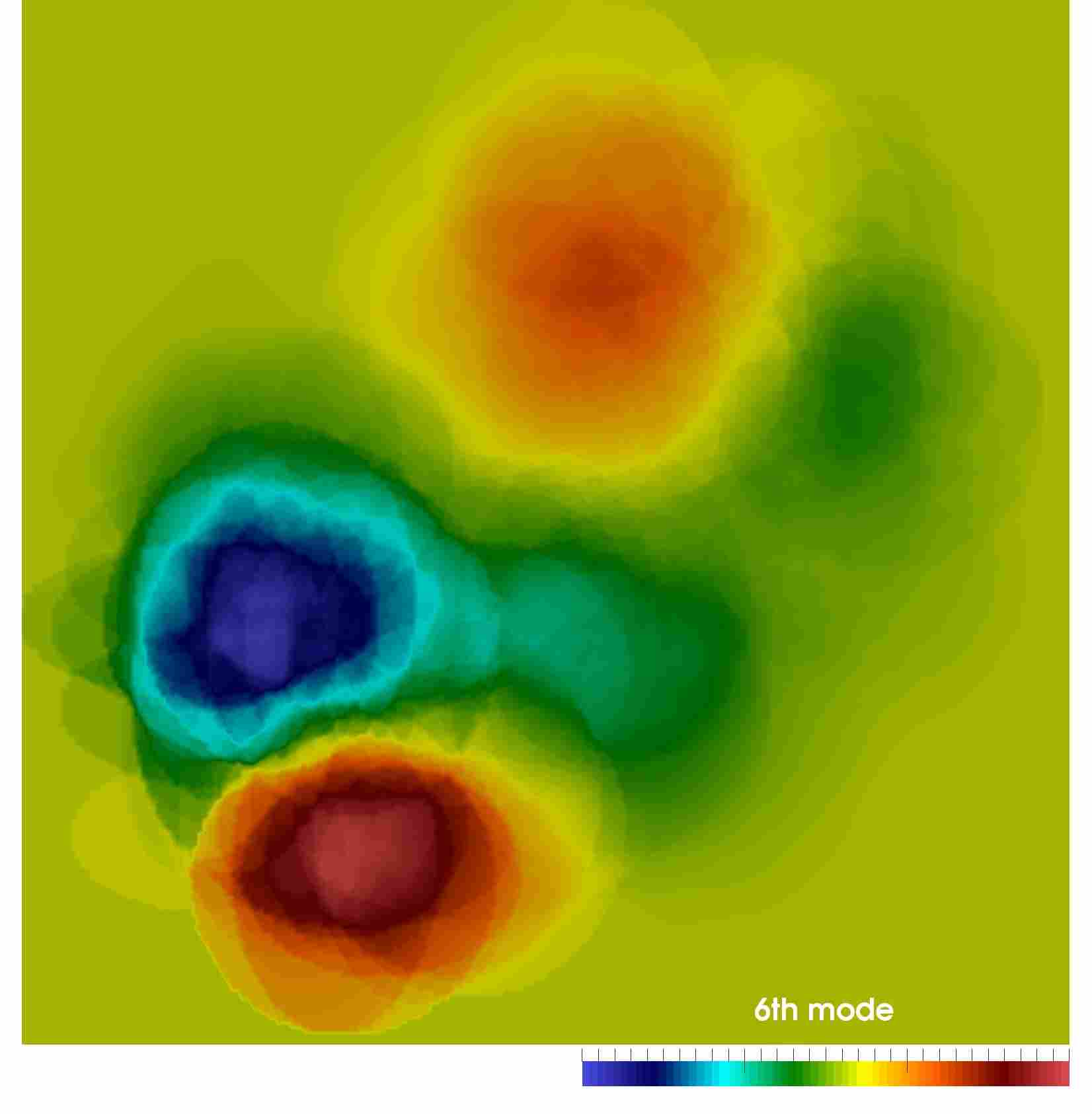}\\
\caption{\emph{Darcy flow test case}: first six POD modes for the natural smooth extension without transportation.}
\label{Fig:Poisson_no_transport_modes}
\end{figure}

\begin{figure}
\centering
\includegraphics[width=0.32\textwidth]{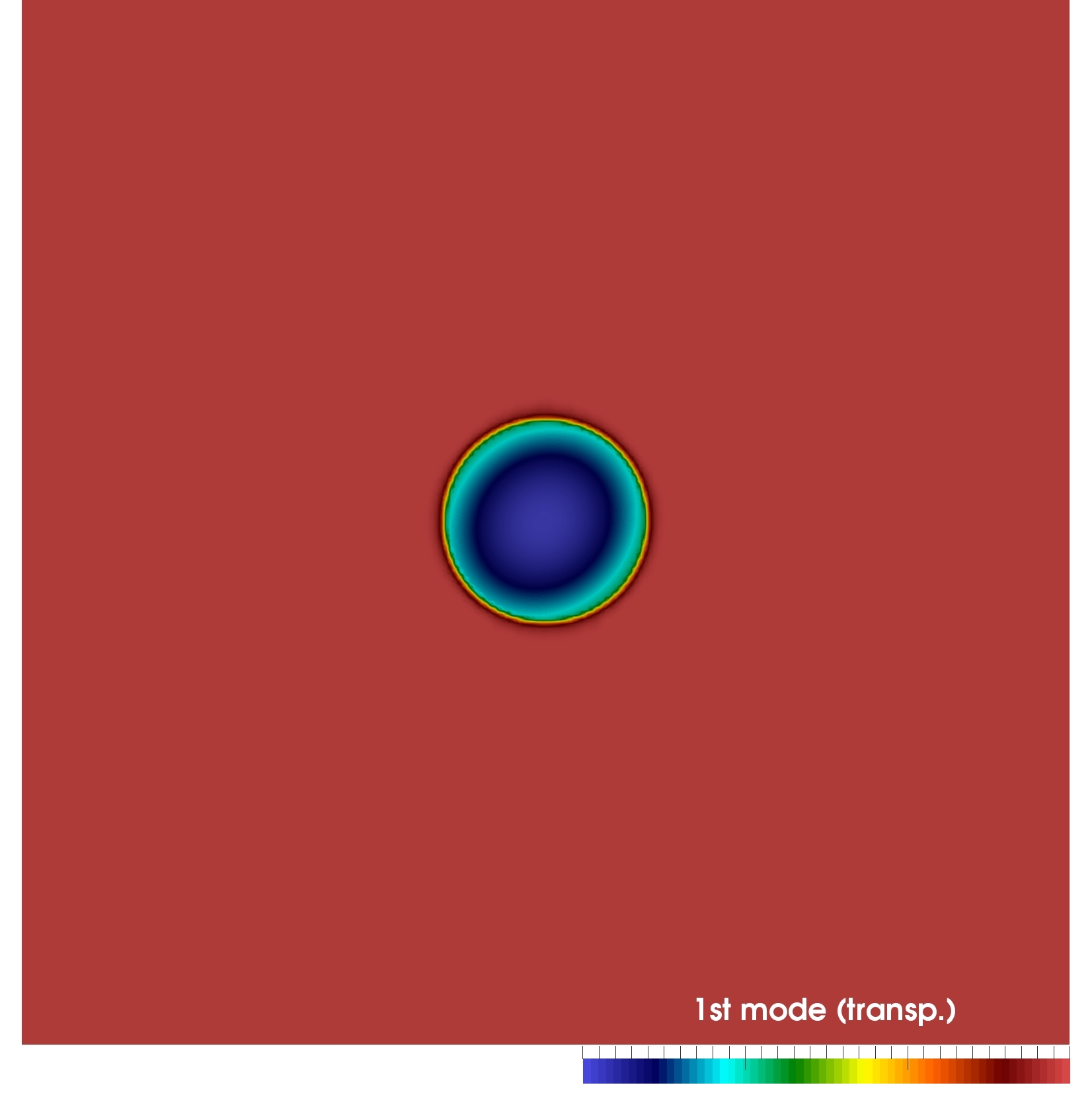}
\includegraphics[width=0.32\textwidth]{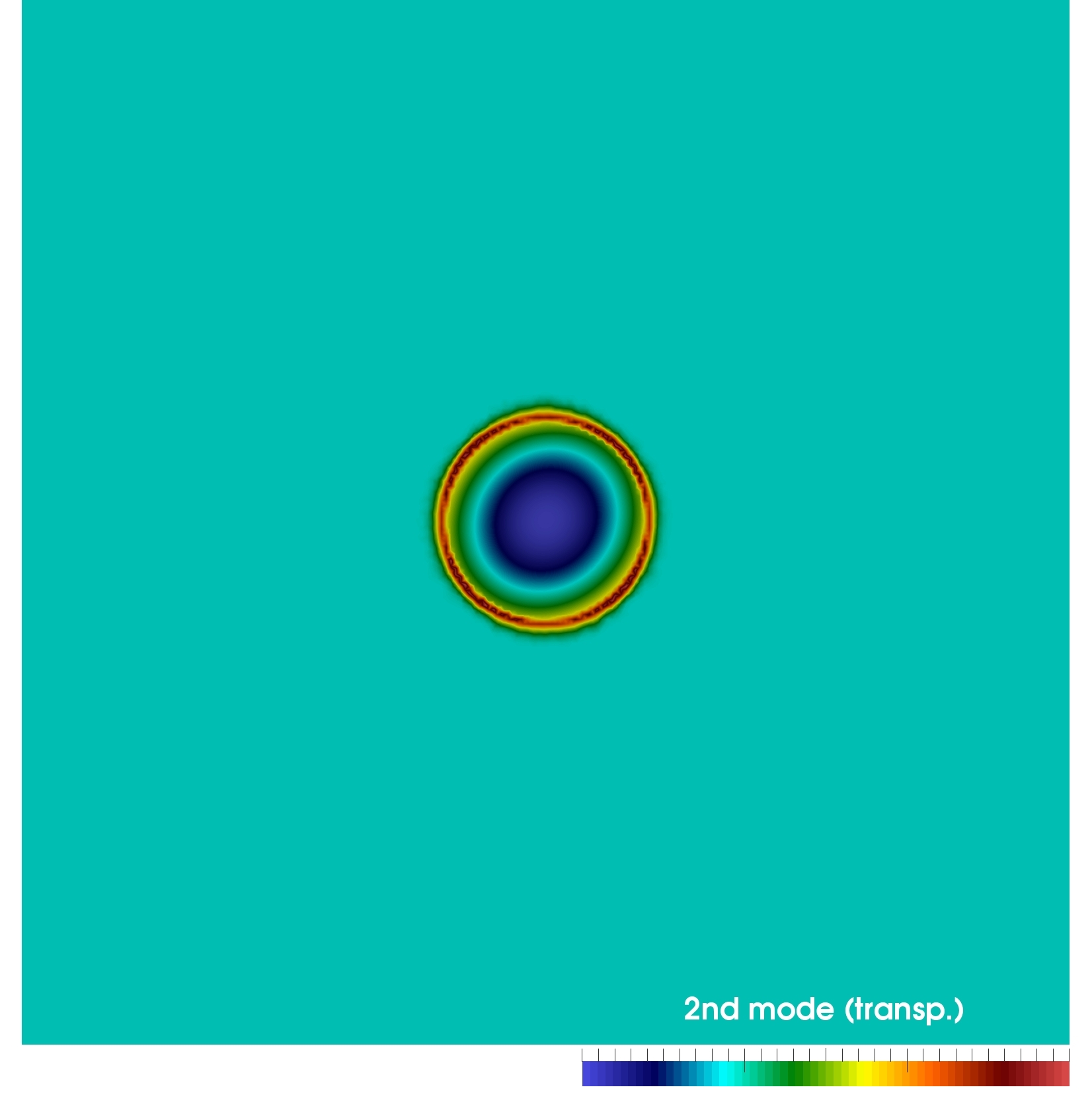}
\includegraphics[width=0.32\textwidth]{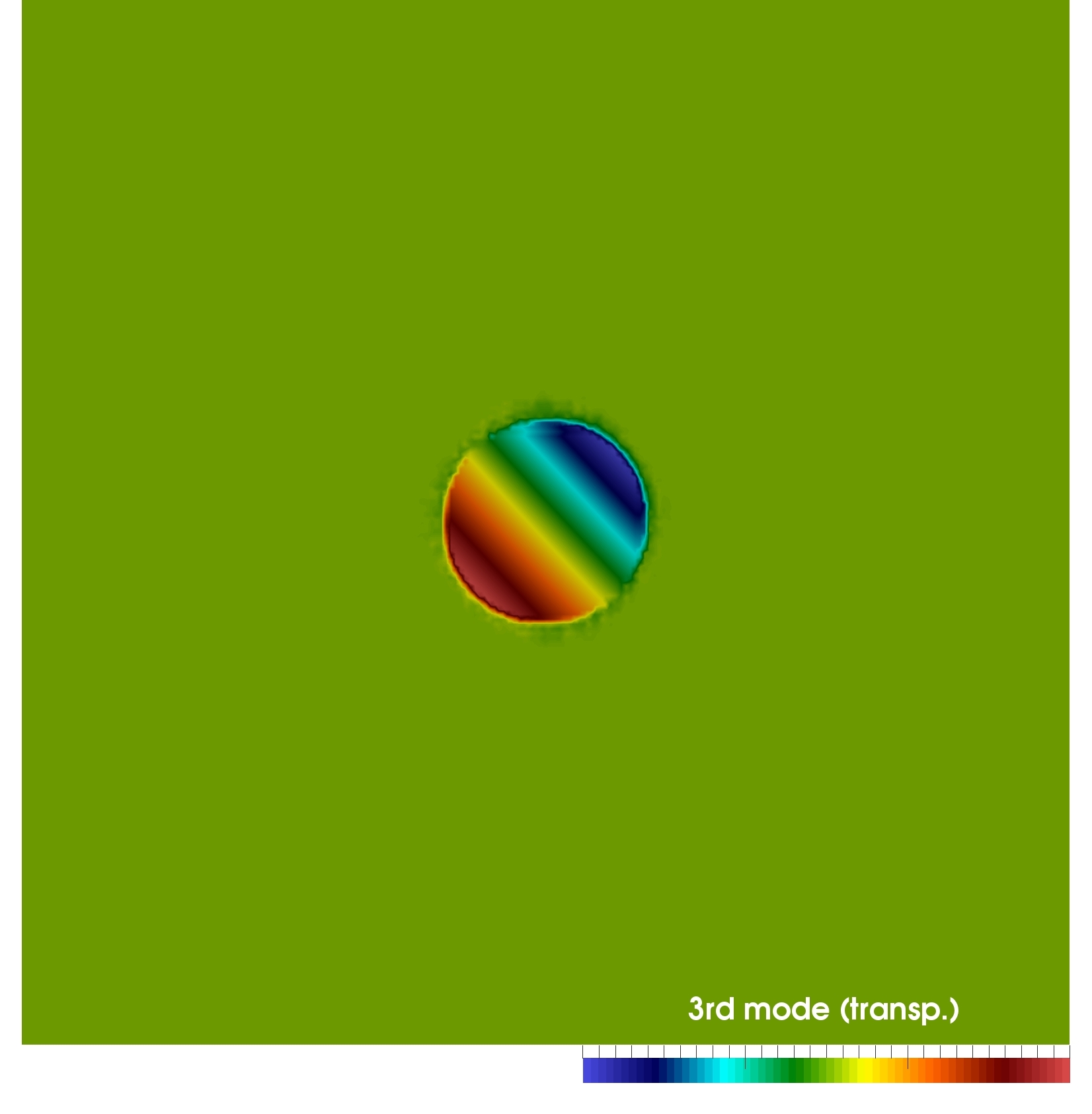}\\
\includegraphics[width=0.32\textwidth]{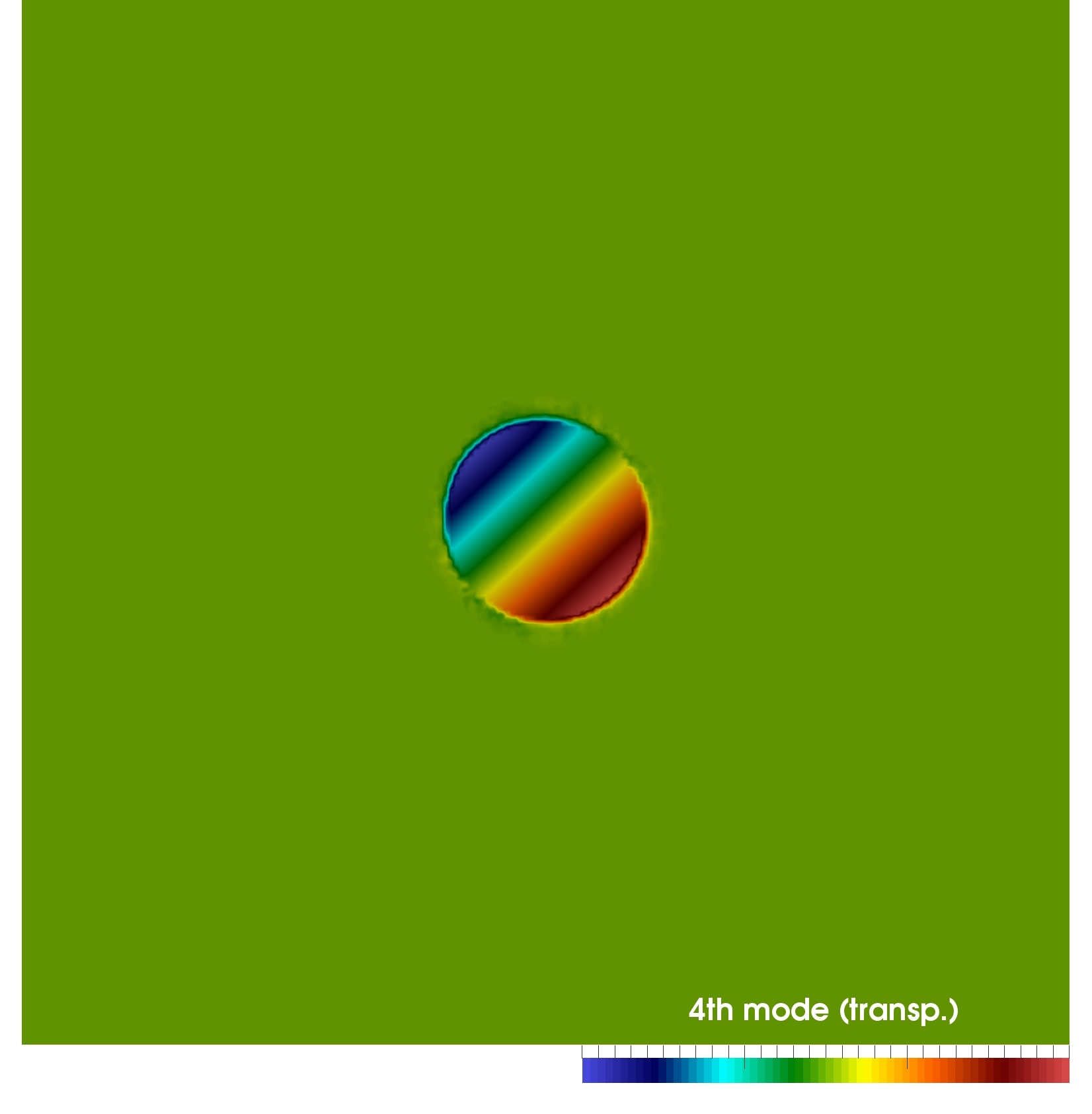}
\includegraphics[width=0.32\textwidth]{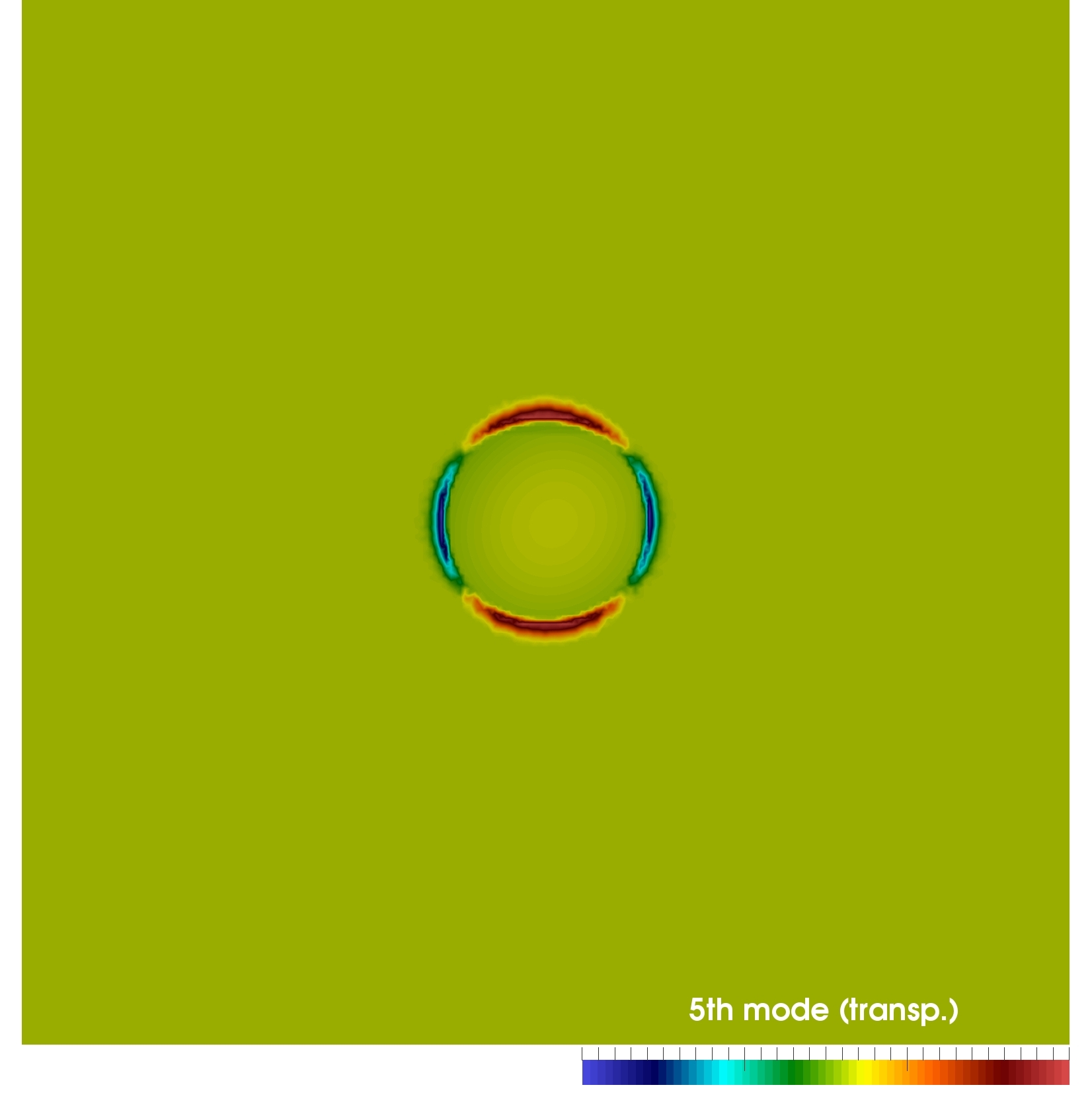}
\includegraphics[width=0.32\textwidth]{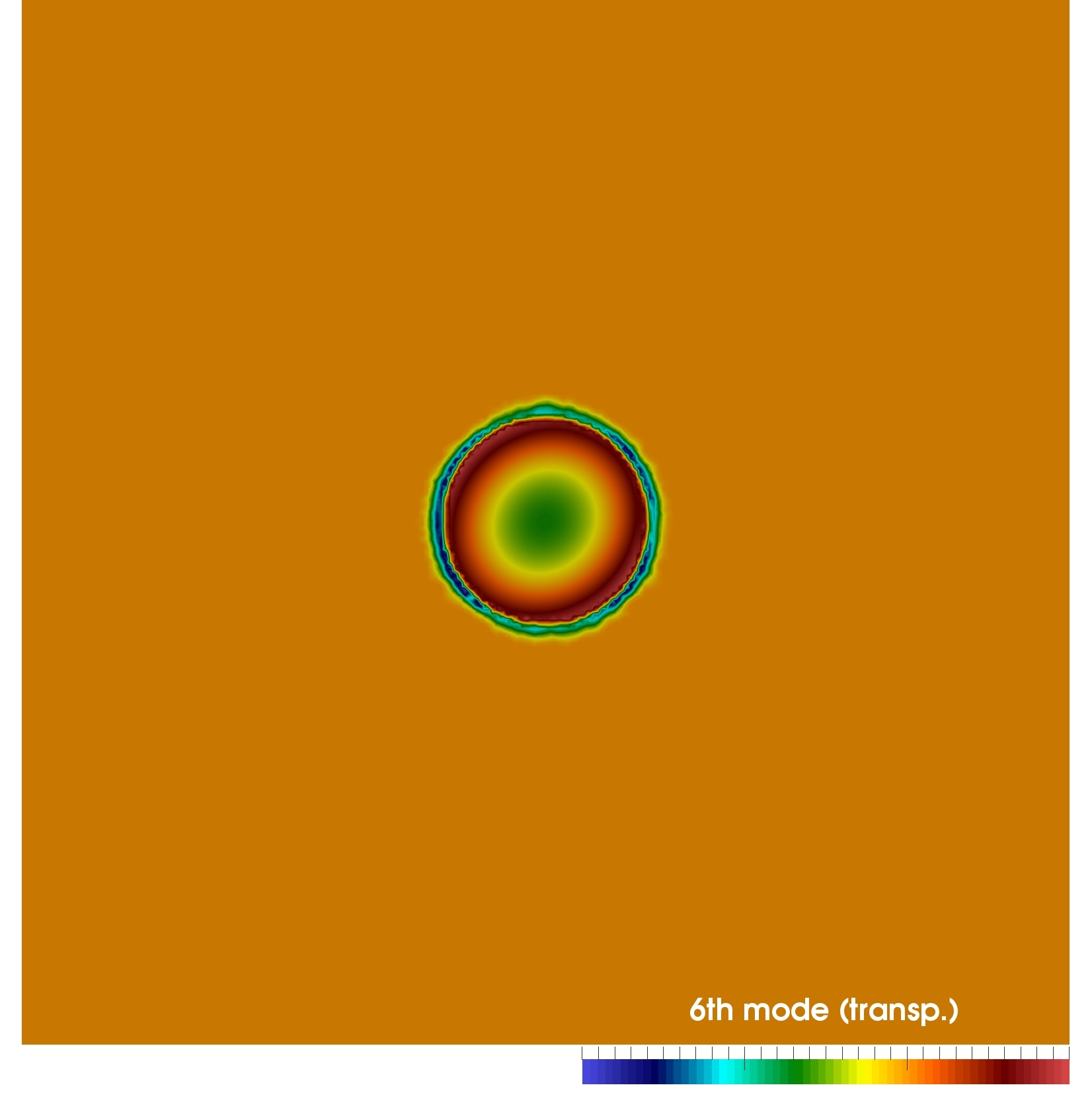}\\
\caption{\emph{Darcy flow test case}: first six POD modes for the natural smooth extension with transportation.}
\label{Fig:Poisson_transport_modes}
\end{figure}

The results of the offline stage, run over a training set of 400 snapshots, are summarized in Figure \ref{fig:Poisson eigs}, where POD eigenvalues (normalized to the maximum eigenvalue) are plotted against the number of modes. The three extension methods introduced in Section \ref{sec:extension} are considered, as well as their combination to the snapshots transportation procedure introduced in Section \ref{sec:transport}. Results show that small reduced basis spaces cannot be obtained relying only on extension procedures; indeed, Figure \ref{fig:Poisson eigs} shows a very slow eigenvalue decay, with negligible differences among the different extension options (zero extension faring slightly worse than the others). Moreover, the training set is not large enough to provide a representative reduced basis, as the smallest eigenvalue is only 4 orders of magnitude less than the largest one.
A faster decay is obtained instead thanks to the transportation. In particular, comparing to the cases with extension only, an improvement of almost three orders of magnitude is obtained for $N=100$, more than four for $N=200$; furthermore, the transport POD drops below numerical precision at $N=300$. Such improvement is instrumental in keeping a low space dimension in the ROM, say $N \approx 100$. Finally, we plot the first six basis functions for the case of natural smooth extension without transportation (Figure \ref{Fig:Poisson_no_transport_modes}) and with transportation (Figure \ref{Fig:Poisson_transport_modes}). From a qualitative point of view, the improvement of the transportation procedure is justified by looking at the extent of the support of the resulting basis functions. Indeed, as expected, basis functions generated through the transportation procedure are compactly supported in a circular neighborhood of the origin, while POD modes obtained without transport are characterized by larger (and varying) support.

\begin{figure} \centering
  \includegraphics[width=0.7\textwidth]{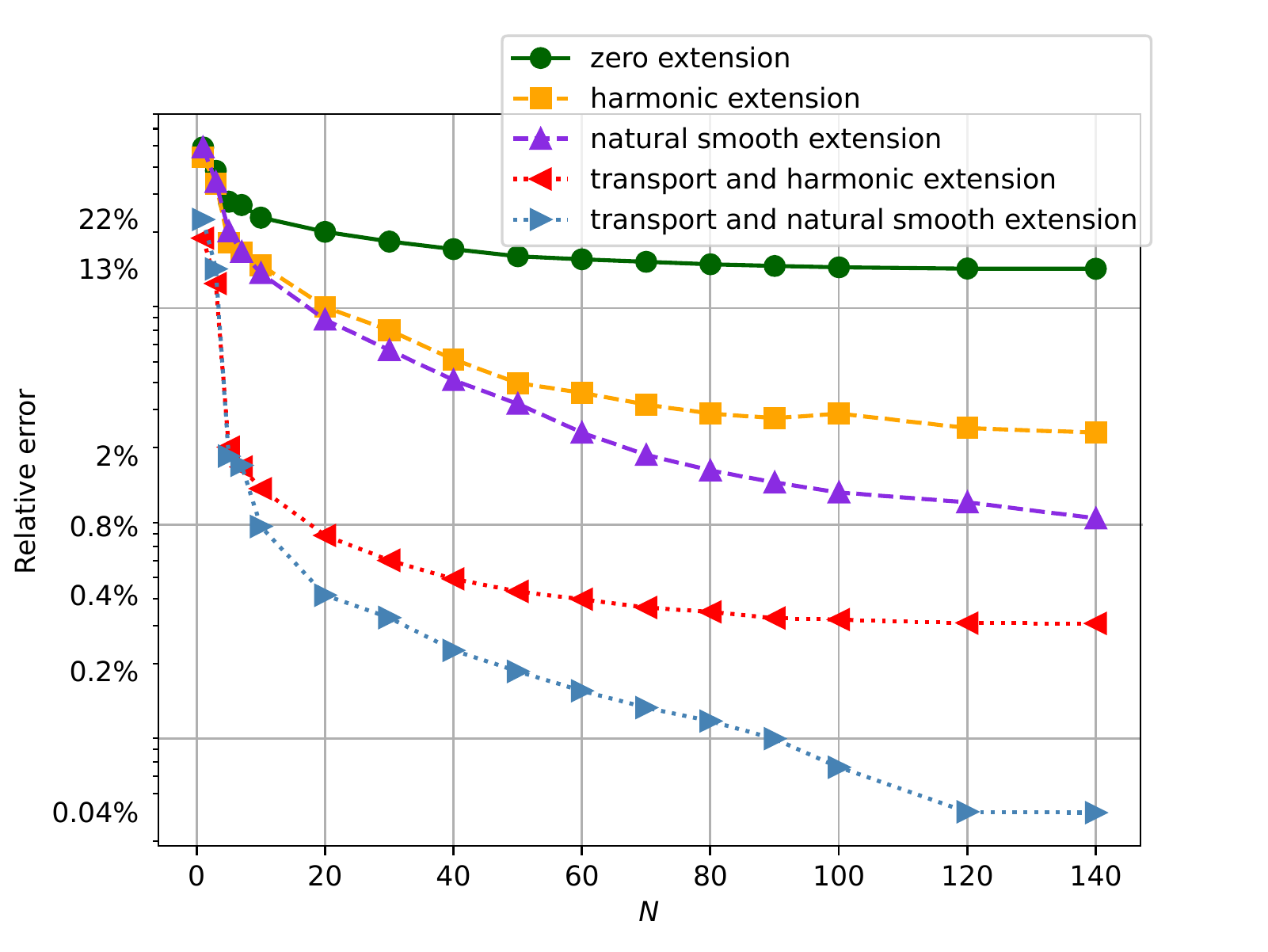}
  \caption{\emph{Darcy flow test case}: error analysis between high fidelity and reduced order approximations.}
  \label{fig:Poisson results}
\end{figure}

\begin{figure}
\centering
\includegraphics[width=0.325\textwidth]{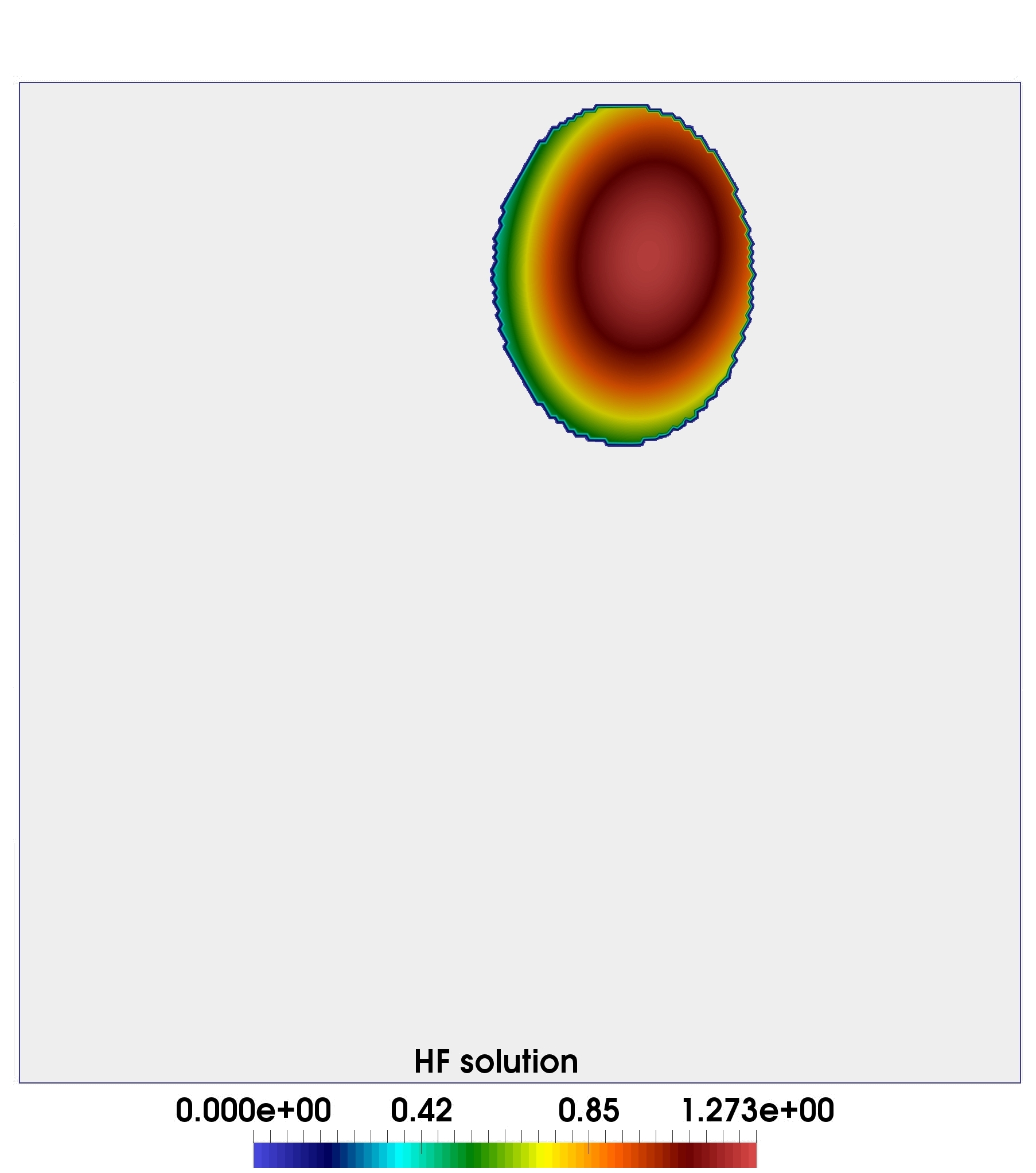}
\includegraphics[width=0.325\textwidth]{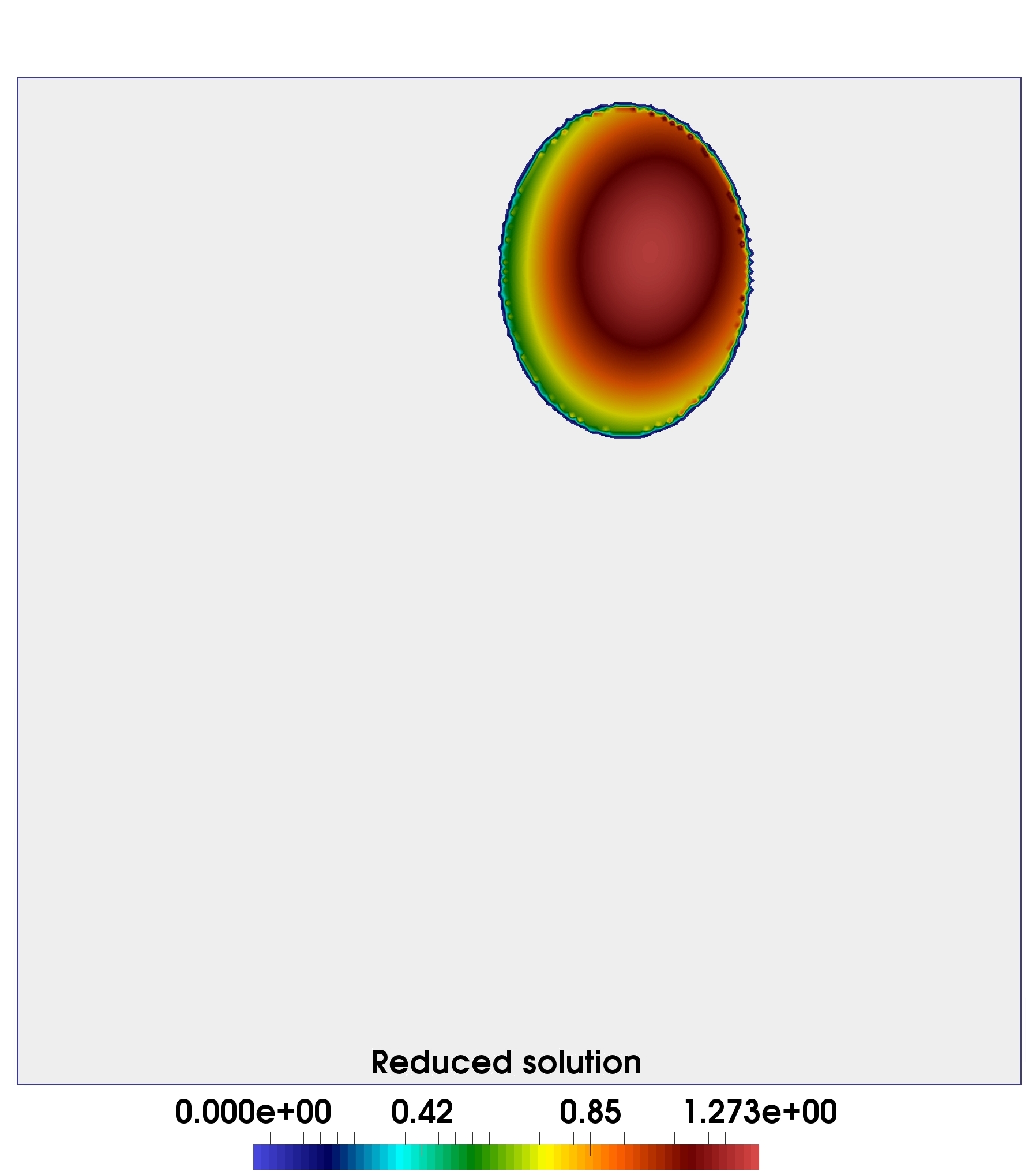}
\includegraphics[width=0.325\textwidth]{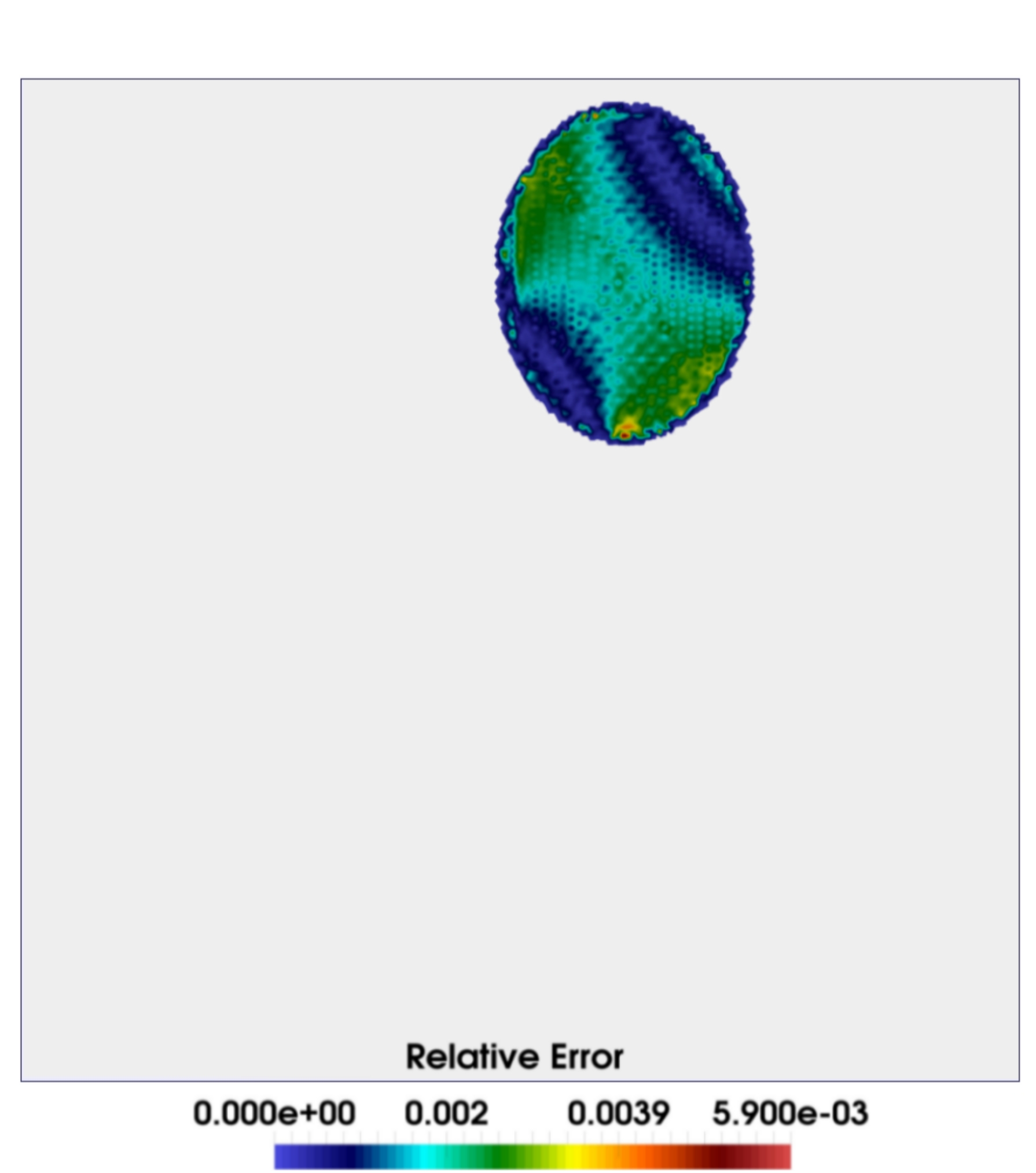}
\caption{\emph{Darcy flow test case}: the high fidelity solution (left), the reduced basis solution (center) for the natural smooth extension with transportation, and the corresponding relative error (right) for the random parameter $\mu = (1.3222,1.7666,0.2514, 0.7365)$.}
\label{Fig:RB-HF sol_Error}
\end{figure}

Figure \ref{fig:Poisson results} shows the results of an error analysis between the reduced order and high fidelity approximations over a testing set of {\blue{$30$}} parameter values, in the $L^2$ norm. In particular, the average of the relative error over the testing set is plotted against the reduced basis size. The reduced solution obtained from the zero extension is inaccurate even for $N=140$, being affected by relative errors of the order of $10^{-1}$. A non-zero extension is beneficial, resulting in relative errors that are of the order of $10^{-2}$ for the maximum value of $N$. Furthermore, the combination with inverse transportation allows to further improve results, up to errors of $10^{-4}$ for $N=140$ in the case of POD basis obtained from transport and natural smooth extension. Thus, the pivotal role of snapshots transportation can be inferred from these results, being capable of improving the results of almost three orders of magnitude compared to the simplest zero extension.
Nonetheless, all methods reach a plateau after which no further improvements are shown. We claim that this is due to integration errors occurring on $\partial\mathcal{D}(\mu)$ and Nitsche weak imposition of Dirichlet boundary conditions, as the maximum values of the error are consistently attained on parts of the the boundary (see Figure \ref{Fig:RB-HF sol_Error} for a representative case). 
\begin{figure} \centering
  \includegraphics[width=0.5\textwidth]{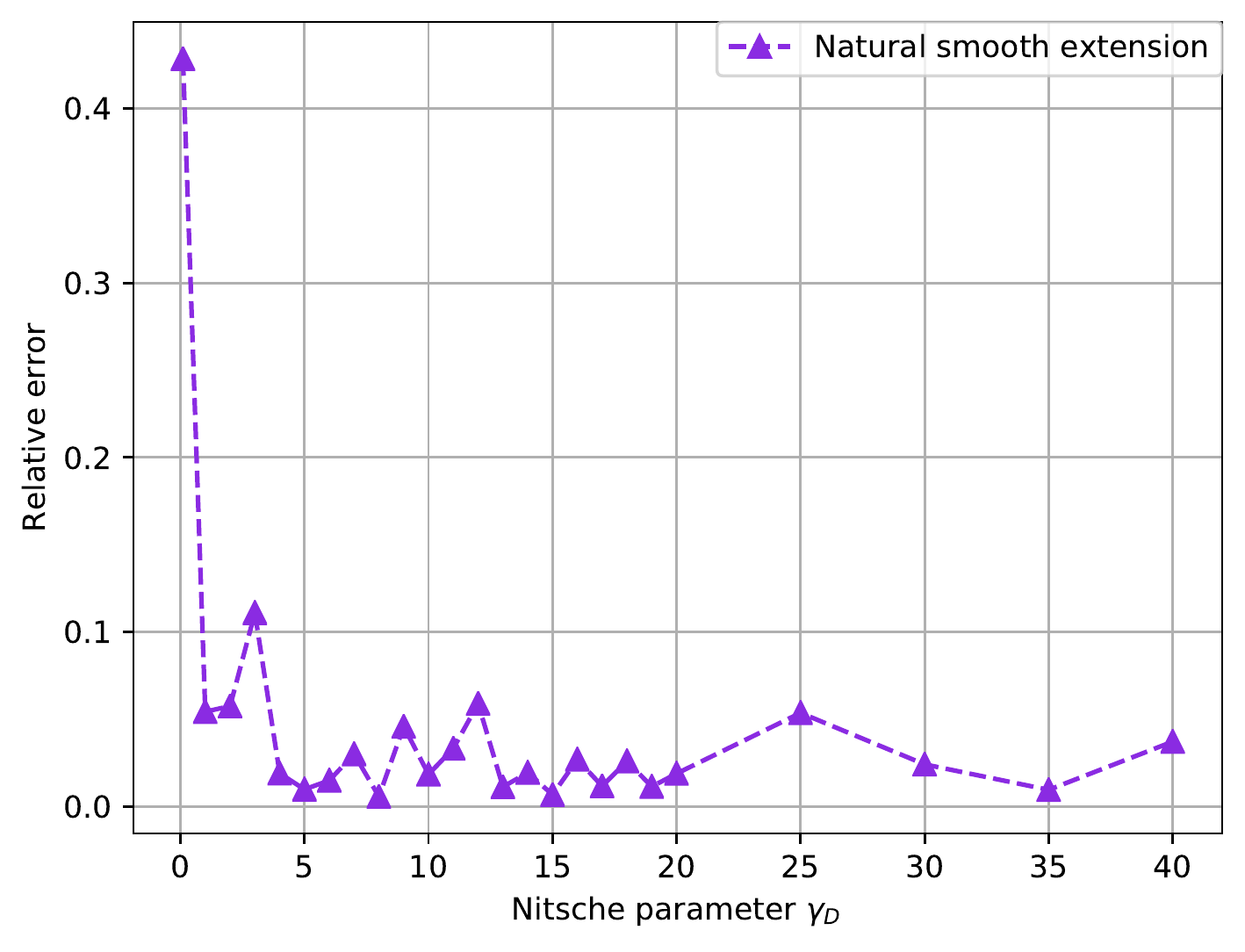}
  \caption{{\blue{\emph{Darcy flow test case}: dependence of the relative error (over a random testing set) on the value of the Nitsche parameter in the range $[0.1, 40]$.}}}
  \label{fig:gamma_D_dependence}
\end{figure}
{\blue{Further investigation on this issue shows slightly improved relative errors for some specific Nitsche parameters values. In particular, we tested our algorithm for several values of $\gamma_D$ in the range $[0.1, 40]$. For each coefficient $\gamma_D$ we trained the ROM using $400$ snapshots and we queried the reduced solver for a fixed number of 120 basis components and 30 random geometrical parameter test cases. We employed the best extension, namely the natural smooth extension, without transportation (in order to separate the beneficial effects of transportation from those of changing $\gamma_D$). As seen in \autoref{fig:gamma_D_dependence}, the average relative error shows a mild dependence on the chosen value of $\gamma_D$, at least in the range $[5, 40]$. Yet, slightly improved results (compared to those of \autoref{fig:Poisson results}) can be obtained for the choice of $\gamma_D=8$, resulting in a relative error of $0.0056$, instead of the relative error $0.0112$ for $\gamma_D=10$ as in \autoref{fig:Poisson results}. Even though a careful tuning of $\gamma_D$ allows to improve the accuracy, such issue should be further investigated in future publications.
}}  

\subsection{Steady Stokes problem}\label{sec:num_exp_vector}
This experiment considers the Stokes flow around a circular cylinder, embedded in a cavity ${\mathcal{B}}= [-2, 2] \times [-1, 1]$. The parametrization, which allows to change the $y$ coordinate of the center of the cylinder, is inspired by the numerical example in \cite[Section 4]{KaratzasStabileNouveauScovazziRozza2018}.
The embedded domain is thus parametrized through $\mu\in\mathcal{K}=[-0.5,0.5]$ according to the following level set expression:
$$
\phi(x,y;\mu) = \left(x-\frac{3}{2}\right)^2 +\left(y-\mu\right)^2 - R^2,
$$
where the parameter $\mu$ describes the $y$ coordinate of the center of the embedded circular domain. The radius $R=0.2$, while the viscosity $\nu$ is set to $1$. A constant velocity in the $x$ direction, $u_{\text{in}} = 1$ is applied (strongly) at the left side of the domain, and a Neumann boundary condition is applied on the right. No slip, i.e. homogeneous Dirichlet boundary conditions are applied (strongly) on the top and bottom edges, as well as (through Nitsche weak imposition) on the boundary of the embedded cylinder.
The results for the test problem have been obtained with a mesh size of $h = 0.0350$ for the background mesh ${{\mathcal{B}}_h}$, 
using 15022 triangles for the discretization and $P^1/P^1$ finite elements in space with stabilization terms as described in Section \ref{sec:HFStokes}.

In this test case we will only employ the natural smooth extension among all the possible extension options presented in Section \ref{sec:extension}, as it was the one which gave best results in Section \ref{sec:num_exp_scalar}.
For what concerns the transportation method of Section \ref{sec:transport}, we choose the reference parameter $\overline\mu=0$, which corresponds to the center of the parameter range $\mathcal{K}$. The expression $\tilde{y} = \boldsymbol{\tau}(y;\mu)=y+\mu(1-y^2)$ is chosen so that the {\blue{center-line}} $y = 0$ of the domain is mapped to $\tilde{y} = \mu$, while keeping top and bottom walls fixed. Parameter range $\mathcal{K}$ has been chosen so that the inverse mapping being $\boldsymbol{\tau}^{-1}(y;\mu)=\frac{1}{2\mu}(1-(4\mu ^2 - 4\mu y + 1)^{\frac{1}{2}})$ is well defined.

\begin{figure} \centering
  \includegraphics[width=0.7\textwidth]{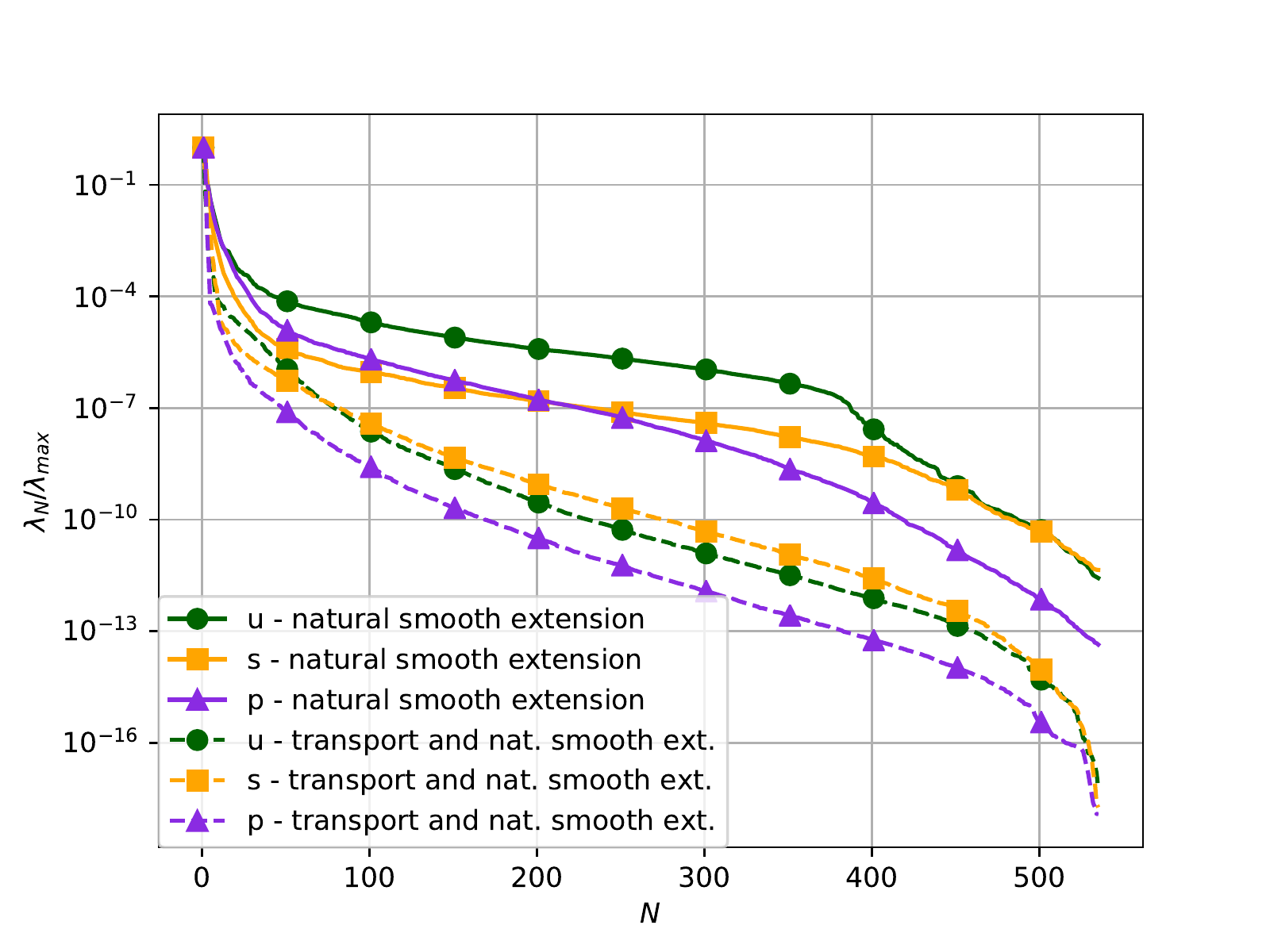}
  \caption{\emph{Stokes flow test case}: The POD eigenvalues decay (normalized to the maximum eigenvalue) for a training set of 600 snapshots is visualized against the number of modes.}
  \label{fig:Stokes_eigs_Cavity}
\end{figure}
\begin{figure} \centering
\includegraphics[width=0.47\textwidth]{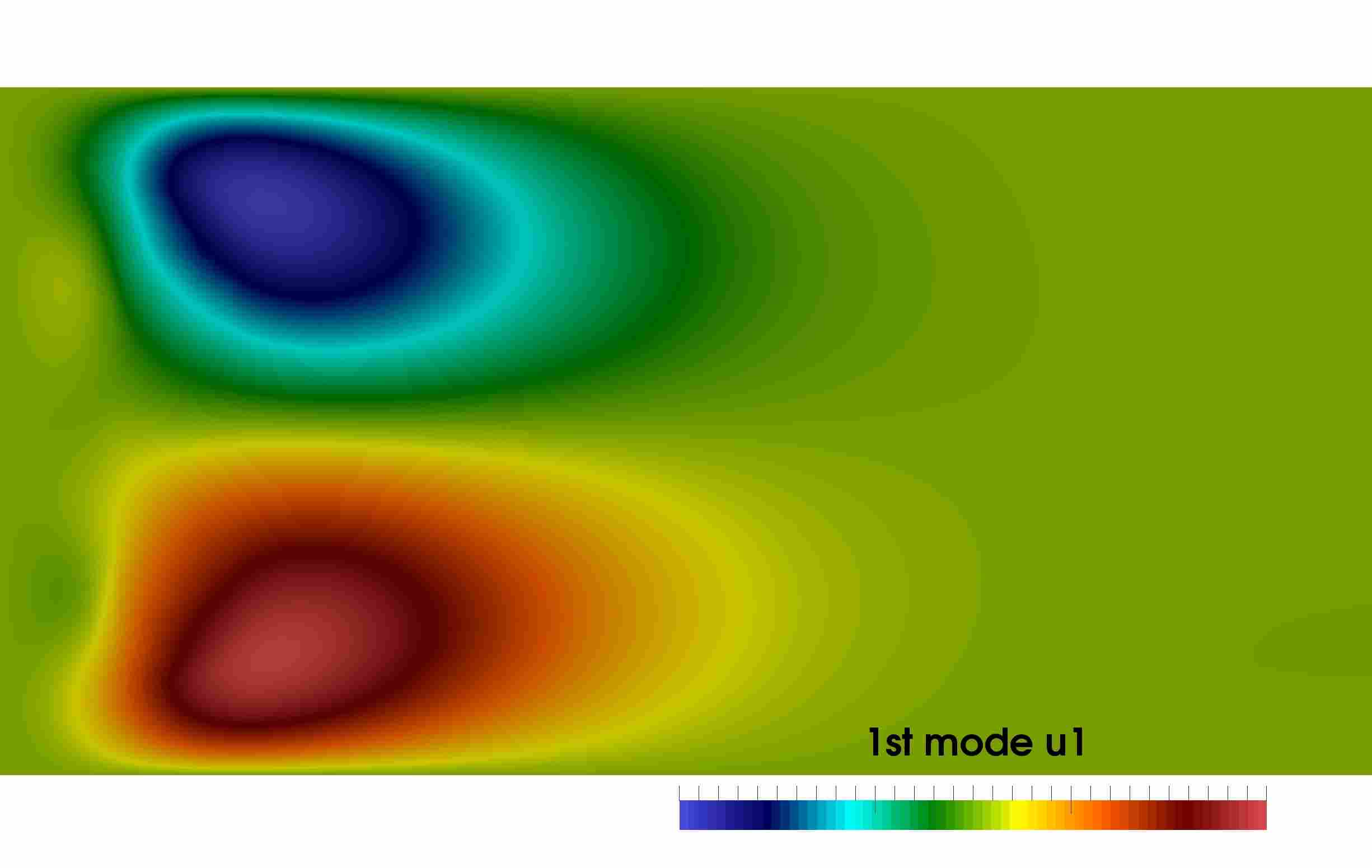}
\includegraphics[width=0.47\textwidth]{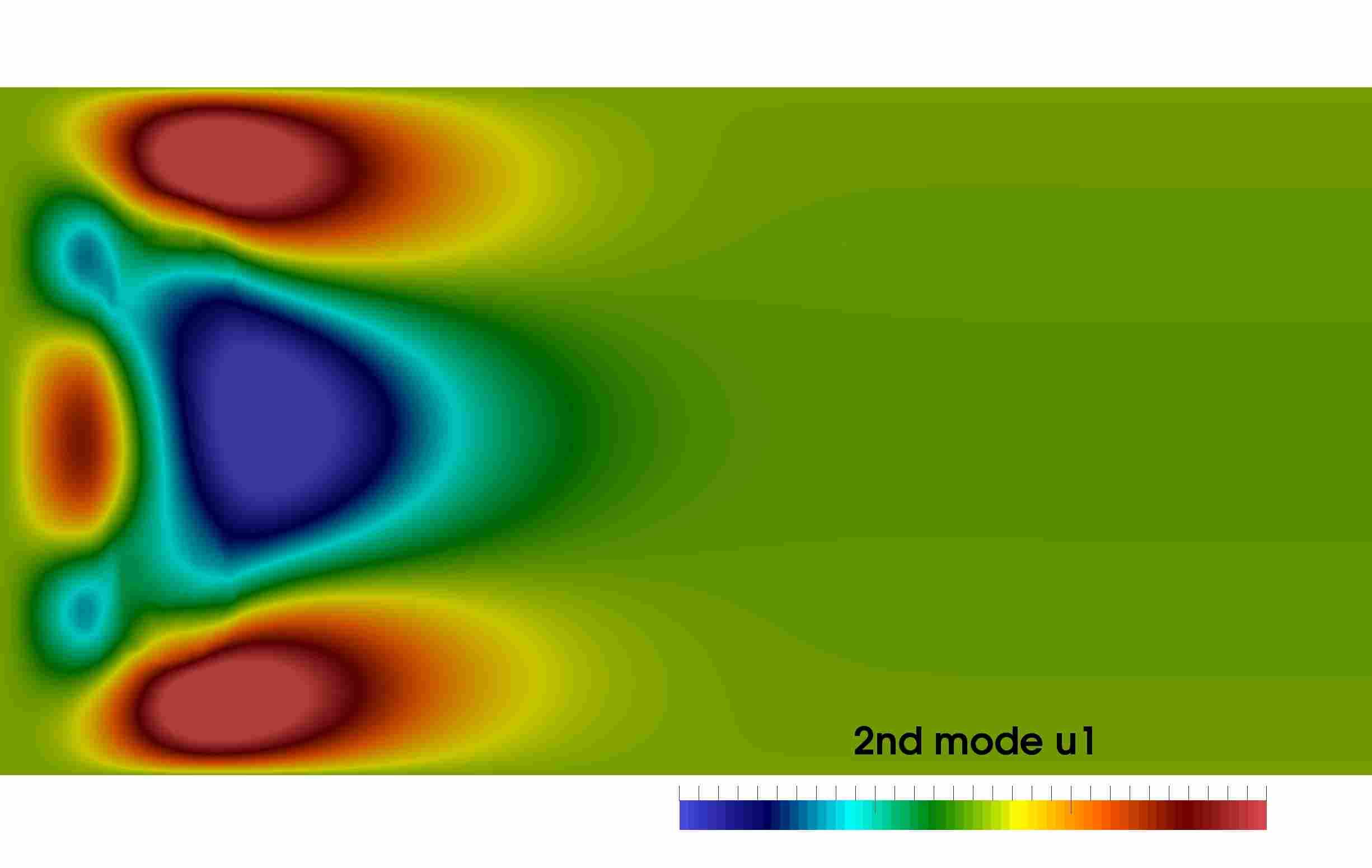}\\
\includegraphics[width=0.47\textwidth]{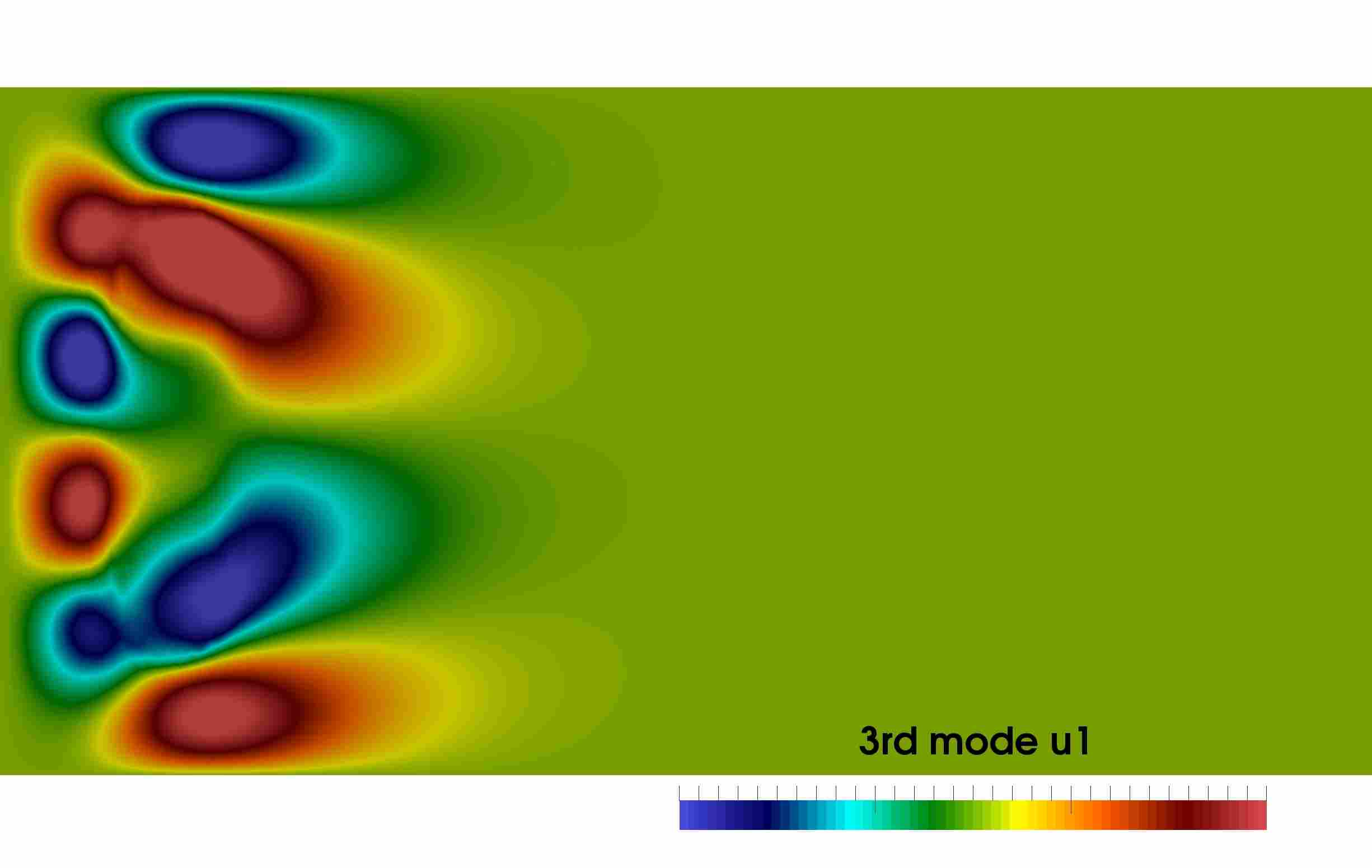}
\includegraphics[width=0.47\textwidth]{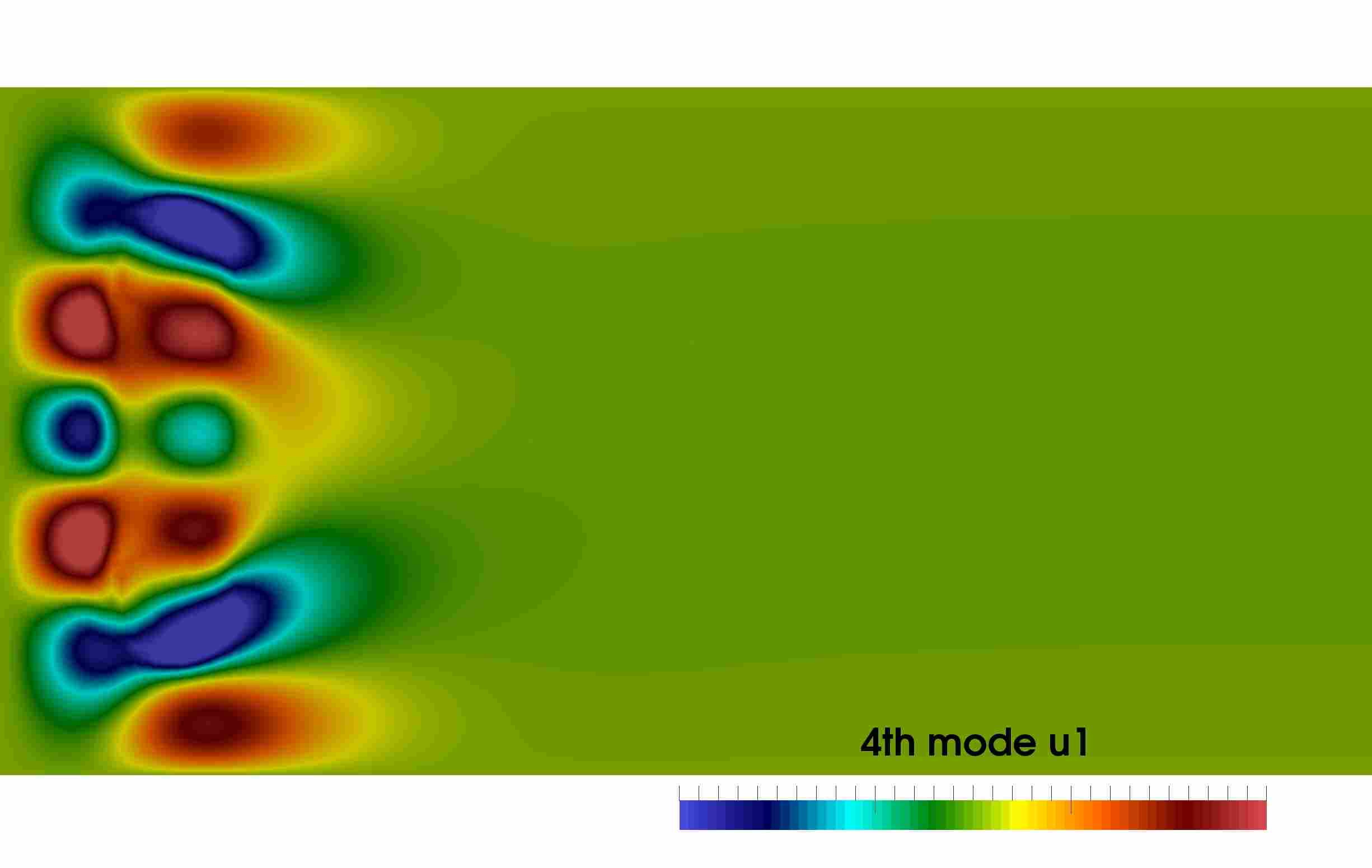}\\
\includegraphics[width=0.47\textwidth]{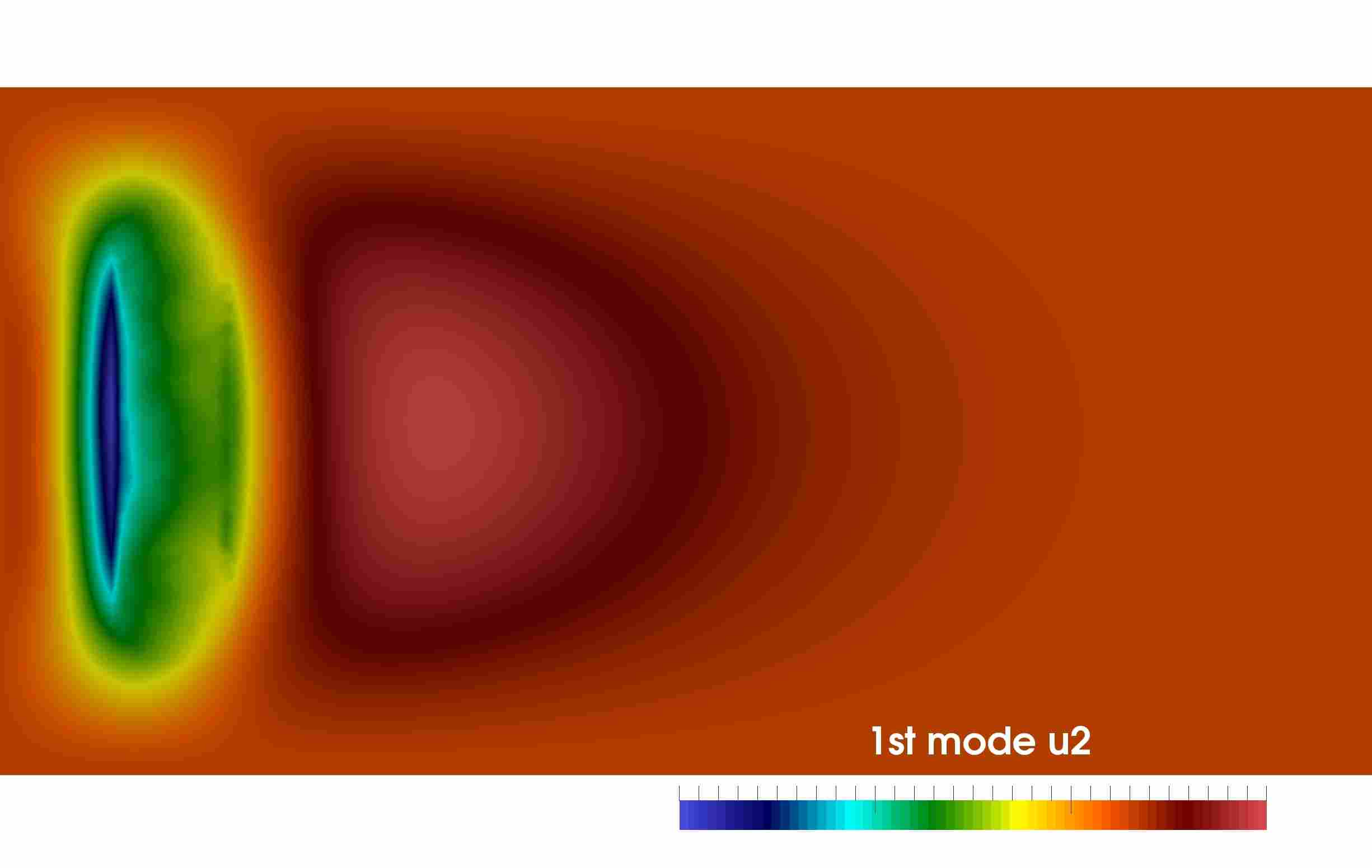}
\includegraphics[width=0.47\textwidth]{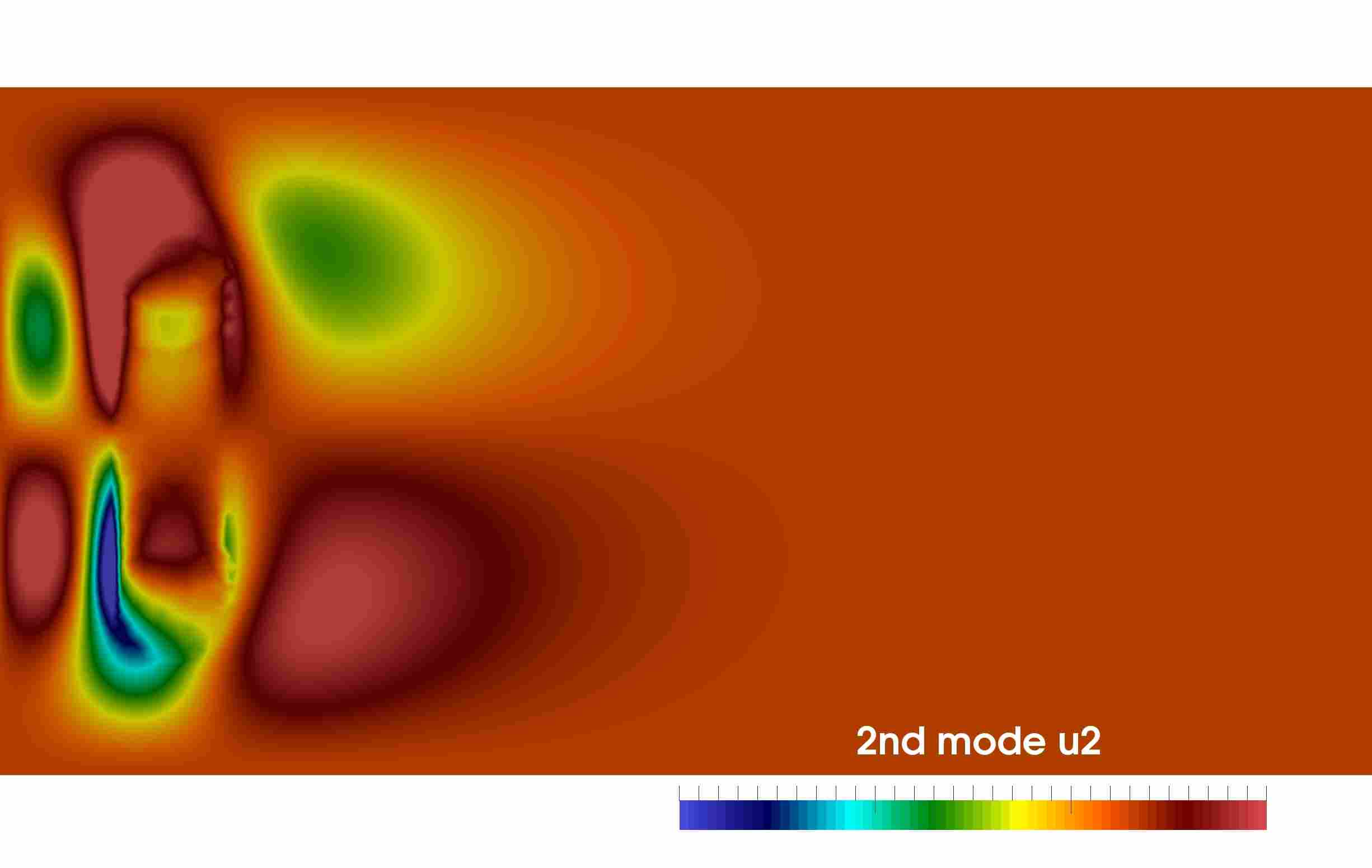}\\
\includegraphics[width=0.47\textwidth]{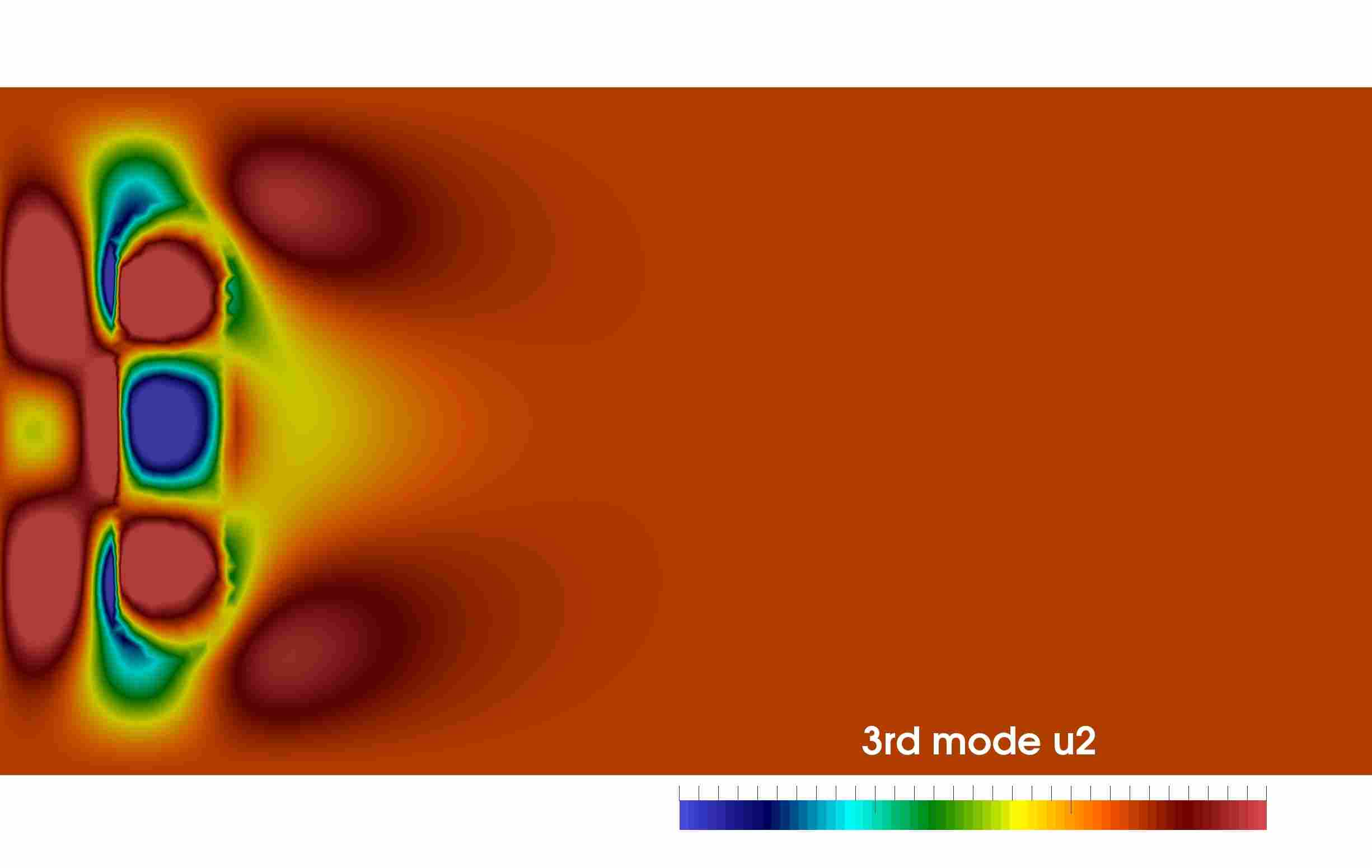}
\includegraphics[width=0.47\textwidth]{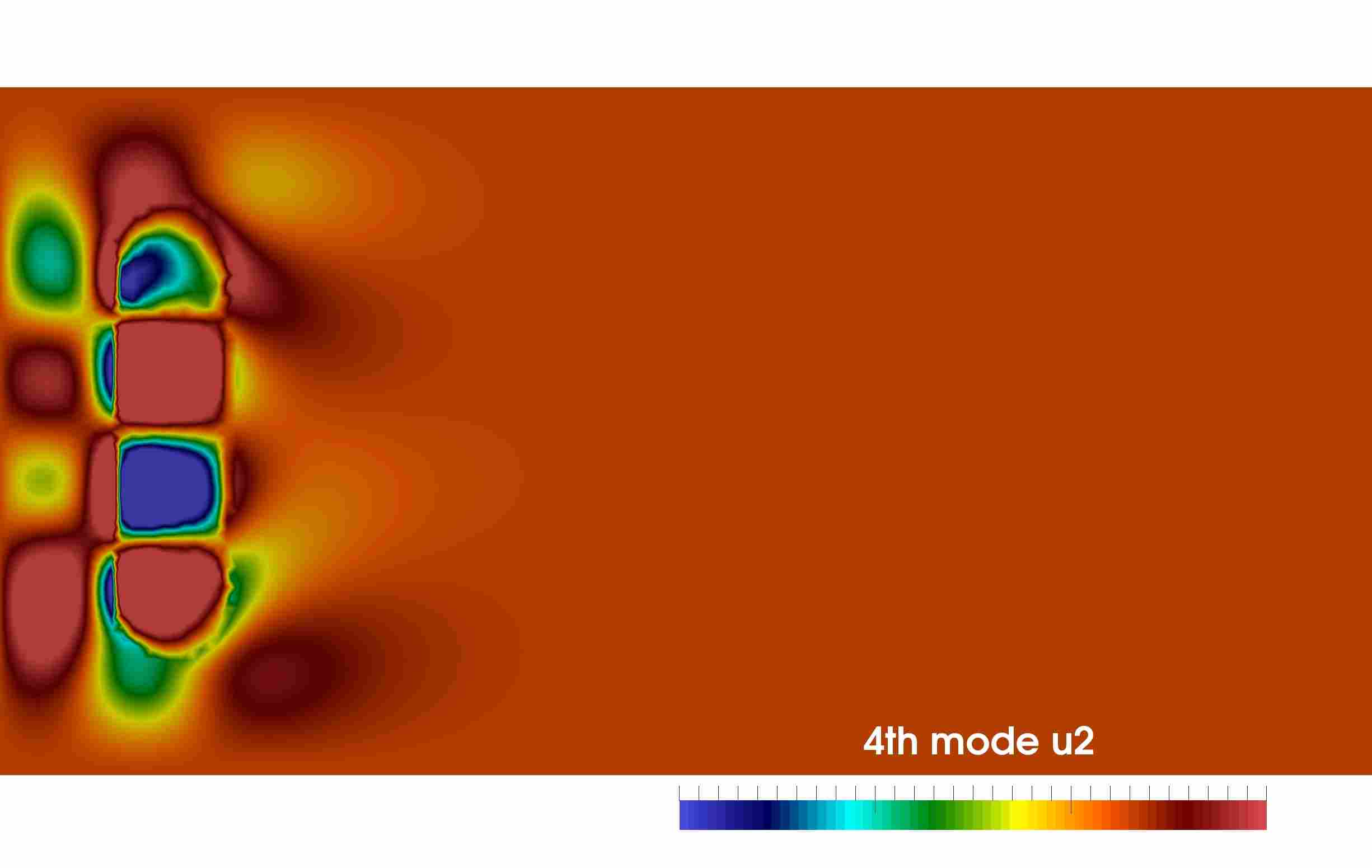}
\caption{\emph{Stokes flow test case}: first four velocity POD modes (plotted component by component) for the natural smooth extension without transportation.}
\label{Fig:Stokes_u_no_transport_modes}
\end{figure}
\begin{figure}
\centering
\includegraphics[width=0.47\textwidth]{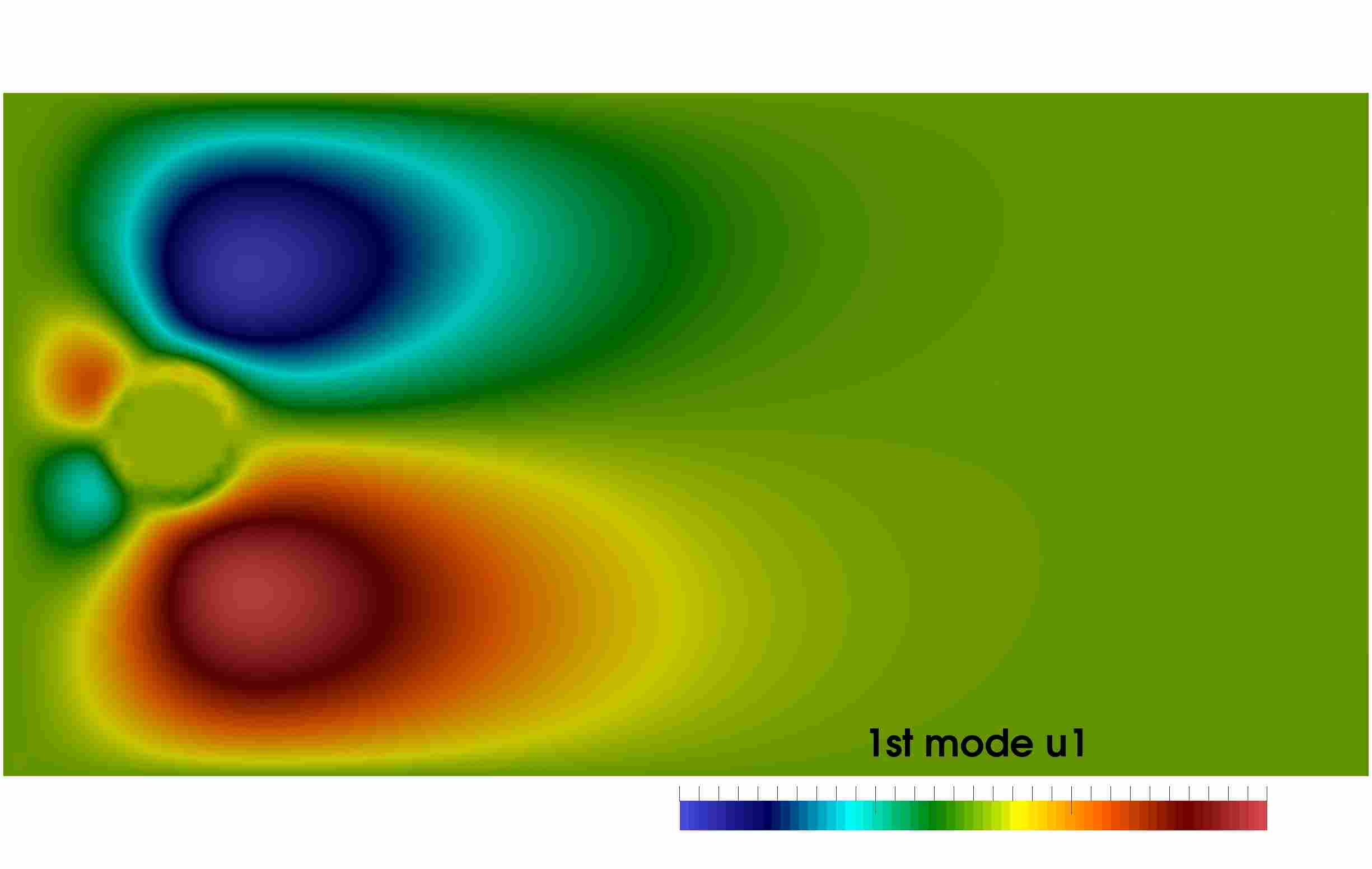}
\includegraphics[width=0.47\textwidth]{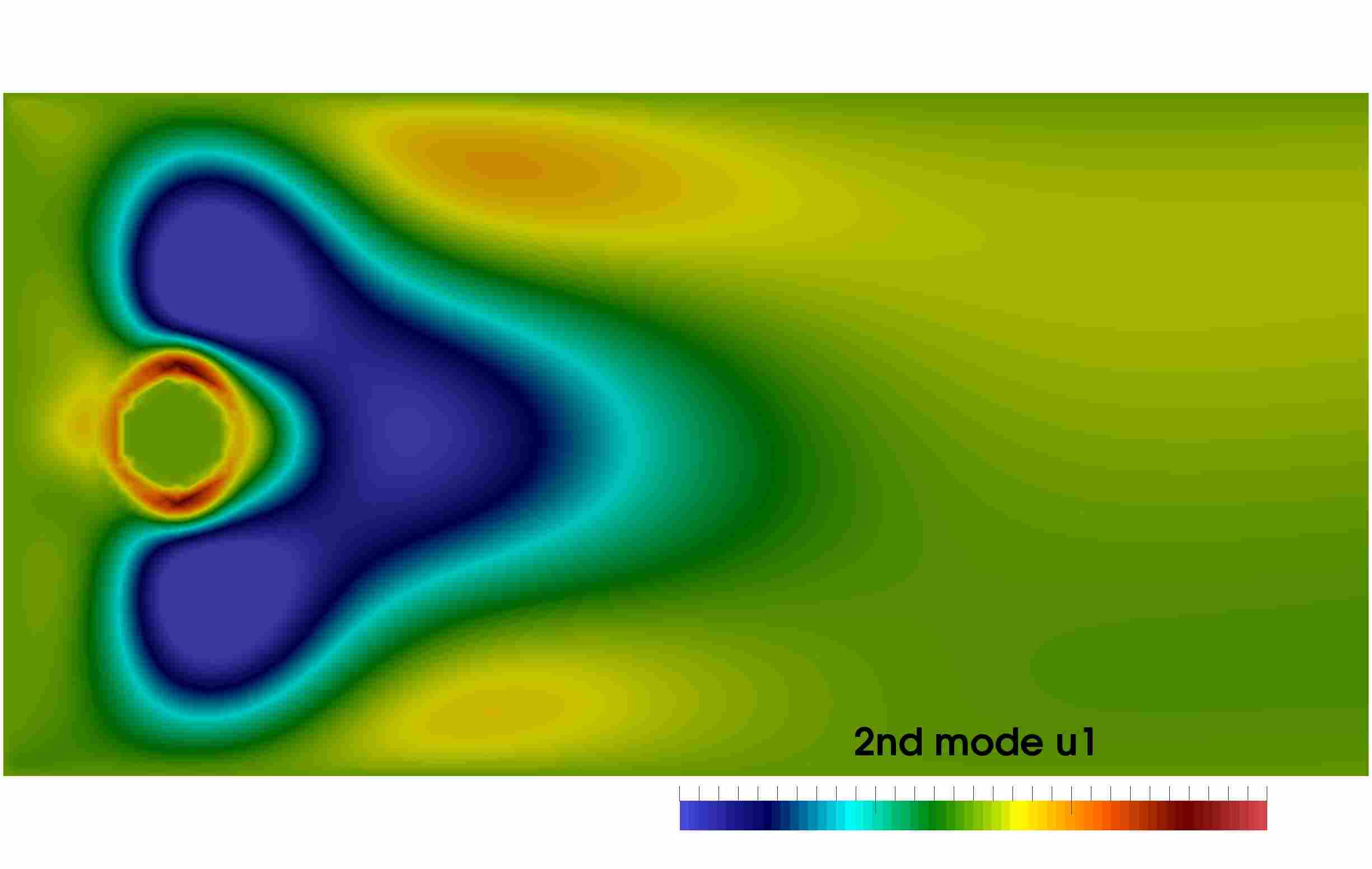}\\
\includegraphics[width=0.47\textwidth]{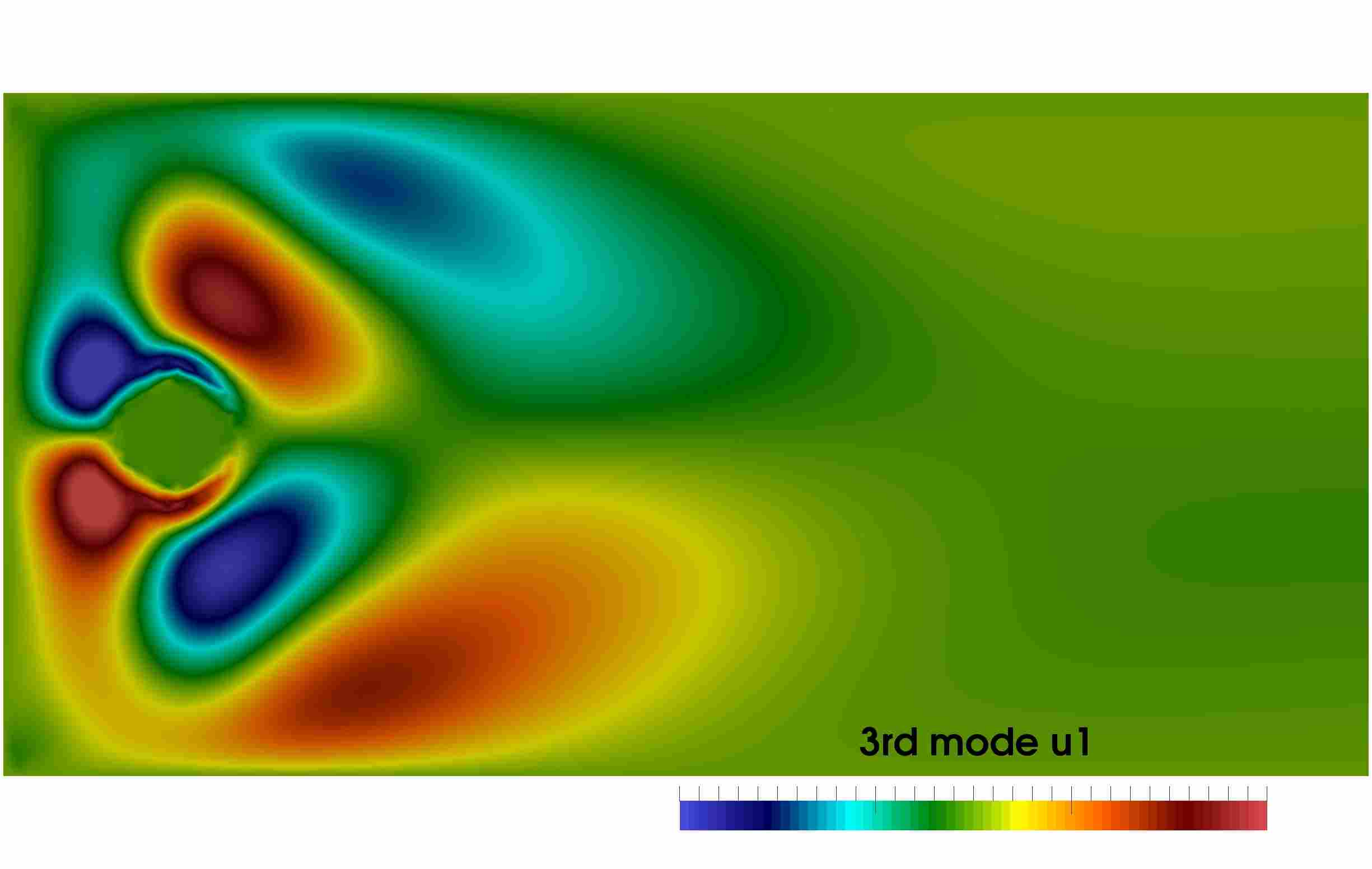}
\includegraphics[width=0.47\textwidth]{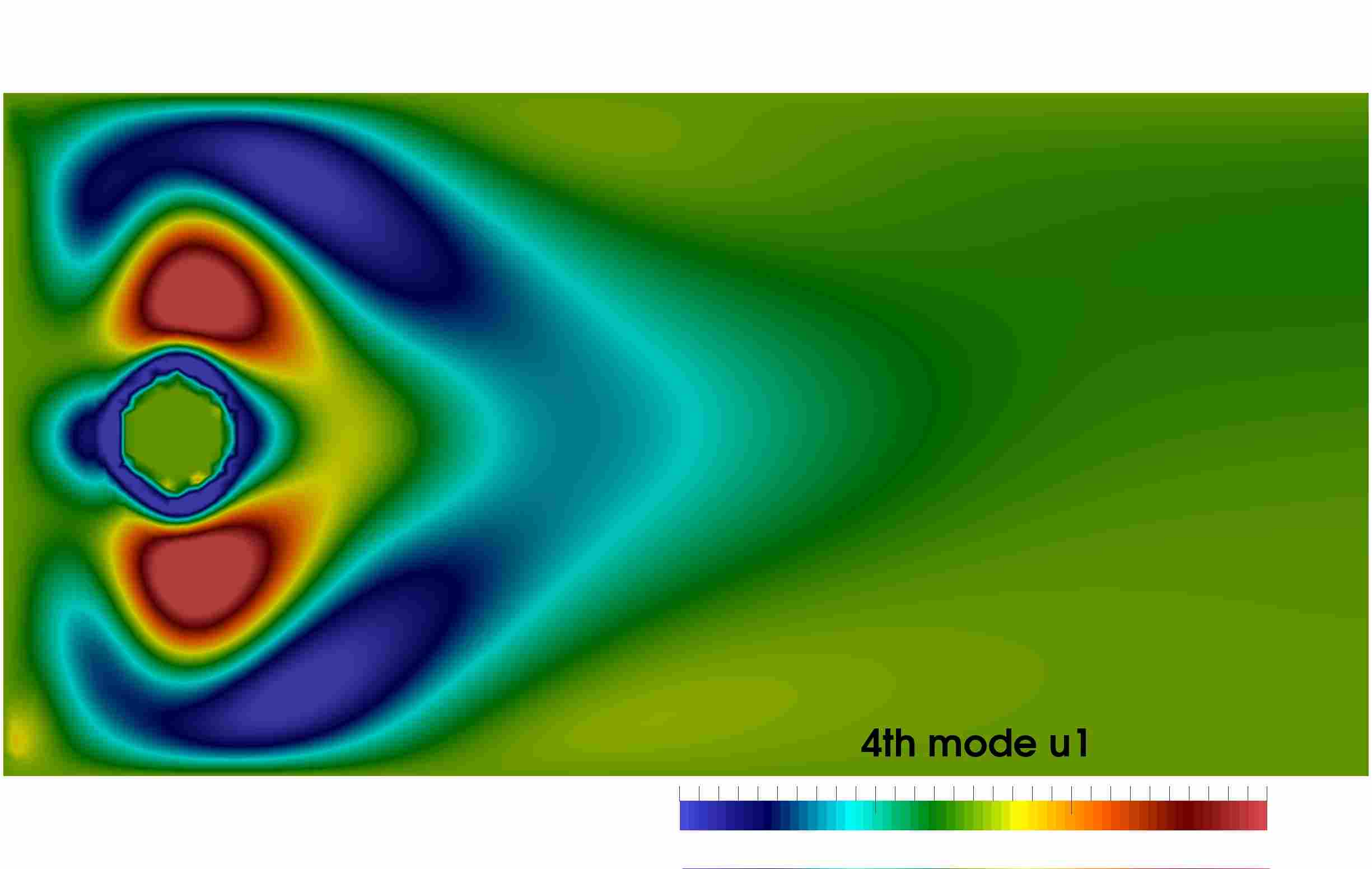}\\
\includegraphics[width=0.47\textwidth]{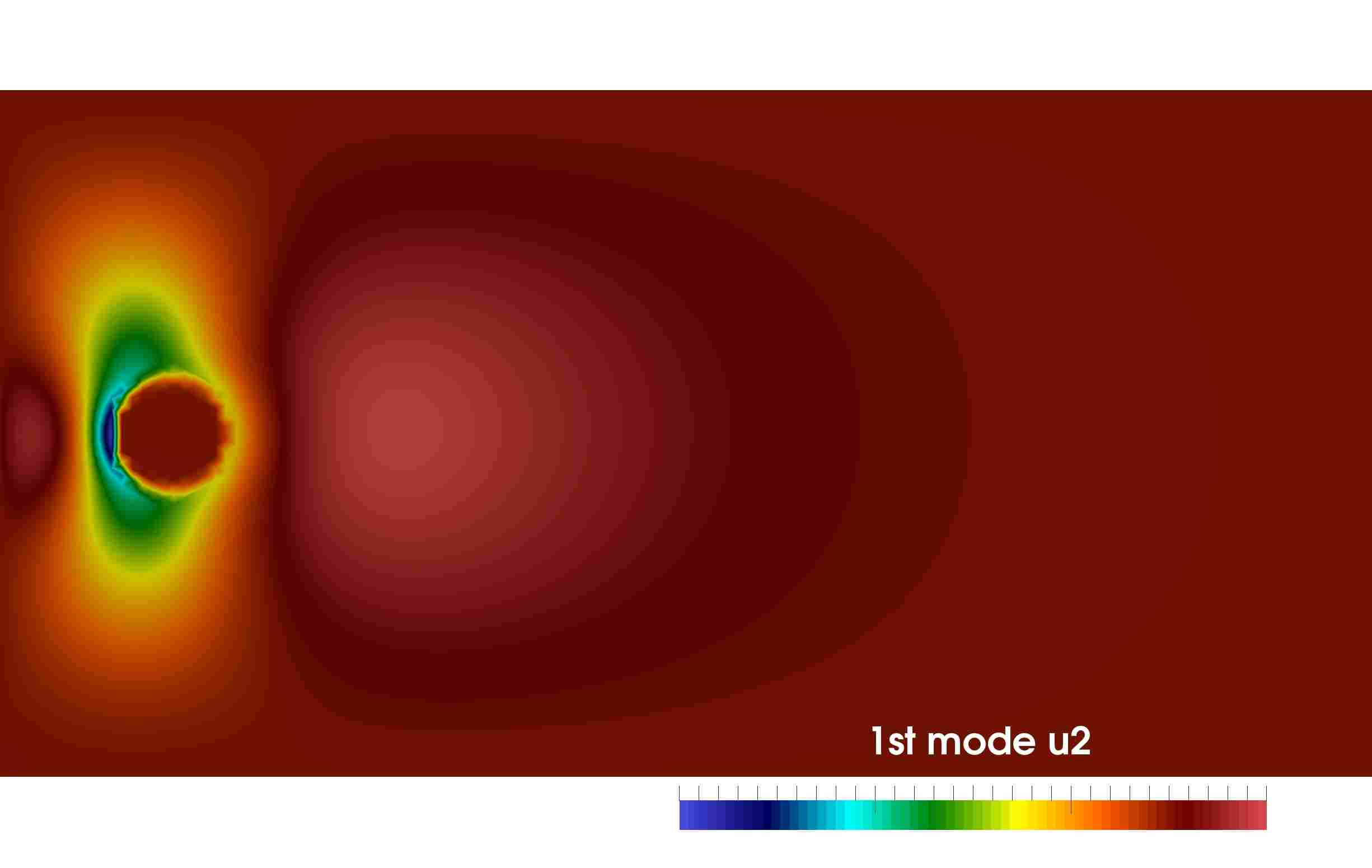}
\includegraphics[width=0.47\textwidth]{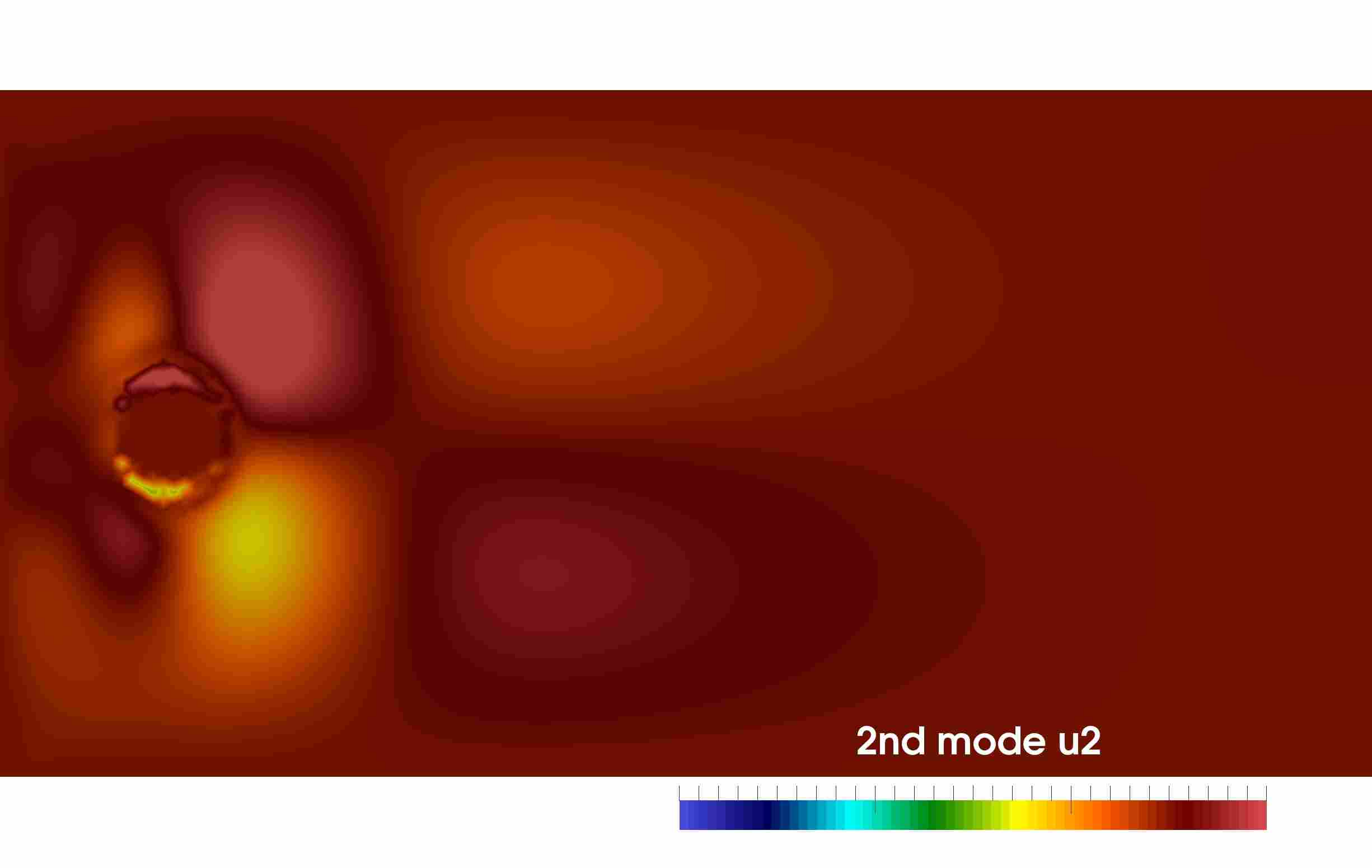}\\
\includegraphics[width=0.47\textwidth]{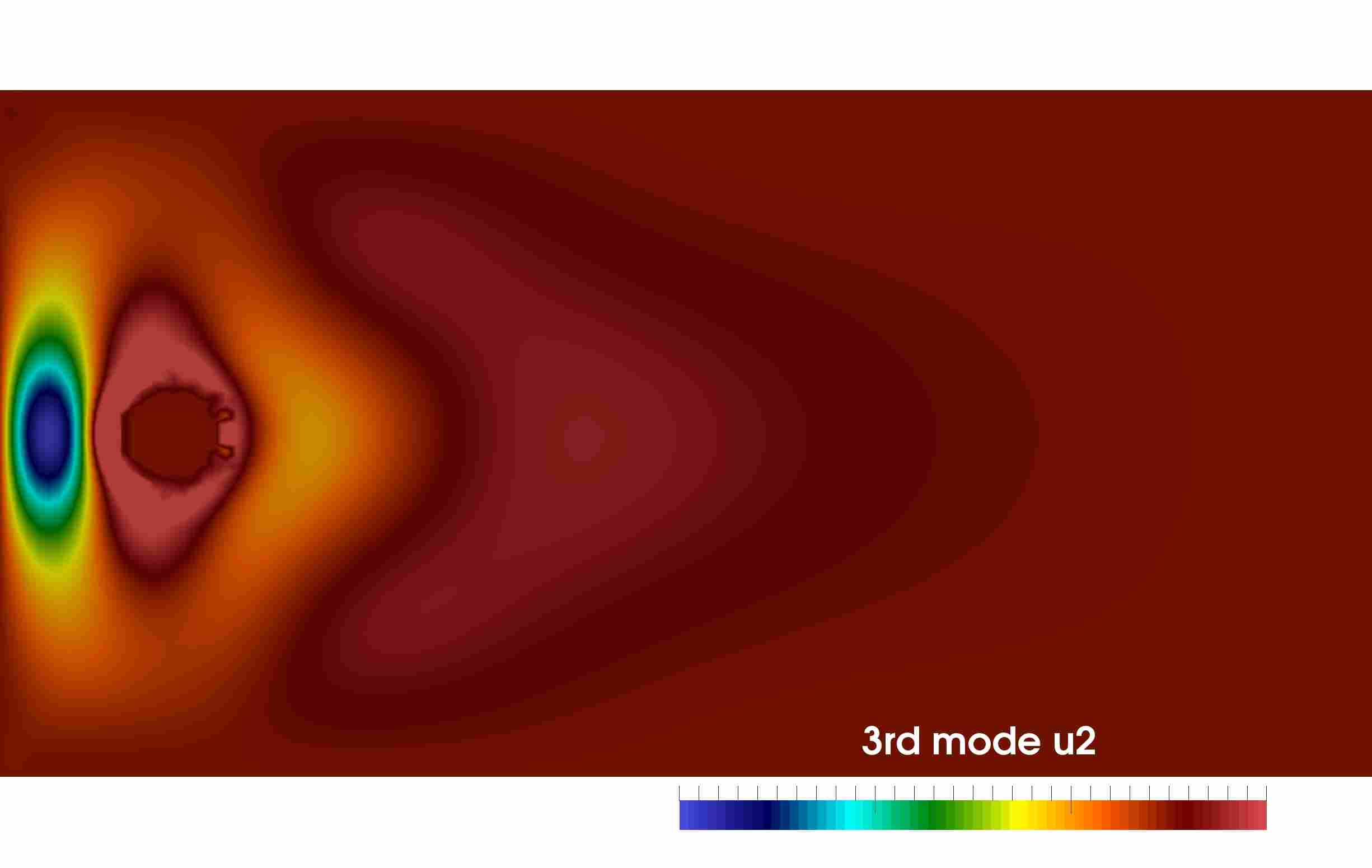}
\includegraphics[width=0.47\textwidth]{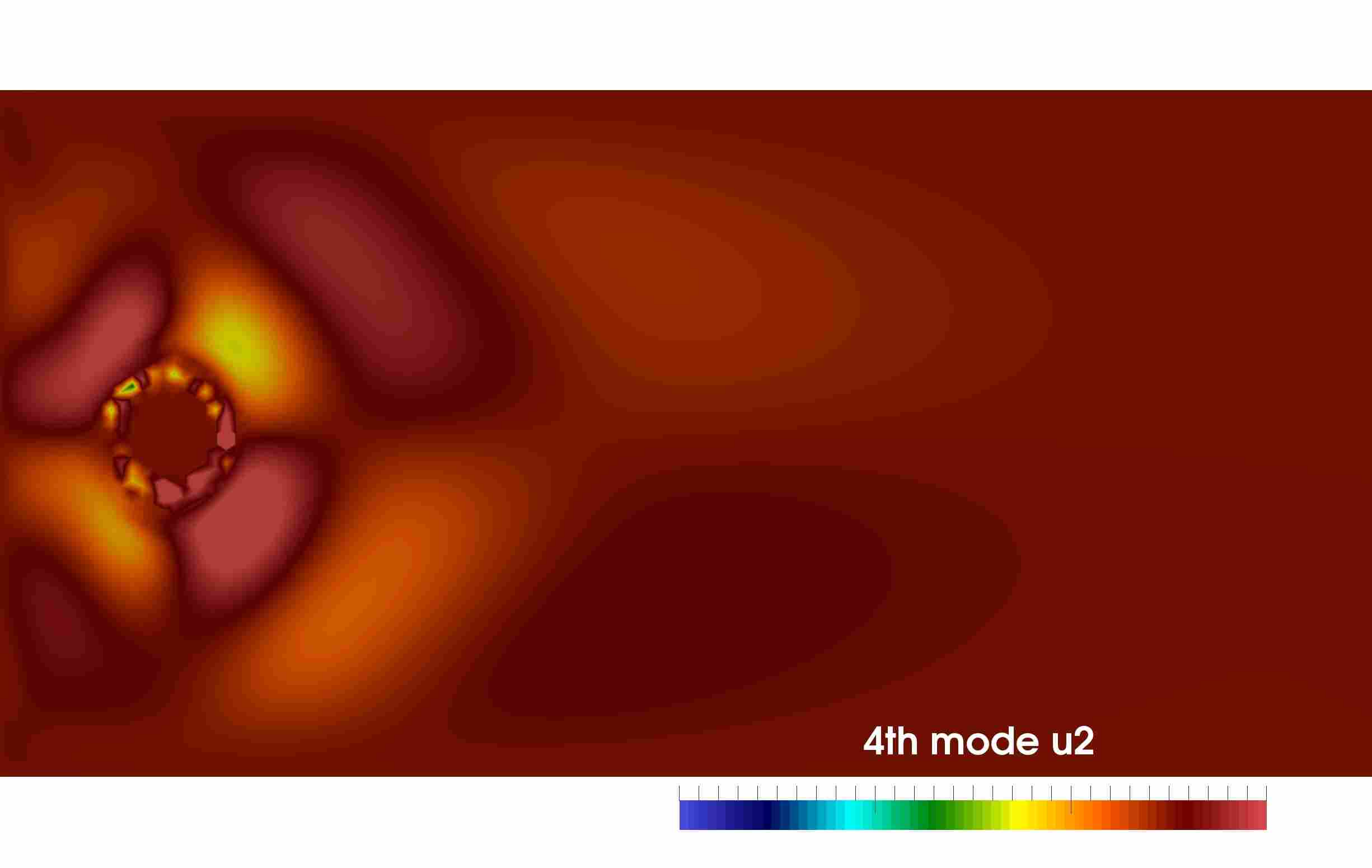}
\caption{\emph{Stokes flow test case}: first four velocity POD modes (plotted component by component) for the natural smooth extension with transportation.}
\label{Fig:Stokes_u_transport_modes}
\end{figure}
\begin{figure}
\centering
\includegraphics[width=0.47\textwidth]{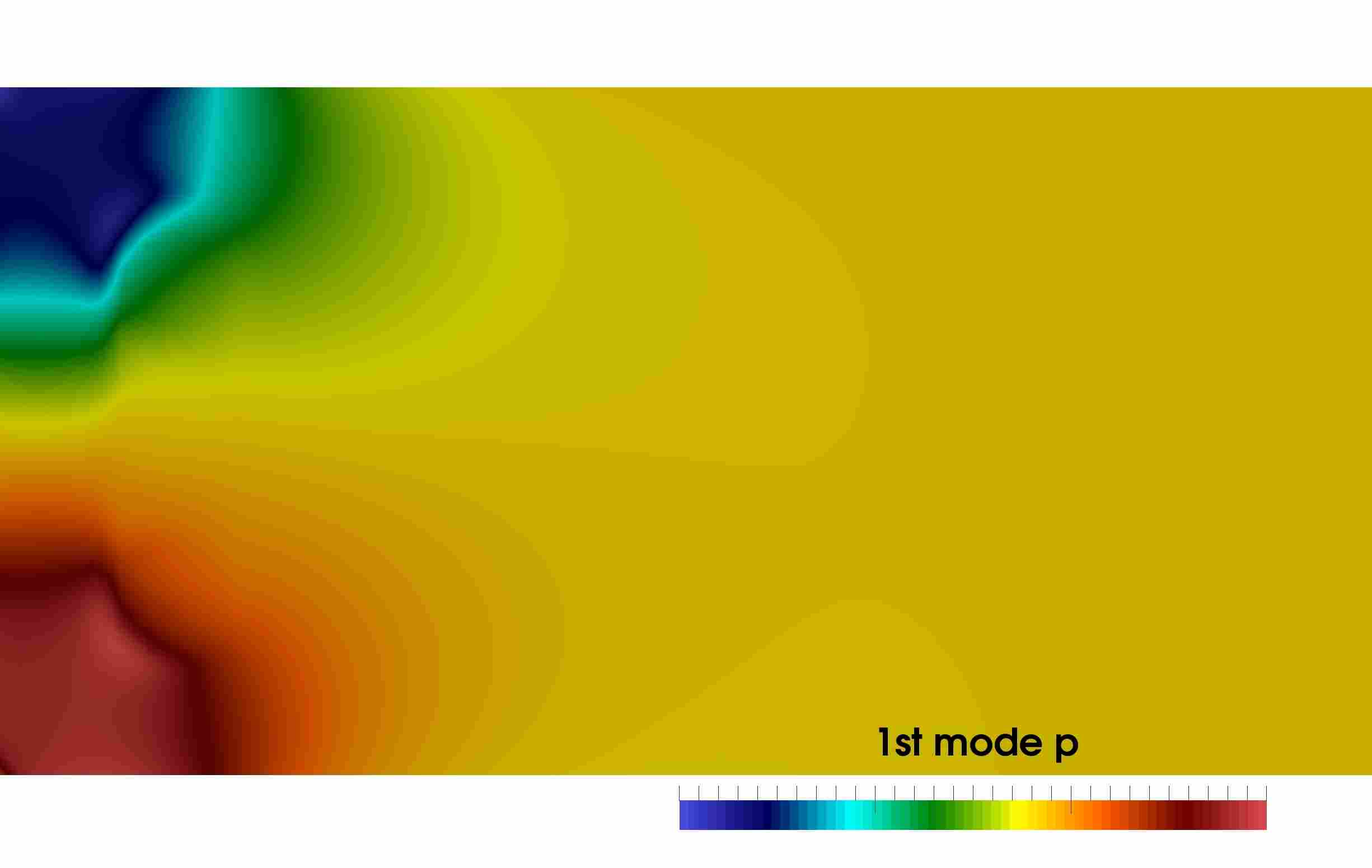}
\includegraphics[width=0.47\textwidth]{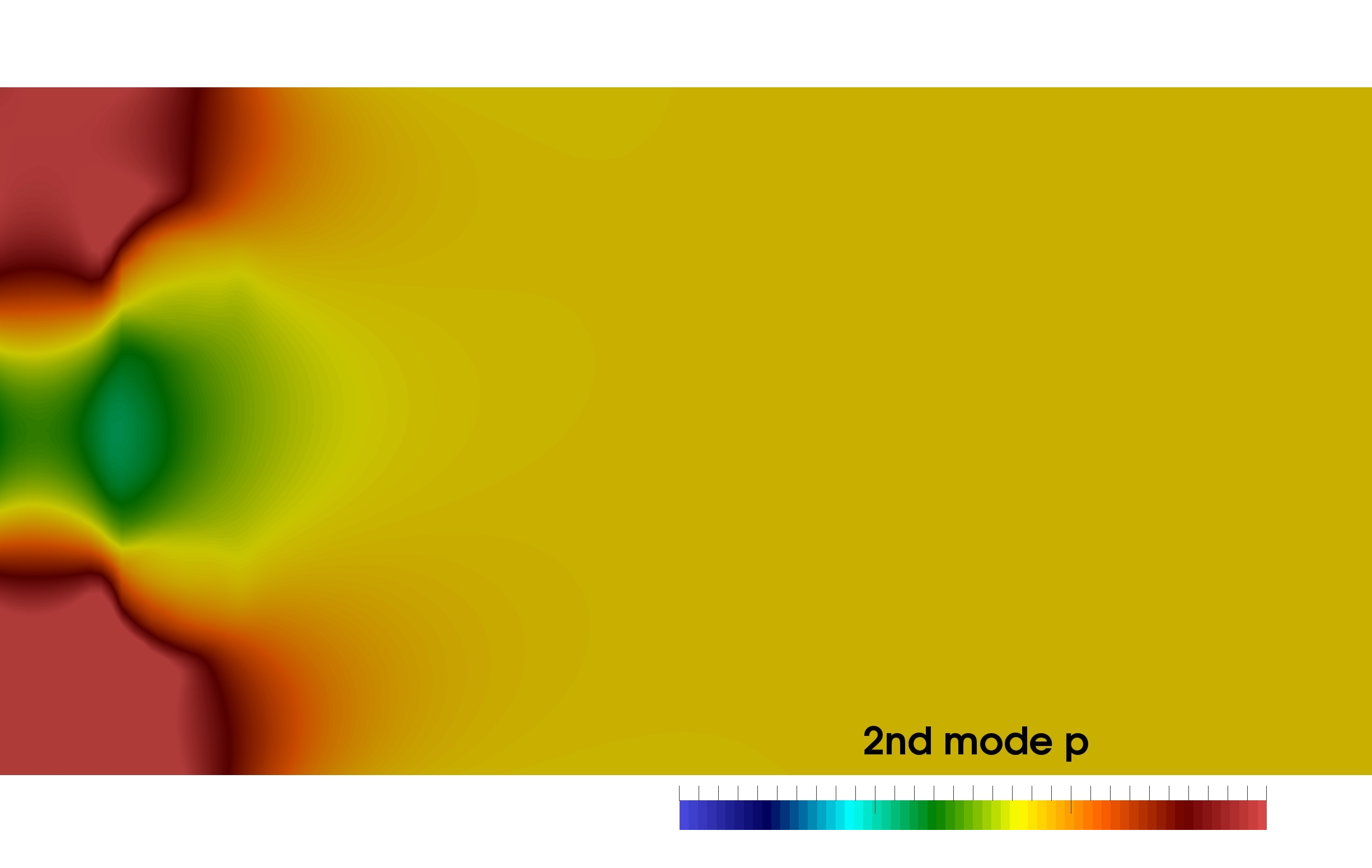}\\
\includegraphics[width=0.47\textwidth]{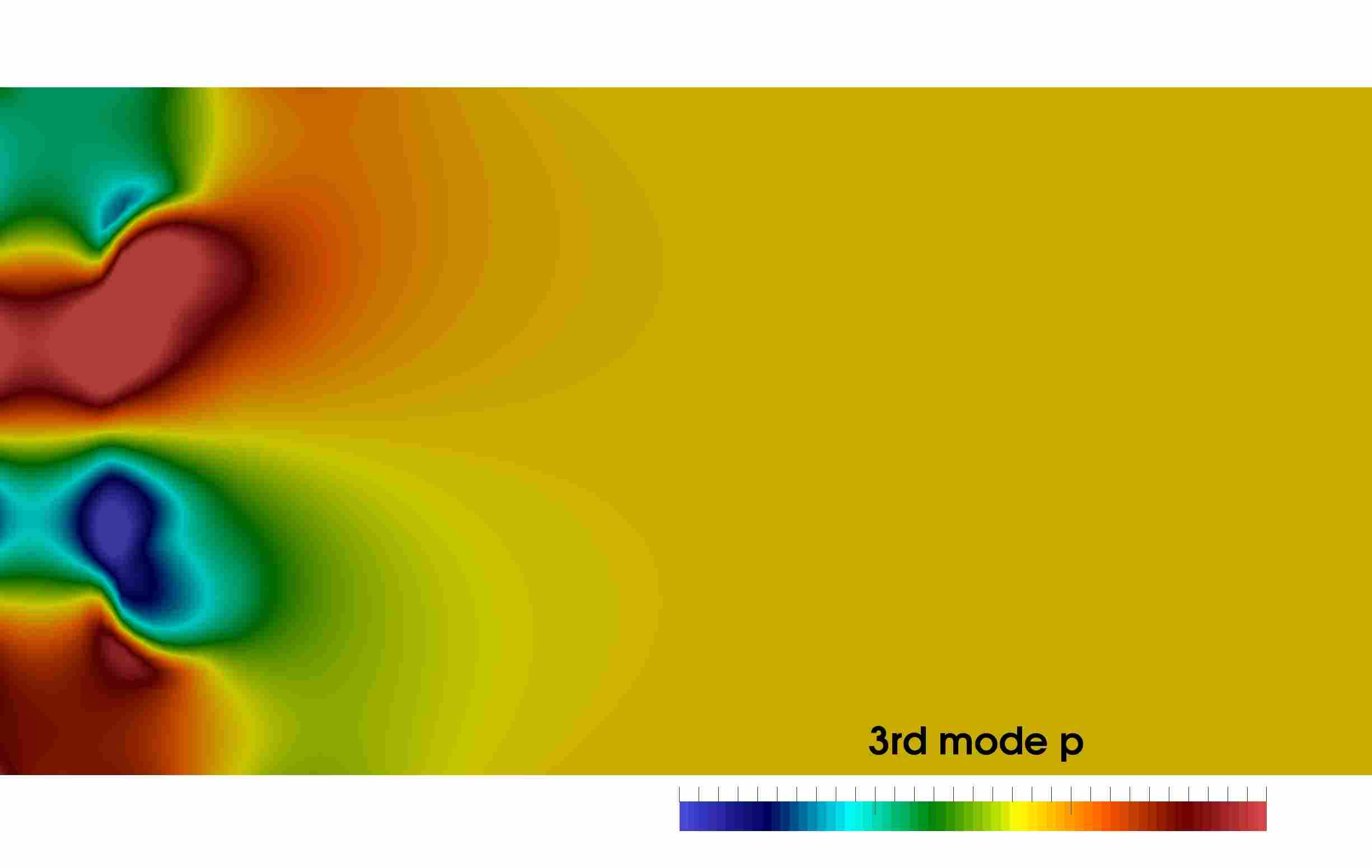}
\includegraphics[width=0.47\textwidth]{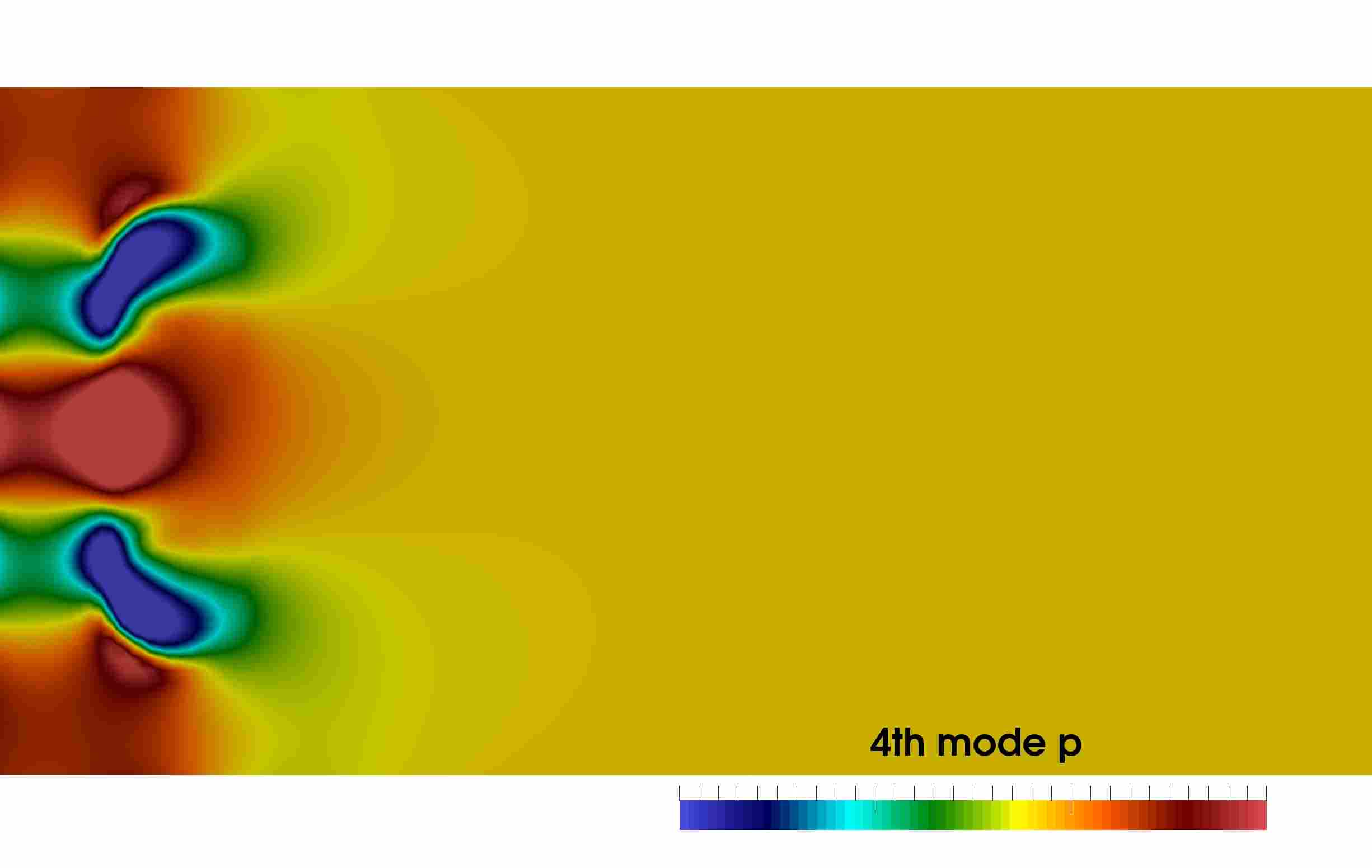}
\caption{\emph{Stokes flow test case}: first four pressure POD modes for the natural smooth extension without transportation.}
\label{Fig:Stokes_p_no_transport_modes}
\end{figure}
\begin{figure}
\centering
\includegraphics[width=0.47\textwidth]{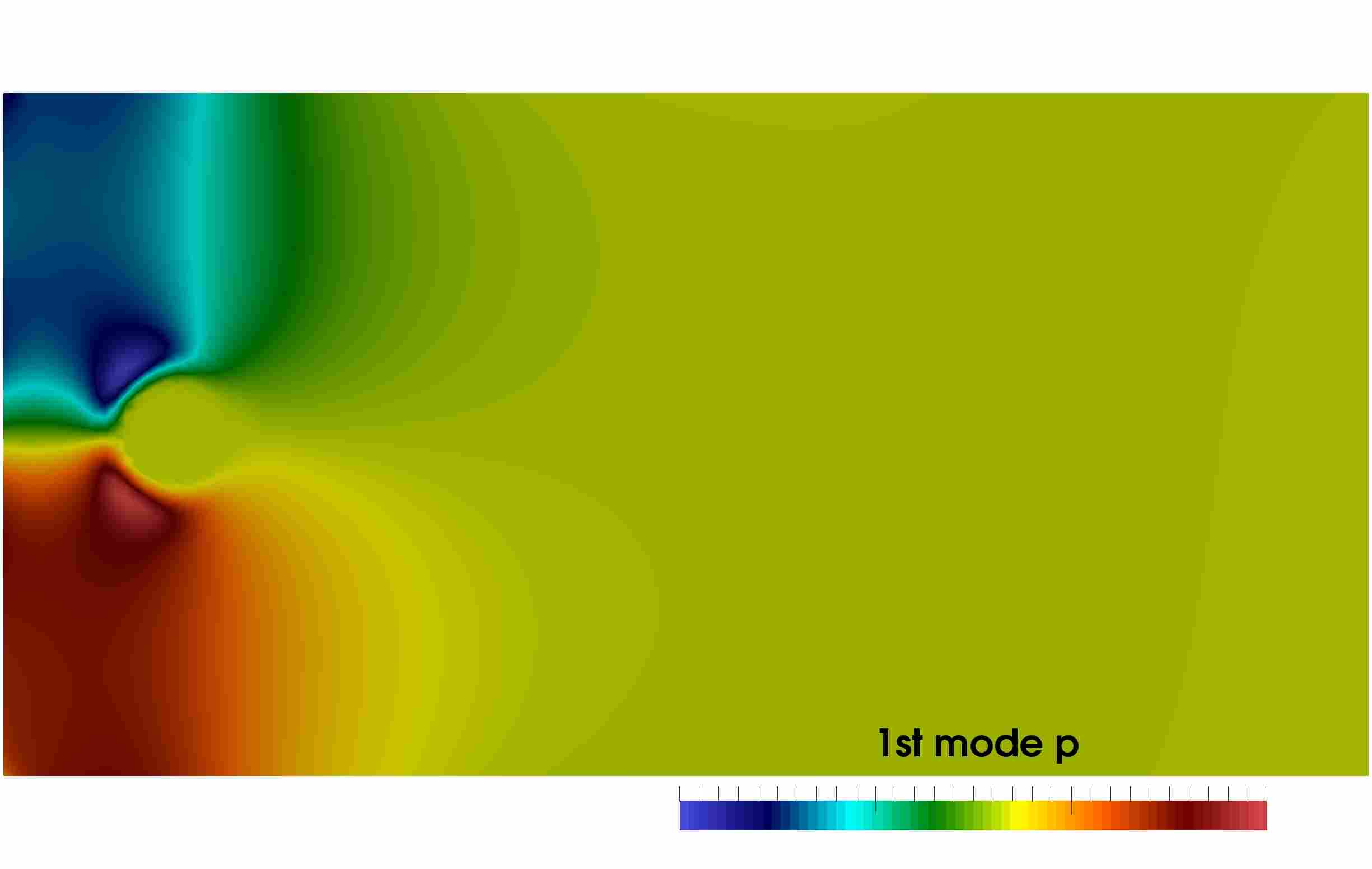}
\includegraphics[width=0.47\textwidth]{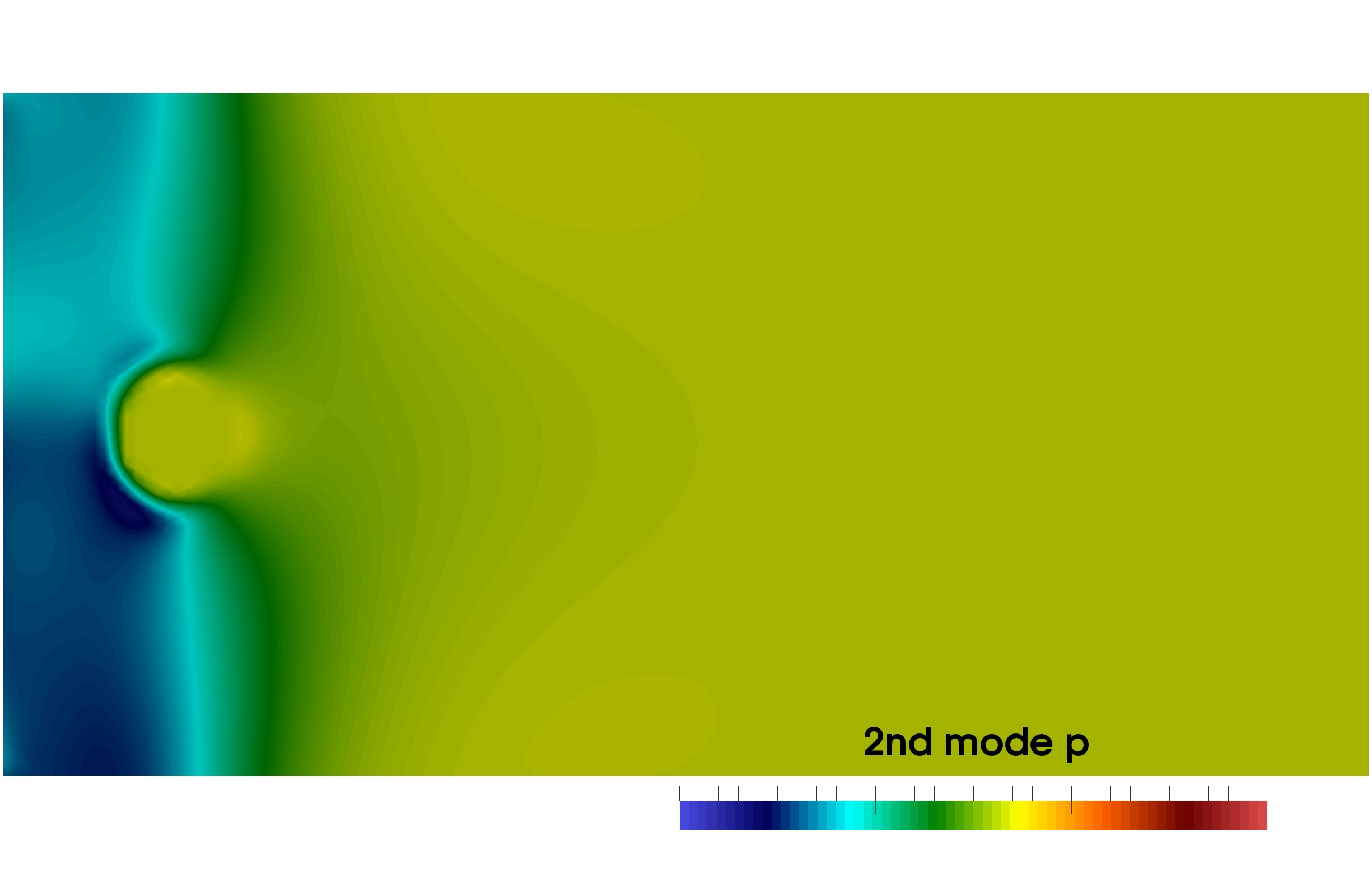}\\
\includegraphics[width=0.47\textwidth]{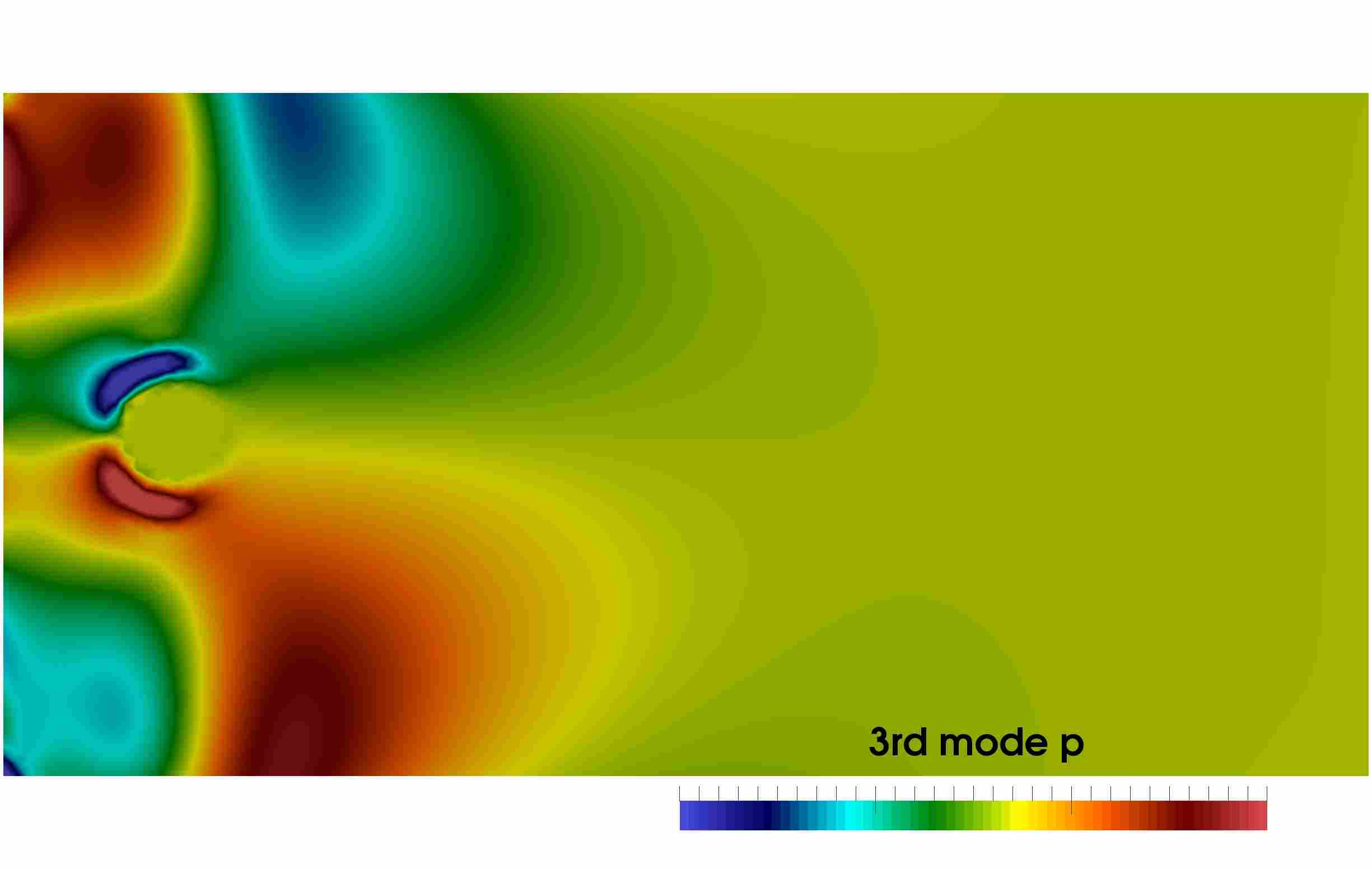}
\includegraphics[width=0.47\textwidth]{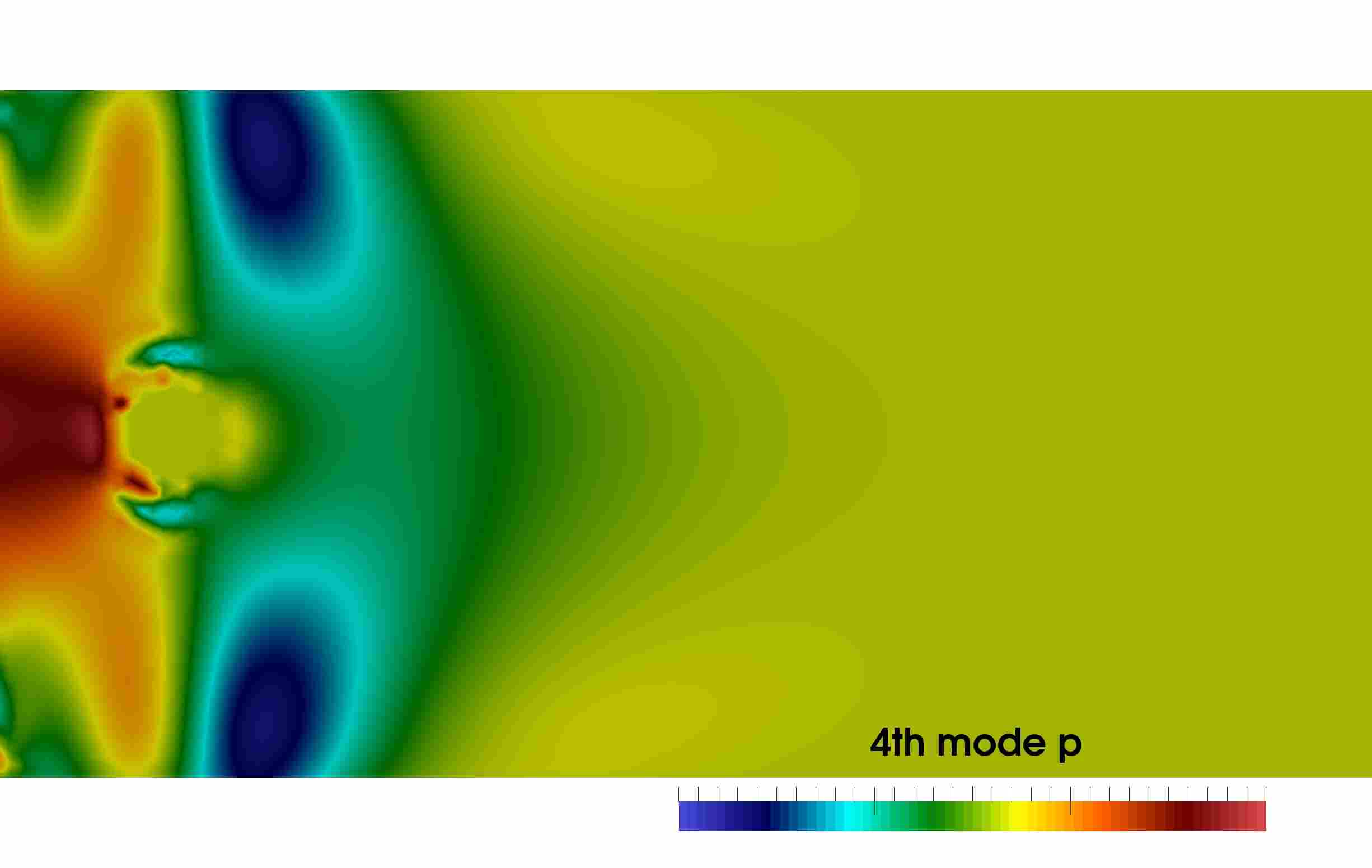}
\caption{\emph{Stokes flow test case}: first four pressure POD modes for the natural smooth extension with transportation.}
\label{Fig:Stokes_p_transport_modes}
\end{figure}
%

The results of the offline stage, run over a training set of 600 snapshots, are summarized in Figure \ref{fig:Stokes_eigs_Cavity}, where POD eigenvalues (normalized to the maximum eigenvalue) are plotted against the number of modes, for the velocity $u$, supremizer $s$ and pressure $p$.
As in the scalar Darcy flow experiments, results show that smaller reduced basis spaces can be obtained by combining extension and transportation procedures; indeed, Figure \ref{fig:Stokes_eigs_Cavity} shows a slower eigenvalue decay for the case without transportation. The difference among the two options are less marked here than in the scalar case due to the simpler parametrization, which only affects one geometrical quantity rather than four. Nonetheless, the qualitative difference between the resulting POD modes is still very appreciable comparing Figures \ref{Fig:Stokes_u_no_transport_modes}-\ref{Fig:Stokes_p_no_transport_modes} (velocity and pressure, respectively, without transportation) to Figures \ref{Fig:Stokes_u_transport_modes}-\ref{Fig:Stokes_p_transport_modes} (velocity and pressure, respectively, with transportation). In particular, close to the inlet, POD basis functions without transportation are characterized by peaks/sinks in the $y$ direction, their number increasing with the basis index. Such behavior is required for such basis to capture the moving circle. In contrast, no such phenomenon is present when applying transportation to the snapshots.

\begin{figure}
\centering
\includegraphics[width=0.7\textwidth]{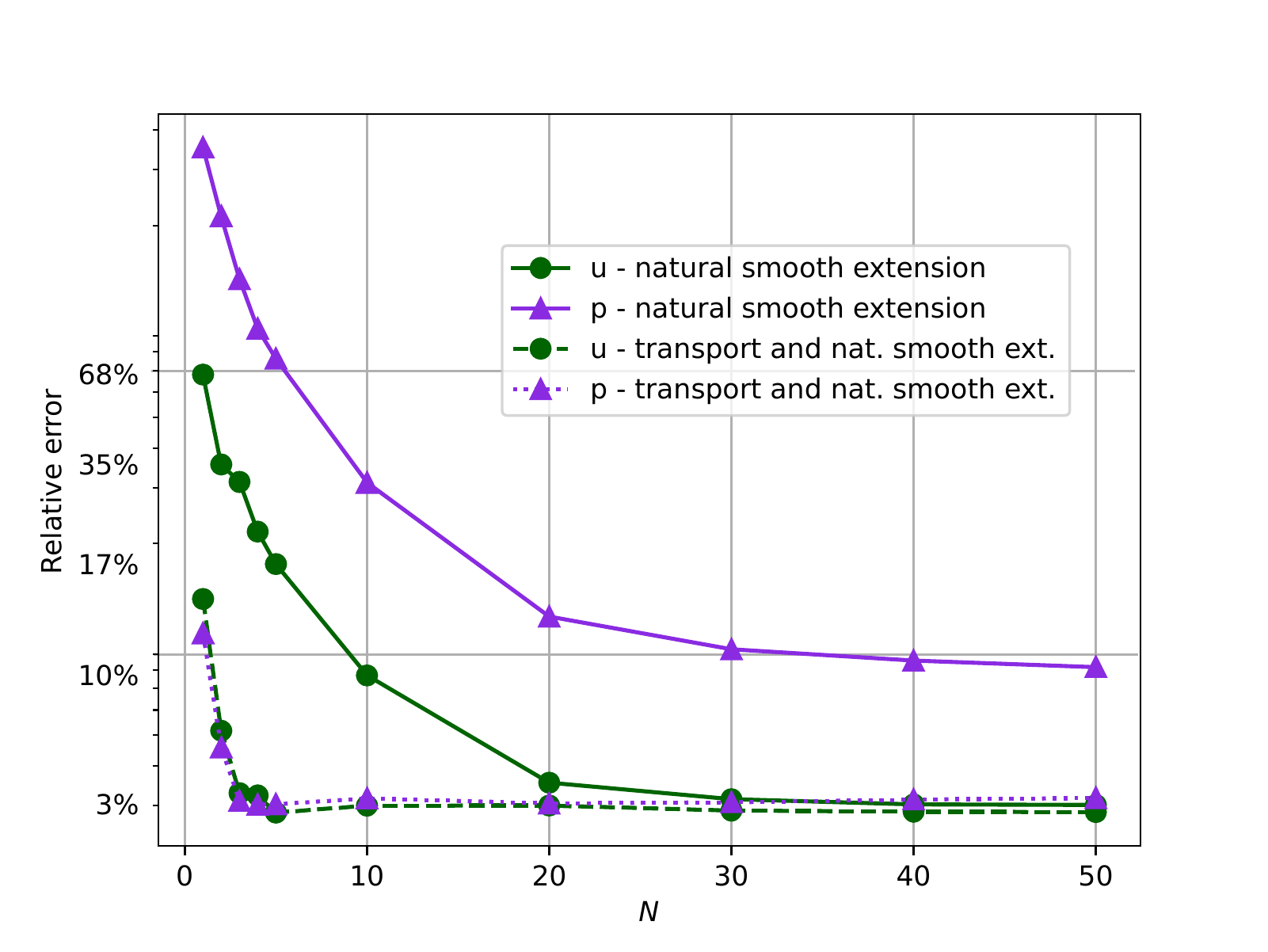}
\caption{\emph{Stokes flow test case}: error analysis between high fidelity approximations and reduced order.}
\label{fig:Stokes_R_errors_Cavity}
\end{figure}

\begin{figure}
\centering
\includegraphics[width=0.47\textwidth]{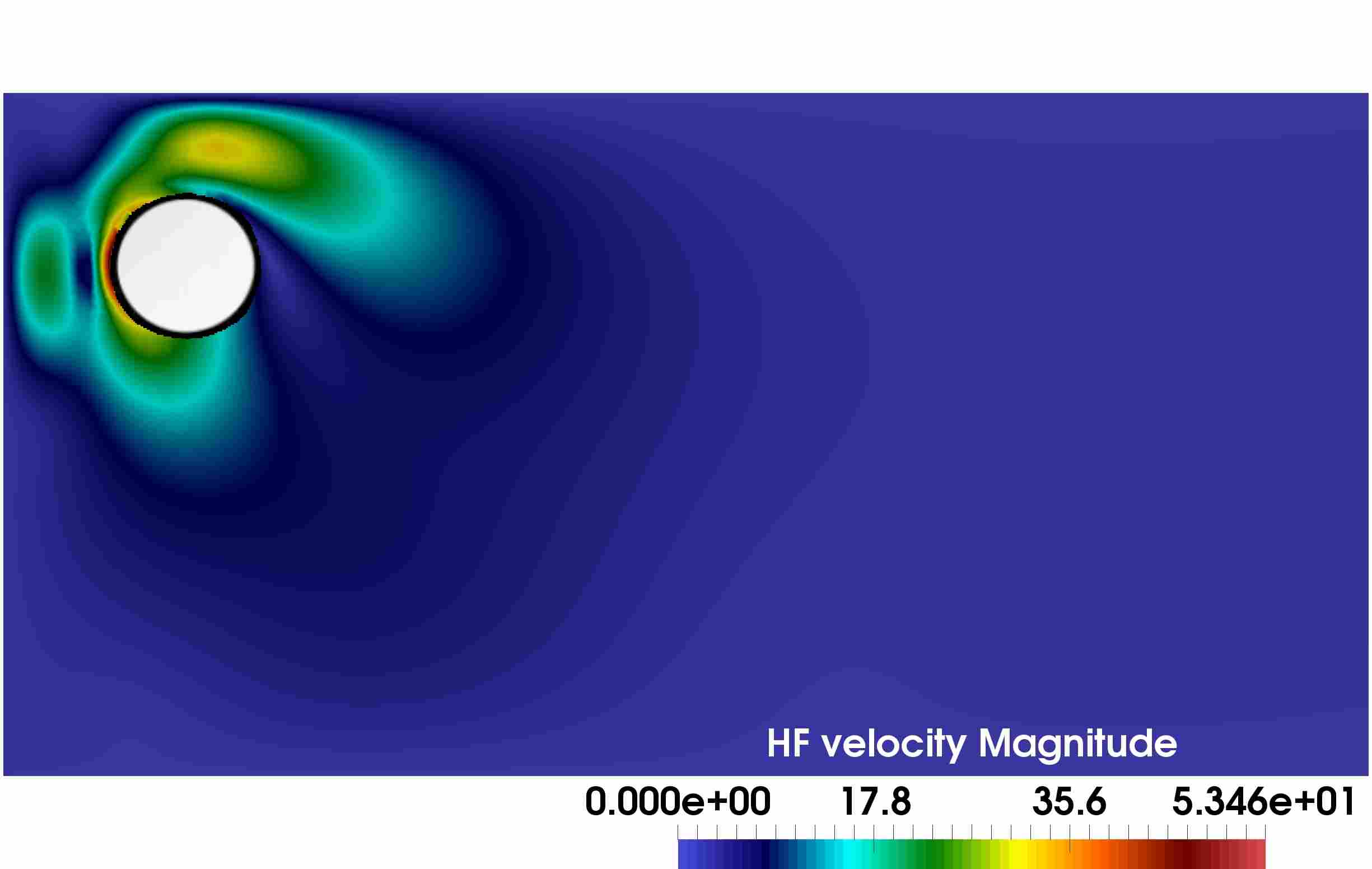}
\includegraphics[width=0.47\textwidth]{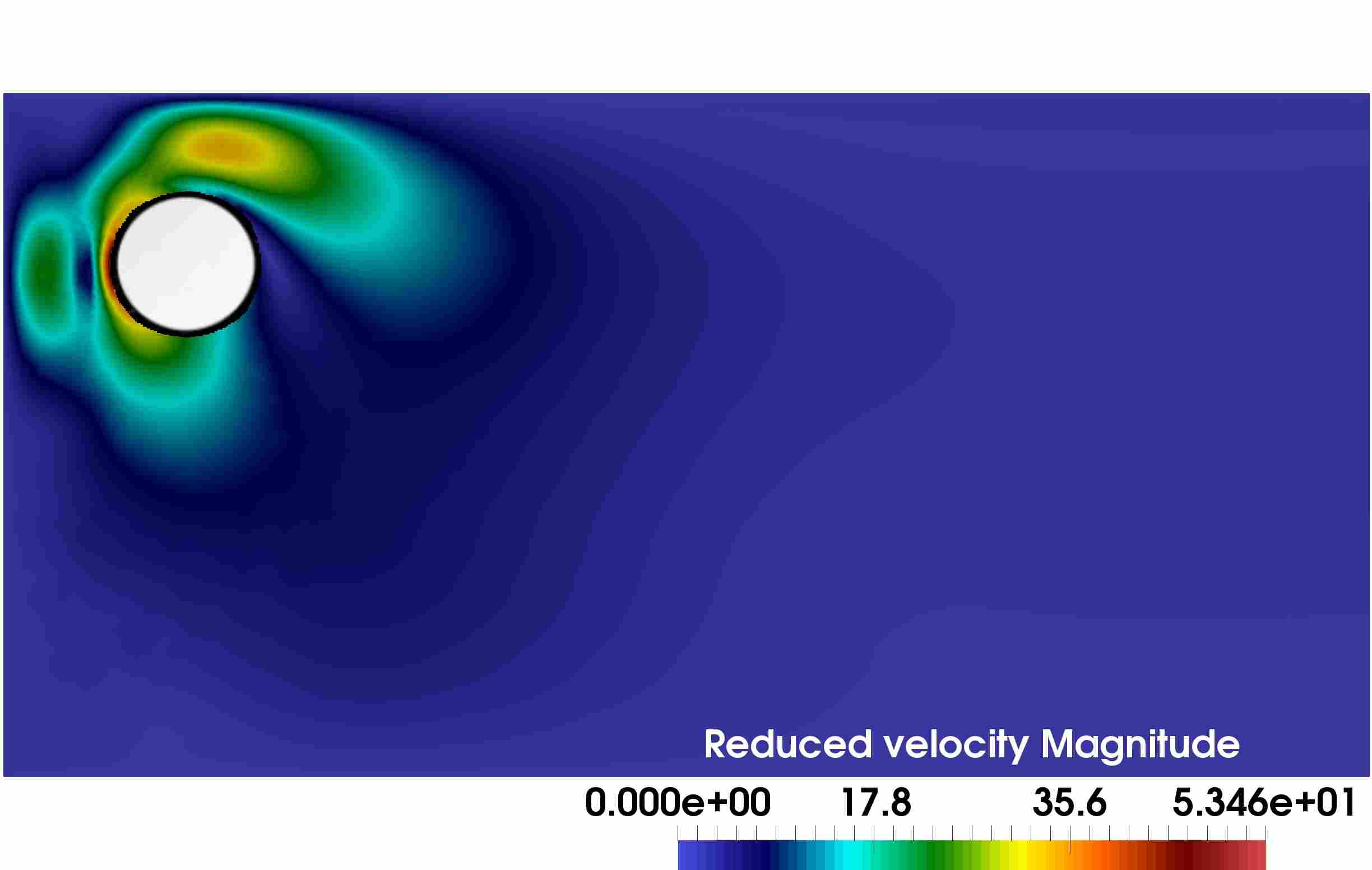}\\
\includegraphics[width=0.47\textwidth]{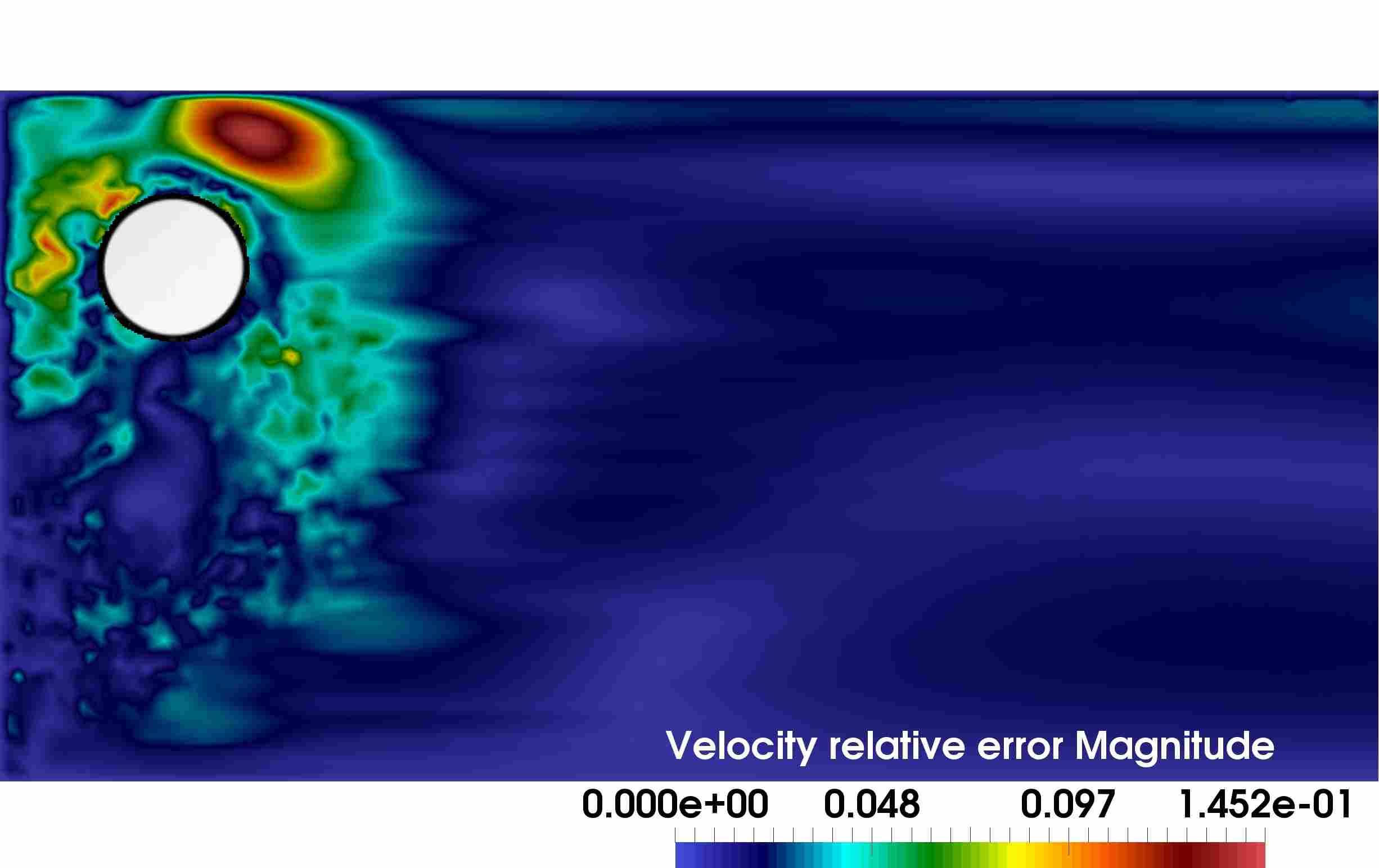}
\caption{\emph{Stokes flow test case}: the high fidelity velocity (top left), the reduced order velocity (top right) for the natural smooth extension with transportation, and the corresponding relative error (bottom) for the random parameter $\mu = 0.4998$.}
\label{Fig:Stokes_u_solution}
\end{figure}
\begin{figure}
\centering
\includegraphics[width=0.47\textwidth]{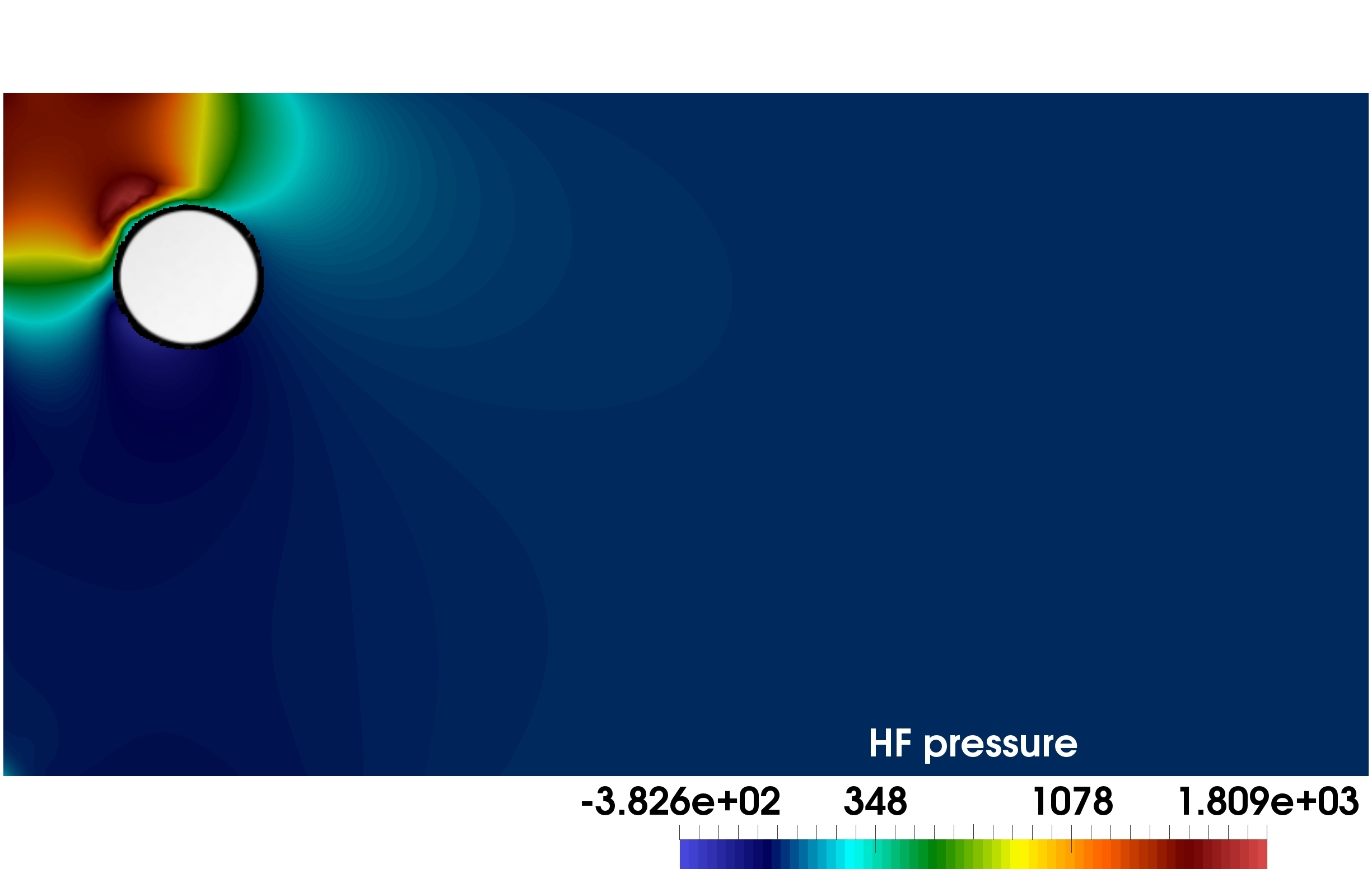}
\includegraphics[width=0.47\textwidth]{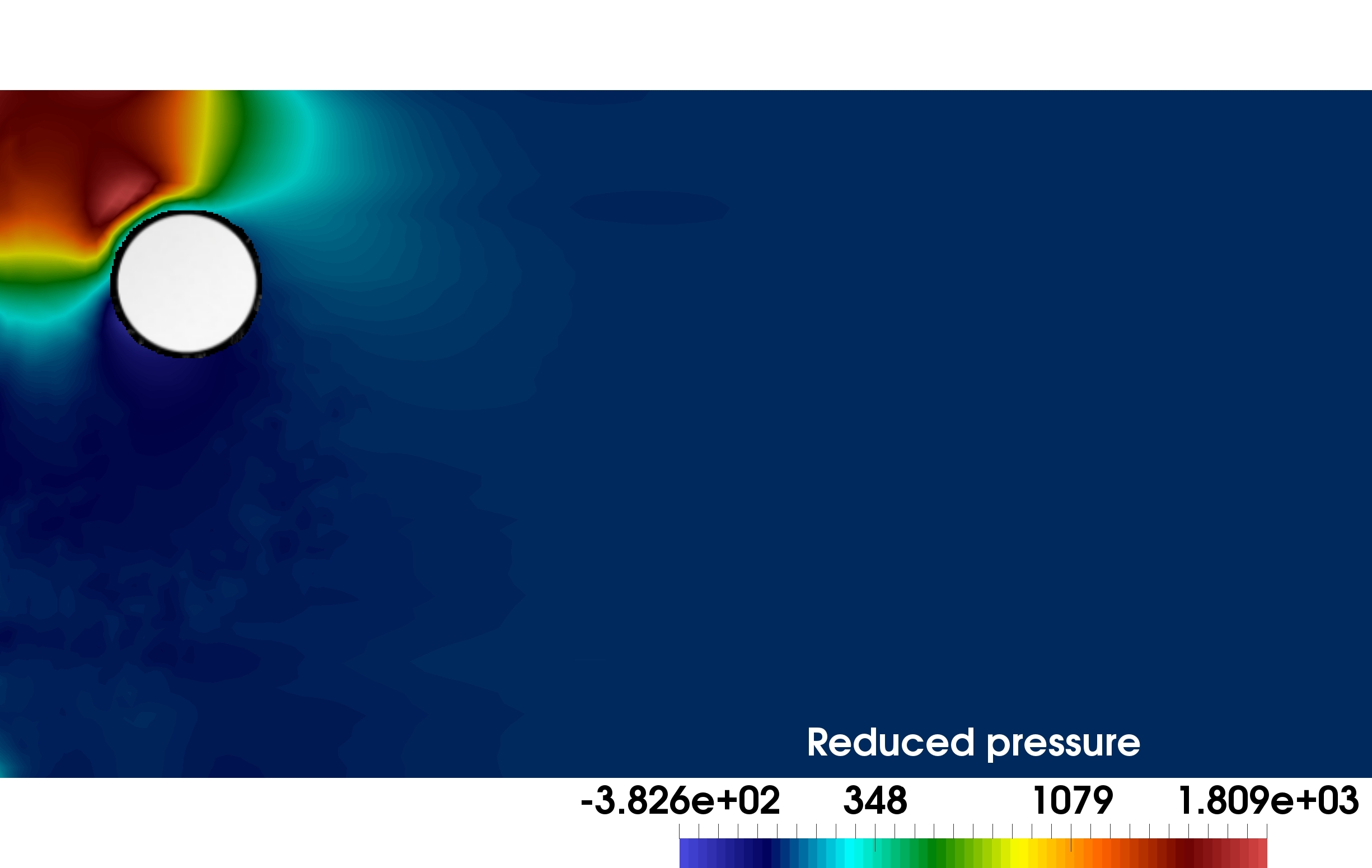}
\includegraphics[width=0.47\textwidth]{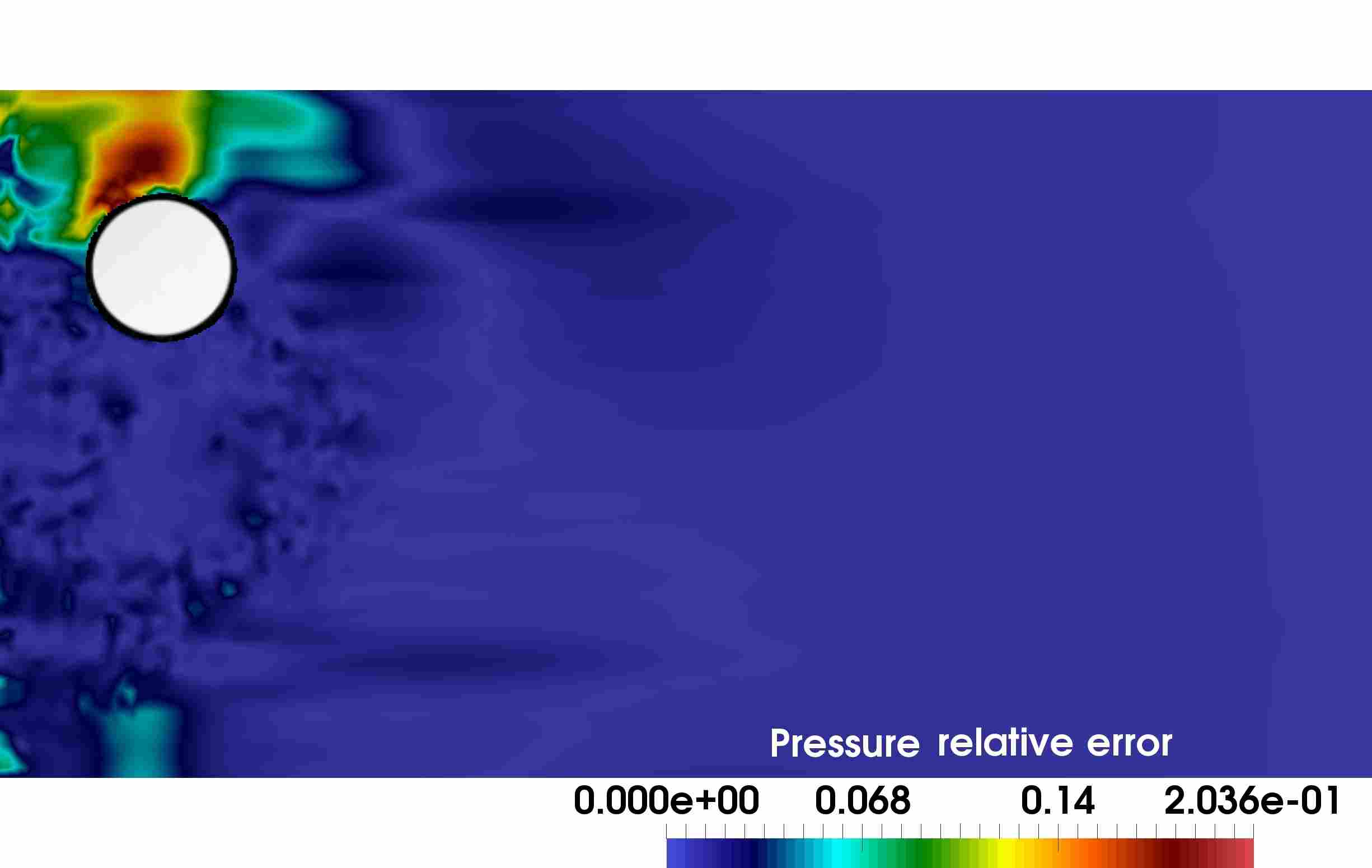}
\caption{\emph{Stokes flow test case}: the high fidelity pressure (top left), the reduced order pressure (top right) for the natural smooth extension with transportation, and the corresponding relative error (bottom) for the random parameter $\mu = 0.4998$.}
\label{Fig:Stokes_p_solution}
\end{figure}

Figure \ref{fig:Stokes_R_errors_Cavity} shows the $L^2$ relative errors results, for velocity and pressure, of an error analysis between the reduced order and high fidelity approximations over a testing set of $100$ parameter values. In particular, the average of the relative error over the testing set is plotted against the reduced basis size. In the case of natural smooth extension, the reduced pressure solution obtained is still inaccurate even for $N=50$ even with the addition of supremizers, being affected by relative errors of the order of $10^{-1}$. Reduced velocity solution is more accurate, resulting in errors of the order of $10^{-2}$ for $N\ge20$. In contrast, the combination to the transportation procedure results in relative errors that is of the order of $10^{-2}$ starting from $N=3$ for both velocity and pressure. Thus, again the role of snapshots transportation can be inferred from these results, being capable of improving the results of almost an order of magnitude compared to the natural smooth extension for a very low number of basis functions.
Unfortunately, both methods again reach a plateau after which no further improvement is shown. As in the elliptic case, we claim that this is due to weak imposition of boundary conditions, as errors tend to concentrate in a neighborhood of the embedded circle (especially for the pressure, which will thus affect the downstream velocity; see Figures \ref{Fig:Stokes_u_solution}-\ref{Fig:Stokes_p_solution} for a representative case); further work on this topic is forthcoming.

\section{Conclusions and future developments}\label{sec:conclusions}
In this work a POD-Galerkin ROM based on CutFEM high fidelity simulations was presented for linear PDE problems (elliptic and Stokes), characterized by a geometrical parametrization with (possibly) large variations. A CutFEM discretization naturally allows to use a level set description of the parametrized geometry. In our opinion, this results in a simpler and more versatile high fidelity method when compared to a FE formulation with pull back to a reference domain.

A separation between construction and evaluation phases is sought. While the evaluation of the ROM follows a standard Galerkin projection, a careful adaptation to embedded methods has been necessary for the construction stage, especially for what concerns the definition of snapshot on a common background mesh through a suitable combination of extension to the background mesh (in order to have all snapshots defined on a common mesh) and transportation on the background mesh itself (in order to enforce a rapid decrease of the Kolmogorov $n$-width). Indeed, thanks in particular to the transportation step, the developed ROM is able to reproduce the high fidelity solution in an accurate manner, with relative errors of the order of $10^{-4}$ for the elliptic case and $10^{-2}$ for the Stokes case.

Even though the proposed ROM follows a construction-evaluation paradigm, it is not offline-online separable in the usual sense \cite{HeRoSta16}, yet. Such separability, which is required for the efficient evaluation of the reduced order system, is a perspective to be thoroughly investigated in future. Suitable hyper reduction techniques, such as the empirical interpolation method \cite{BARRAULT2004667,stabile_geo_}, should be employed for this goal. Certification of the error is another topic of remarkable interest, in view of the application of greedy algorithms during the generation of the reduced basis space. Furthermore, a deep investigation on the role of Nitsche terms on the accuracy of the resulting ROM is needed, as numerical results show a plateau due to the weak imposition of Dirichlet boundary conditions.

As a further future development, we mention the extension of the proposed ROM to nonlinear problems in fluid dynamics, such as the Navier-Stokes equations. Moreover, a prototypical case for the application of embedded methods are fluid-structure interaction problems, which is another future perspective of this work. This will introduce additional complexities such as  additional equations to account for the structural dynamics, as well as fluid structure-interaction coupling. 

\section*{Acknowledgments}
We acknowledge the support by European Union Funding for Research and Innovation -- Horizon 2020 Program -- in the framework of European Research Council Executive Agency: Consolidator Grant H2020 ERC CoG 2015 AROMA-CFD project 681447 ``Advanced Reduced Order Methods with Applications in Computational Fluid Dynamics'' and  FSE project - European Social Fund - HEaD "Higher Education and Development" SISSA operazione 1, Regione Autonoma Friuli Venezia-Giulia. We also acknowledge INdAM-GNCS projects ``Tecniche di Riduzione di Modello per Applicazioni Mediche'', ``Metodi numerici avanzati combinati con tecniche di riduzione computazionale per PDEs parametrizzate e applicazioni'' and HFRI and GSRT under grant agreement No 1115. Numerical simulations were carried out using the {\emph{ngsxfem}} extension of {\emph{ngsolve}} \cite{ngsxfem,ngsolve} for the high fidelity simulations, and RBniCS \cite{BaRoSa15} for the reduced order ones. We acknowledge developers and contributors to both libraries. We thank in particular Prof. Christoph Lehrenfeld, Institute for Computational and Applied Mathematics, G\"{o}ttingen, for useful instructions regarding the {\emph{ngsxfem}} library.

\bibliographystyle{bst/elsarticle-num} 
\bibliography{bibfile_sissa}
\end{document}